\documentclass[12pt]{amsart} 
\title{\texttt{mcode.sty} Demo}
\usepackage{epsfig}
\usepackage{siunitx} 
\usepackage{color,graphics}
\usepackage{booktabs}
\usepackage{epstopdf}
\usepackage{soul}
\usepackage[all]{xy}
\usepackage{graphicx}
\usepackage{xcolor}
\usepackage{pagecolor}
\usepackage{amsmath, amssymb, amsthm}
\usepackage[vlined, ruled]{algorithm2e} 
\usepackage{fixmath, bm, amssymb, amsfonts, latexsym}
\usepackage{xcolor}
\usepackage{caption}
\usepackage{subcaption}
\usepackage{setspace}
\usepackage{ragged2e}
\usepackage{multirow}
\usepackage{hyperref}

\newtheorem{theorem}{Theorem}[section]

\newtheorem{lemma}[theorem]{Lemma}

\newtheorem{remark}[theorem]{Remark}

\newcommand{\baa}{\begin{eqnarray*}}
\newcommand{\eaa}{\end{eqnarray*}}
\newcommand{\ba}{\begin{equation}}
\newcommand{\ea}{\end{equation}}

\newcommand{\R}{\mathbb{R}}
\renewcommand{\epsilon}{\varepsilon}

\pagecolor{white}

\usepackage{float}

\begin{document}
\author{Sudipta Sahu$^\S$, Emanuele Macca$^\dagger$, Rathan Samala$^{\S*}$}
\title{Semi-Implicit Central scheme for Hyperbolic Systems of Balance Laws with Relaxed Source Term}
\thanks{
$^{\S}$ Department of Humanities and Sciences, Indian Institute of Petroleum and Energy-Visakhapatnam,
India-530003 ({sudipta.sahu@iipe.ac.in, rathans.math@iipe.ac.in}),
\newline
$^{\dagger}$  Dipartimento di Matematica ed Informatica, Università di Catania, Viale Andrea Doria 6, 95125
Catania, Italy (emanuele.macca@unict.it),
}
\maketitle
\begin{abstract}
 Quasi-linear hyperbolic systems with source terms introduce significant computational challenges due to the presence of a stiff source term. To address this, a finite volume Nessyahu-Tadmor (NT) central numerical scheme is explored and applied to benchmark models such as the Jin-Xin relaxation model, the shallow-water model, the Broadwell model, the Euler equations with heat transfer, and the Euler system with stiff friction to assess their effectiveness. The core part of this numerical scheme lies in developing a new implicit-explicit (IMEX) scheme, where the stiff source term is handled in a semi-implicit manner constructed by combining the midpoint rule in space, the trapezoidal rule in time with a backward semi-implicit Taylor expansion. The advantage of the proposed method lies in its  stability region and maintains robustness near stiffness and discontinuities, while asymptotically preserving second-order accuracy. The numerical validation further extends to two-dimensional configurations of the Jin-Xin relaxation model and a Jin-Xin-type relaxation system of 2D Euler equation. Theoretical analysis and numerical validation confirm the stability and accuracy of the method, highlighting its potential for efficiently solving the stiff hyperbolic systems of balance laws of 1D and 2D.
\end{abstract}

\bigskip
\noindent 2000 AMS(MOS) Classification:
41A10, 65M06


\pagestyle{myheadings} \thispagestyle{plain} \markboth{Sudipta Sahu,  Emanuele Macca, Rathan Samala}{ Semi-Implicit Central scheme for Hyperbolic Systems of Balance Laws with Relaxed Source Term}

%
\section{Introduction}
Hyperbolic systems of balance laws with stiff source terms appear in a wide spectrum of applications, ranging from fluid dynamics to geophysical modeling. Numerical methods for such systems must effectively reconcile the conflicting requirements of accuracy and efficiency, especially when confronted with stiff source terms that impose severe time-step restrictions under explicit discretizations. Traditional fully explicit schemes suffer from restrictive Courant-Friedrichs-Lewy (CFL) conditions \cite{Shu1997,Castro2020, Macca2024CATMOOD} while fully implicit approaches, although unconditionally stable, become computationally prohibitive due to the need to solve large and nonlinear systems at each time step \cite{Quinpi}. To bridge this gap, semi-implicit methods have been developed, which treat the non-stiff convective fluxes explicitly and the stiff source terms implicitly. This dual approach enhances computational efficiency while maintaining stability and accuracy; for relaxation settings leading to \emph{linearly implicit} stages we refer to~\cite{PareschiRusso}, while more general semi-implicit formulations are discussed in~\cite{Boscarino-Filbet,Macca2024,MaccaRussoBumi,Thomann2025}.

In recent years, central schemes have gained widespread acceptance for approximating hyperbolic conservation laws. Starting from the seminal work of Nessyahu and Tadmor (NT) \cite{NT}, many extensions have been proposed to incorporate source terms by means of non-splitting strategies \cite{russo,Liotta,pareschi2001central}. However, these approaches have commonly been formulated in the finite difference framework, and their adaptation to a finite volume setting, particularly with a view toward improving the integration of source terms, remains an area of active research.

The work presented here is motivated by the need for an efficient semi-implicit finite volume scheme that not only relaxes the strict CFL time-step restrictions but also attains at least second-order accuracy in both the convective and source term discretizations. To this end, we propose a novel method that employs a modified trapezoidal rule to discretize the source term in finite volume framework inspired from \cite{Macca} which is proposed in finite difference framework. By incorporating this quadrature strategy within the finite volume framework, our method maintains consistency and stability even in the presence of stiff relaxation.
Our proposed approach, based on the central schemes, may be regarded as an alternative to implicit-explicit (IMEX) Runge-Kutta schemes \cite{Boscarino-Filbet,ImexBosca}.  Note that the finite volume discretization that applies central differencing to approximate the convective fluxes and employs a modified trapezoidal rule for the semi-implicit treatment of the stiff source terms circumvents the need for Riemann solvers while still efficiently capturing the dynamics associated with stiff source contributions \cite{Toro2009}. This strategy offers a viable alternative to conventional IMEX RK methods by relaxing the CFL restrictions and accurately resolving slow, source-dominated processes, as encountered in models such as the Jin-Xin system \cite{jin1995runge,natalini1996convergence}, Broadwell model \cite{broadwell1964shock} and Euler equation with heat transfer \cite{Macca}, as well as the shallow water system and the Euler system with stiff friction \cite{SW_friction, Euler_friction}. Moreover, the scheme is numerically extended to multidimensional settings and applied to two-dimensional Jin-Xin and Euler Jin-Xin-type relaxation systems \cite{SW_friction}, demonstrating its robustness in higher-dimensional regimes.

Furthermore, an important key aspect of the proposed scheme is its asymptotic-preserving (AP) property. This feature guarantees that as the stiffness parameter tends to zero, the numerical scheme naturally transitions to a consistent discretization of the reduced (equilibrium) model without the need for modifying the time step or spatial discretization \cite{SB}. Moreover, we show that the proposed method is also \emph{asymptotic-accurate} (AA), in the sense that it preserves the designed second-order accuracy uniformly as $\varepsilon \to 0$, yielding a second-order consistent discretization of the equilibrium limit without order degradation. In other words, the scheme retains its stability and accuracy in the stiff limit, ensuring that the correct asymptotic behavior is captured even when the relaxation terms dominate. This property is particularly important in practical applications where explicit schemes would otherwise require extremely small time steps to maintain accuracy. The AP characteristic of our approach is established through both rigorous theoretical analysis and comprehensive numerical experiments, which confirm that the scheme converges to the appropriate equilibrium model as the stiffness parameter vanishes \cite{jin2010asymptotic,wanner1996solving}. In the AP literature, Abgrall and Torlo~\cite{AbgrallTorlo2020} propose high-order asymptotic-preserving IMEX schemes for kinetic models within a residual-distribution (RD) framework, where the temporal order is increased through a deferred-correction procedure (DeC). In this paper, we remain within a finite-volume central-scheme setting and obtain second-order AP/AA accuracy by integrating the stiff source through a modified trapezoidal rule with only two implicit stages, without relying on Riemann solvers or on RD/DeC machinery.

The remainder of the paper is organized as follows. In Section~2, we detail the governing equations of systems of balance laws with relaxation term with models like the Jin-Xin, Broadwell, and Euler equations with heat transfer, as well as the shallow water equations and the Euler equations with stiff friction in 1D, Jin-Xin and Euler Jin-Xin-type relaxation systems in 2D, which further demonstrate the applicability of the proposed approach.  In Section~3, we introduce the derivation of a semi-implicit finite volume scheme based on the midpoint rule in space and the trapezoidal rule with a backward semi-implicit Taylor expansion in time for the source term, while adapting the Nessyahu–Tadmor numerical approximation for the convective term for one and two dimension. Section~4 is dedicated to the theoretical analysis of the scheme, including consistency and stability proofs. Section~5 presents numerical experiments that validate the proposed method against five benchmark problems, including multidimensional test cases, demonstrating its superior performance and efficiency. Finally, in Section~6 some conclusions are drawn and future perspectives are discussed.

%
\section{Governing equations}
\label{sec:Governing_equations}

This paper deals with the design of second-order accurate methods for hyperbolic systems of balance laws with relaxation term. Nonetheless, let us consider the conservative form of the hyperbolic equation 
\begin{equation}\label{equ:scal_equ}
u_t(x,t) + f(u(x, t))_x = 0,
\end{equation}
with initial condition $u(x, 0) = u_0(x),$
where $u:\R\times[0,+\infty) \rightarrow\R$ is the unknown and $f:\R\rightarrow\R$ is the flux function.

In many applications it is necessary to include source terms; hence, we consider the scalar balance law
\begin{equation}
\label{equ:scalar_balance}
u_t(x,t) + f(u(x, t))_x = s(u(x,t)),
\end{equation}
with initial condition $u(x, 0) = u_0(x),$
where $s:\R\to\R$ is the source term that, in general, could represents  physical effects like relaxation or reaction.

In a more generic framework, the systems of balance law with relaxation term can be written as
\begin{equation}\label{equ:hyp_law}
U_t + F(U)_x = S(U, \varepsilon),
\end{equation}
with initial condition $U(x,0) = U_0(x)$, where $U : \mathbb{R} \times [0, +\infty) \rightarrow \mathbb{R}^d$ is the unknown vector field; $F : \mathbb{R}^d \rightarrow \mathbb{R}^d$ is the flux function, and $S : \mathbb{R}^d \times (0,1] \rightarrow \mathbb{R}^d$ is the source term that contains the stiff parameter $\varepsilon$. 
In the following subsections, we briefly introduce the models that will be considered in our article. 
Henceforth, in all models, we shall employ a simplified notation, whereby variables such as $u(x,t)$ will be denoted simply by $u$, in order to streamline the presentation.

\subsection{Jin-Xin Model}

The relaxation model developed by Jin-Xin is a mathematical framework used to describe fluid flows with relaxation effects \cite{jin1995runge}. This model incorporates the notion of a relaxation time, accounting for the finite time required for a fluid to reach its equilibrium state after an external disturbance. Relaxation terms capture the gradual adjustment of fluid properties, such as density, velocity, and temperature, towards their equilibrium values. The 1D Jin-Xin relaxed model can be written as
\begin{equation}\label{Xin_Jin_model}
\begin{cases}
u_t + v_x = 0,\\
v_t + u_x = -\dfrac{1}{\varepsilon}\left(v - au\right),
\end{cases}
\end{equation}
where $\varepsilon >0$ is the relaxation parameter. Equation \eqref{Xin_Jin_model} represents a particular relaxation system associated with a linear transport equation, which is widely referred to in the literature as the Jin-Xin relaxed model. In this work, it is employed as a prototypical example within the general relaxation framework introduced by Jin and Xin. System \eqref{Xin_Jin_model} may be written in the form \eqref{equ:hyp_law} by introducing the constitutive variable $U = [u,\, v]^T$, the flux $F(U) = [v,\, u]^T$, and the source term 
\[
S(U,\varepsilon) = -\dfrac{1}{\varepsilon} \begin{bmatrix} 0 \\ v- au \end{bmatrix}.
\]
When $\varepsilon \rightarrow 0$, we have $v = au$, and system \eqref{Xin_Jin_model} relaxes to the scalar conservation law
\begin{equation}\label{X_J_relax} 
u_t + a u_x = 0.
\end{equation}

\subsection{Shallow water Model}

We consider the shallow water system with a stiff relaxation source term, which can be written as
\begin{equation}\label{ShallowWater_model}
\begin{cases}
h_t + (hu)_x = 0, \\[2mm]
(hu)_t + \left(hu^2 + \tfrac{1}{2} h^2\right)_x =- \dfrac{1}{\varepsilon}\left(hu - \tfrac{1}{2} h^2\right),
\end{cases}
\end{equation}
where $h$ denotes the water height above the bottom topography and $hu$ represents the momentum. The system \eqref{ShallowWater_model} can be compactly expressed in the balance-law form \eqref{equ:hyp_law}, where the vector of conservative variables is
$$
U = \begin{bmatrix} h \\ hu \end{bmatrix},
$$
the flux function is
$$
F(U) = \begin{bmatrix} hu \\ hu^2 + \tfrac{1}{2} h^2 \end{bmatrix},
$$
and the relaxation source term is given by
$$
S(U,\varepsilon) = -\frac{1}{\varepsilon} \begin{bmatrix} 0 \\ hu - \tfrac{1}{2} h^2 \end{bmatrix}.
$$
In the asymptotic limit $\varepsilon \to 0$, the relaxation source enforces the algebraic constraint
$$
hu = \tfrac{1}{2} h^2,
$$
so that the system \eqref{ShallowWater_model} reduces to the scalar  Burgers' conservation law
\begin{equation}\label{ShallowWater_relax}
h_t + \left(\tfrac{1}{2} h^2\right)_x = 0.
\end{equation}
This represents the effective dynamics of the water height in the strong-relaxation regime.

\subsection{Broadwell Model}
The Broadwell model \cite{broadwell1964shock, PareschiRusso} is used to study the behavior of conserved quantities in the presence of relaxation processes. In this model, the conserved variables are represented by the vector 
\[
U = [\rho,\, m,\, z]^T,
\]
where $\rho$ denotes density, $u$ represents velocity, $m = \rho u$ denotes momentum and $z$ is an additional state variable. The flux function is given by 
\[
F(U) = \begin{bmatrix} m \\ z \\ m
\end{bmatrix},
\]
while the relaxation source term is defined as 
\[
S(U,\varepsilon) = \dfrac{1}{2\varepsilon} \begin{bmatrix} 0 \\ 0 \\ \rho^2 + m^2 - 2 \rho z \end{bmatrix}.
\]
Explicitly, the Broadwell model reads
\begin{equation}\label{Broadwell_model}
\begin{cases}
\rho_t + m_x = 0,\\[2mm]
m_t + z_x = 0,\\[2mm]
z_t + m_x = \dfrac{1}{2\varepsilon}\left(\rho^2 + m^2 - 2\rho z\right).
\end{cases}
\end{equation}
When $\varepsilon \rightarrow 0$, the relaxation yields $z = (\rho^2 + m^2)/{2 \rho}$, for $\rho \neq 0$ and system \eqref{Broadwell_model} relaxes to a 1D system of conservation laws:
\begin{equation}\label{B_W_relax}
\begin{cases}
\rho_t + m_x = 0,\\[2mm]
m_t + \left(\dfrac{\rho ^2 + m^2}{2 \rho}\right)_x = 0.
\end{cases}
\end{equation}

\subsection{Euler equations with heat transfer Model}

The Euler equations with heat transfer provide a fundamental description of compressible fluid dynamics while incorporating thermal energy exchange between the fluid and its surroundings. These equations extend the classical Euler equations by including an additional source term based on Fourier's law of heat conduction. In one dimension and neglecting gravitational acceleration, the Euler equations with heat transfer are given by
\begin{equation}\label{Euler_heat}
\begin{cases}
\rho_t + (\rho u)_x = 0,\\[2mm]
\left(\rho u\right)_t + \left(\rho u^2 + p\right)_x = 0,\\[2mm]
(\rho E)_t + \left(\rho u E + p u \right)_x = -K \rho \left(T - T_0\right),
\end{cases}
\end{equation}
where $\rho$ is the density, $u$ is the velocity, $E = e + \dfrac{u^2}{2}$ is the total energy per unit volume with $e$ denoting the internal energy, $p$ is the pressure, $T$ is the temperature, $K = 1/\varepsilon$ is the thermal conductivity, and $T_0$ is the temperature of the thermal bath. For a $\gamma$-law gas, the pressure is given by $p=(\gamma-1)\rho e$ and the internal energy is proportional to the absolute temperature, namely $e=\mathcal{R}T/(\gamma-1)=c_v T$, where $\mathcal{R}$ is the specific gas constant. These equations allow one to analyze complex fluid dynamics phenomena, including the propagation of shock and temperature waves, as well as the interactions between fluid flow and thermal effects.
\textcolor{red}{The} system \eqref{Euler_heat} can write explicitly, where $U = \left[\rho, \;\rho u, \;\rho E\right]^T$, the flux function
\[
F(U) = \begin{bmatrix}
\rho u\\
\left(1-\dfrac{\gamma -1}{2}\right)\dfrac{(\rho u)^2}{\rho}+(\gamma -1)\rho E \\
\dfrac{(\rho u)}{\rho}\left( \gamma \rho E - \left(\dfrac{\gamma -1}{2}\right)\dfrac{(\rho u)^2}{\rho}\right)
\end{bmatrix},
\]
and the source term with the relaxation parameter $\varepsilon$
\[
S(U, \varepsilon) = \begin{bmatrix}
0\\
0\\
\dfrac{1}{\varepsilon}\left(\rho T_0 -\dfrac{1}{c_v} \left(\rho E - \dfrac{(\rho u)^2}{2\rho}\right)\right)
\end{bmatrix}.
\]
When $\varepsilon \rightarrow 0$, the relaxation yields $\rho E = c_v \rho T_0 + (\rho u)^2 / {2 \rho}$, for $\rho \neq 0$ and system \eqref{Euler_heat} relaxes to a 1D system of conservation laws:
\begin{equation}\label{EulerHeat_relax}
\begin{cases}
\rho_t + (\rho u)_x = 0,\\[2mm]
(\rho u)_t + \left(\dfrac{(\rho u)^2}{2 \rho} + \left(\gamma -1\right)c_v \rho T_0\right)_x = 0.
\end{cases}
\end{equation}




\subsection{Euler system with stiff friction}
We consider the one-dimensional compressible Euler equations with a linear friction term. The primary unknowns are the mass density $\rho$, the velocity $u$, the specific internal energy $e$, and the specific total energy
$$
E = e + \tfrac{1}{2} u^2,
$$
so that the total energy per unit volume is $\rho E = \rho e + \tfrac{1}{2}\rho u^2$. The governing equations in conservative form read
\begin{eqnarray} \label{eulerGas:stiff_Friction}
\begin{cases}
    \begin{aligned}
        &\rho_t + (\rho u)_x = 0, \label{eq:mass_rhoE}\\[0.3em]
       &(\rho u)_t +  (\rho u^2 + p)_x = - \alpha \rho u, \label{eq:mom_rhoE}\\[0.3em]
       &(\rho E)_t +  \big( (\rho E + p) u \big)_x = - \alpha \rho u^2, \label{eq:energy_rhoE}
    \end{aligned}       
\end{cases}
\end{eqnarray}
where $\alpha>0$ is the friction coefficient. Closure is provided by the perfect gas law in conservative variables:
$$
p = (\gamma - 1)\,\rho e = (\gamma - 1)\,\rho\left(E - \tfrac{1}{2}u^2\right),
$$
with $\gamma>1$ the ratio of specific heats.

To study the strong-friction regime, we introduce the small parameter $\varepsilon = 1/\alpha$, so that large friction corresponds to $\varepsilon \ll 1$. The momentum equation \eqref{eq:mom_rhoE} can be written as
$$
(\rho u)_t + (\rho u^2 + p)_x = - \frac{1}{\varepsilon}\, \rho u.
$$
Here, the right-hand side is of order $1/\varepsilon$, much larger than the terms on the left-hand side, which are $O(1)$. 
As $\varepsilon \to 0$, the first two terms become negligible, and we obtain
$$
\rho u = O(\varepsilon).
$$
\noindent \textbf{Chapman-Enskog expansion and well-preparedness.}
Let $\varepsilon = 1/\alpha$ and denote $m:=\rho u$. We assume that the solution admits a Chapman-Enskog expansion in $\varepsilon$ of the form
\[
m = m^{(0)} + \varepsilon m^{(1)} + O(\varepsilon^2),
\qquad \rho,\; p = O(1)\ \text{as}\ \varepsilon\to 0.
\]
Substituting into the momentum equation
\[
m_t + \Big(\frac{m^2}{\rho} + p\Big)_x = -\frac{1}{\varepsilon} m
\]
and collecting terms of equal order yields, at leading order $O(1/\varepsilon)$, $m^{(0)}=0$, hence $m=O(\varepsilon)$ in the strong-friction regime. At the next order $O(1)$ we obtain
\[
p_x = -\,m^{(1)},
\qquad\text{so that}\qquad
m = -\varepsilon p_x + O(\varepsilon^2),
\quad
u = -\frac{\varepsilon}{\rho}p_x + O(\varepsilon^2)= -\frac{1}{\alpha\rho}p_x + O(\varepsilon^2).
\]
Therefore $m_t=O(\varepsilon)$ and $\big(m^2/\rho\big)_x=O(\varepsilon^2)$, whereas $(1/\varepsilon)m=O(1)$, which justifies neglecting the inertial terms $(\rho u)_t$ and $(\rho u^2)_x$ in the leading-order balance as $\varepsilon\to 0$.
The above expansion is consistent provided $m(x,0)=O(\varepsilon)$ (\emph{well-prepared} data); for \emph{ill-prepared} data $m(x,0)=O(1)$, a rapid initial relaxation layer on the time scale $t=O(\varepsilon)$ drives the solution toward the equilibrium relation $m\approx -\varepsilon p_x$, after which the limiting diffusive dynamics governs the evolution. 

This equilibrium closure can be interpreted as a Darcy-type law, expressing the balance between the pressure gradient and the friction force; inserting it into the continuity equation yields the corresponding diffusion-type limit in the strong-friction regime.

Substituting the above relation into the continuity equation yields
$$
\rho_t + (\rho u)_x \approx \rho_t - \left(\frac{1}{\alpha}\,p_x\right)_x.
$$
Since $p$ depends on $\rho$ and $\rho E$ through the equation of state, the strong-friction limit leads to a coupled parabolic-type evolution for $\rho$ and $\rho E$. In this limit the momentum equation becomes redundant and the evolution of mass and total energy is governed by diffusive relations driven by pressure gradients.

The conservative form \eqref{eulerGas:stiff_Friction} with $\rho E$ emphasizes that the energy equation describes the evolution of total energy per unit volume. In the strong-friction regime ($\alpha\gg1$) the kinetic part of $\rho E$ is rapidly dissipated by the friction source, and the remaining dynamics reduce to a slower diffusive evolution for the thermodynamic fields. This observation motivates asymptotic-preserving numerical schemes that treat stiff friction consistently while conserving the appropriate conservative variables.

\noindent\textbf{Isentropic Euler system with stiff friction.}
When the flow is assumed to be completely isentropic, the thermodynamic relation reduces to the simple equation of state
$$
p = k \rho^{\gamma}.
$$
where the pressure $p$ is expressed solely in terms of the density $\rho$. Unlike the full Euler system, in this setting the pressure does not depend on the total energy $\rho E$. As a result, the energy balance equation in \eqref{eulerGas:stiff_Friction} becomes redundant and can be removed from the system.
With this simplification, the governing equations reduce to the following system:
\begin{equation}
\label{eulerGas_isentropic:stiff_friction}
\begin{cases}
\rho_t + (\rho u)_x = 0, \\
(\rho u)_t + \big(\rho u^2 + k \rho^\gamma\big)_x = -\alpha \rho u,
\end{cases}
\end{equation}
which is commonly referred to as the isentropic Euler system with stiff friction.

\subsection{Two-Dimensional Governing Equations} 
We consider the two-dimensional hyperbolic conservation law
\begin{equation}\label{equ:hyp_law2d}
u_t + \partial_x f_1(u) + \partial_y f_2(u) = 0,
\end{equation}
for $(x,y)\in \mathbb{R}^2$ and $t>0$, where $u:\mathbb{R}^2\times[0,\infty)\to\mathbb{R}$ is the unknown solution and $f_1,f_2:\mathbb{R}\to\mathbb{R}$ denote the flux functions in the $x$- and $y$-directions. The equation is supplemented with the initial condition
\begin{equation}
u(x,y,0)=u_0(x,y).
\end{equation}
To incorporate relaxation effects, we consider the two-dimensional hyperbolic relaxation equation
\begin{equation}\label{eq:2d_model}
u_t + \partial_x f_1(u) + \partial_y f_2(u) = \frac{1}{\varepsilon} g(u),
\end{equation}
where $g(u)$ denotes the relaxation source term and $\varepsilon>0$ is the relaxation parameter.
\subsubsection{Two-Dimensional Jin-Xin Model}
The Jin-Xin relaxation system \cite{SW_friction} in two spatial dimensions is given by
\begin{equation}\label{Xin_Jin_2D_model}
\begin{cases}
u_t + v_x + w_y = 0,\\[2mm]
v_t + u_x = -\dfrac{1}{\varepsilon}(v - a u),\\[2mm]
w_t + u_y = -\dfrac{1}{\varepsilon}(w - b u),
\end{cases}
\end{equation}
where $u(x,y,t)$ is the primary unknown, $v(x,y,t)$ and $w(x,y,t)$ are auxiliary variables, $\varepsilon > 0$ is the relaxation parameter, and $a, b$ are given constants. In the relaxation limit ($\varepsilon \to 0$), the auxiliary variables satisfy $v = a u$ and $w = b u$, and the system reduces to the scalar conservation laws
\begin{equation}
u_t + a u_x + b u_y = 0.
\end{equation}
This reduced equation describes the asymptotic behavior of the system as the relaxation parameter becomes small.
\subsubsection{Two-Dimensional Euler Jin-Xin-Type Relaxation Model}
The two-dimensional Euler Jin-Xin type relaxation system \cite{SW_friction} is given by
\begin{equation}\label{Euler_Xin_Jin2d_model}
\begin{cases}
U_t + V_x + W_y = 0,\\[2mm]
V_t + A U_x = -\dfrac{1}{\varepsilon}\left(V - F_1(U)\right),\\[2mm]
W_t + B U_y = -\dfrac{1}{\varepsilon}\left(W - F_2(U)\right),
\end{cases}
\end{equation}
where $\varepsilon>0$ denotes the relaxation parameter. Here $V$ and $W$ are auxiliary variables introduced to approximate the fluxes in the $x$- and $y$-directions, respectively. The matrices $A$ and $B$ are diagonal matrices with eigenvalues $a_i$ and $b_i$, respectively, for $1\le i\le4$. The conserved variable vector is
\[
U=
\begin{bmatrix}
\rho\\
\rho u\\
\rho v\\
E
\end{bmatrix},
\]
where $\rho$ denotes the density, $u$ and $v$ are the velocity components in the $x$- and $y$-directions, respectively, and $E$ represents the total energy. The Euler flux functions are given by
\[
F_1(U)=
\begin{bmatrix}
\rho u\\
\rho u^2+p\\
\rho uv\\
(E+p)u
\end{bmatrix},
\qquad
F_2(U)=
\begin{bmatrix}
\rho v\\
\rho uv\\
\rho v^2+p\\
(E+p)v
\end{bmatrix}.
\]
The pressure $p$ is determined from the ideal gas equation of state
\[
p=(\gamma-1)\left(E-\frac12\rho(u^2+v^2)\right),
\]
where $\gamma$ is the ratio of specific heats. In the present work, we take $\gamma=1.4$. The matrices $A$ and $B$ are chosen such that the subcharacteristic condition is satisfied,
\[
a_i\ge \max|\Lambda_x(U)|,\qquad
b_i\ge \max|\Lambda_y(U)|,
\]
where $\Lambda_x(U)$ and $\Lambda_y(U)$ denote the eigenvalues of the Jacobians of the Euler fluxes in the $x$- and $y$-directions, respectively. In the relaxation limit $\varepsilon\to0$, the relaxation terms vanish and the auxiliary variables approach their equilibrium values $V=F_1(U)$ and $W=F_2(U)$. Consequently, system \eqref{Euler_Xin_Jin2d_model} reduces to the two-dimensional Euler system of conservation laws
\begin{equation}
U_t + F_1(U)_x + F_2(U)_y = 0.
\end{equation}
In the remainder of the paper, we focus on the development and analysis of semi-implicit schemes for the models introduced above, aiming to capture both the convective behavior and the stiff relaxation or heat transfer or friction effects present in the governing equations. 

Throughout the paper, $\varepsilon>0$ denotes the relaxation time scale; the proposed scheme is formulated for arbitrary $\varepsilon$ and does not require time-step restrictions.
\section{Methodology}
Before introducing the numerical schemes developed in this paper, we define the time domain as $\mathcal{T} = [0, T],\, T > 0,$
and partition it into subintervals $[t^n, t^{n+1}]$ for $n\in\mathbb{N}$.  Time step $\Delta t = t^{n+1} - t^n.$
Next, consider the one‐dimensional computational domain $\Omega$
which we subdivide into $N$ uniform cells
\[
\omega_i = \left[x_{i-1/2},\,x_{i+1/2}\right],
\qquad i = 1,\dots,N.
\]
 Each cell has width $\Delta x = x_{i+1/2} - x_{i-1/2},$
and its center is located at $x_i = \dfrac{1}{2}\left(x_{i-1/2} + x_{i+1/2}\right).$  Here, $x_i=i \Delta x$ are the grid points and $x_{i+1/2}=x_i+\Delta x /2$ are the cell interface values.
To reconstruct a high‐resolution, non‐oscillatory solution, we employ the MinMod slope limiter.  Define the forward and backward differences
\[
\delta^+_i(t^n) = u_{i+1}(t^n) - u_i(t^n),
\qquad
\delta^-_i(t^n) = u_i(t^n) - u_{i-1}(t^n),
\]
and set
\begin{equation}\label{eq:minmod}
 \Delta x\dfrac{\partial u}{\partial x}(x,t^n)\bigg|_{x=x_i}\approx u'_i\;=\; \operatorname{MM}(\delta^+_i(t^n),\;\delta^-_i(t^n)),
\end{equation}
where
\[
\operatorname{MM}(x,y) = 
\begin{cases}
\operatorname{sgn}(x) \cdot\,\min\!\left(|x|,|y|\right),
  & \text{if }\operatorname{sgn}(x)=\operatorname{sgn}(y),\\
0, & \text{otherwise}.
\end{cases}
\]
Within each cell $\omega_i$ we then build the piecewise‐linear interpolant
\begin{equation}\label{eq:linint}
 L_i(x,t^n) \;=\;
u_i(t^n)\;+\;(x - x_i)\,\dfrac{u'_i}{\Delta x},
\quad x_{i-1/2} \leq x \leq x_{i+1/2}.   
\end{equation}
\begin{remark}
In the finite-volume setting, the numerical solution is represented by cell averages over a prescribed spatial partition. At time $t^n$, these averages are defined on the non-staggered cells $\omega_i=[x_{i-\frac12},x_{i+\frac12}]$. In order to improve the accuracy of the approximation, a piecewise linear reconstruction is constructed within each cell based on these averages. In the context of central schemes, however, advancing the solution on the same set of control volumes may lead to insufficient numerical dissipation. The Nessyahu-Tadmor approach therefore considers an evolution on staggered control volumes $[x_i,x_{i+1}]$, which naturally incorporates additional averaging and enhances the stability of the method without relying on upwind information or Riemann solvers. Since the solution obtained after one update is defined on a staggered grid, a subsequent staggered evolution is performed to recover cell averages on the original non-staggered grid. This alternating evolution between non-staggered and staggered control volumes constitutes a fundamental feature of the NT central scheme.
\end{remark}

\subsection{The Nessyahu-Tadmor predictor–corrector scheme}
Let us now consider the one-dimensional relaxation scalar model
\begin{equation}\label{eq:mbalance}
u_t + f(u)_x = s(u, \varepsilon),
\end{equation}
where $s(u, \varepsilon): \mathbb{R} \times (0, 1] \to \mathbb{R}$ may exhibit stiffness. Let $\omega_i=[x_{i-1/2},x_{i+1/2}]$ denote the $i$-th computational cell. 
The cell average of the solution at time $t^n$ is defined by
\[
\bar{u}_i^n=\frac{1}{\Delta x}\int_{\omega_i} u(x,t^n)\,dx.
\]
For simplicity of notation, we drop the bar in the following and denote the cell averages by $u_i^n$. Defining $$u_{i+1/2}^{\,n+1}=\dfrac{1}{\Delta x}\int_{x_i}^{x_{i+1}} u(x, t^{n+1})dx,$$ we integrate \eqref{eq:mbalance} over the spacetime cell $[x_i,x_{i+1}]\times[t^n,t^{n+1}]$ gives
\begin{equation*}
\begin{aligned}
u_{i+1/2}^{\,n+1}=& \dfrac{1}{\Delta x}\int_{x_i}^{x_{i+1}} u(x, t^{n})dx+ \dfrac{1}{\Delta x}\int_{t^n}^{t^{n+1}}\left[f(u(x_i,t)) - f(u(x_{i+1},t))\right]\;dt\\
& + \dfrac{1}{\Delta x}\int_{t^n}^{t^{n+1}}\!\int_{x_i}^{x_{i+1}}s\left(u(x,t), \varepsilon\right)\;dx\,dt.
\end{aligned}
    \end{equation*}
Using the derivation followed in Nessayahu-Tadmor \cite{NT}, we have
\begin{equation*}
\begin{aligned}
u_{i+1/2}^{\,n+1}=& \dfrac{1}{\Delta x} \left( \int_{x_i}^{x_{i+1/2}} u(x, t^{n})dx+ \int_{x_{i+1/2}}^{x_{i+1}} u(x, t^{n})dx \right) + \dfrac{1}{\Delta x}\int_{t^n}^{t^{n+1}}\left[f(u(x_i,t)) - f(u(x_{i+1},t))\right]\;dt\\
& + \dfrac{1}{\Delta x}\int_{t^n}^{t^{n+1}}\!\int_{x_i}^{x_{i+1}}s\left(u(x,t), \varepsilon\right)\;dx\,dt.
\end{aligned}
    \end{equation*}
Upon using linear interpolation defined in \eqref{eq:linint}, we get
\begin{equation*}
\begin{aligned}
u_{i+1/2}^{\,n+1}=& \dfrac{1}{\Delta x} \left( \int_{x_i}^{x_{i+1/2}} L_i(x, t^{n})dx+ \int_{x_{i+1/2}}^{x_{i+1}} L_{i+1}(x, t^{n})dx \right) + \dfrac{1}{\Delta x}\int_{t^n}^{t^{n+1}}\left[f(u(x_i,t)) - f(u(x_{i+1},t))\right]\;dt\\
& + \dfrac{1}{\Delta x}\int_{t^n}^{t^{n+1}}\!\int_{x_i}^{x_{i+1}}s\left(u(x,t), \varepsilon\right)\;dx\,dt.\\
=& \dfrac{1}{2}\,\left(u_i^n + u_{i+1}^n\right)
  + \dfrac{1}{8}\,\left(u'_i - u'_{i+1}\right)
  + \dfrac{1}{\Delta x}\int_{t^n}^{t^{n+1}}\left[f(u(x_i,t)) - f(u(x_{i+1},t))\right]\;dt\\
& + \dfrac{1}{\Delta x}\int_{t^n}^{t^{n+1}}\!\int_{x_i}^{x_{i+1}}s\left(u(x,t), \varepsilon\right)\;dx\,dt.\\
\end{aligned}
\end{equation*}
Approximating the flux difference in time via the midpoint rule, we have the following. 
\begin{equation}\label{eq:mainnum}
\begin{aligned}
u_{i+1/2}^{\,n+1}=& \dfrac{1}{2}\,\left(u_i^n + u_{i+1}^n\right)
  + \dfrac{1}{8}\,\left(u'_i - u'_{i+1}\right)
  + \dfrac{\Delta t}{\Delta x}\,\left[f(u_i^{\,n+1/2}) - f(u_{i+1}^{\,n+1/2})\right]\\
& + \dfrac{1}{\Delta x}\int_{t^n}^{t^{n+1}}\!\int_{x_i}^{x_{i+1}}s\left(u(x,t), \varepsilon\right)\;dx\,dt.\\
\end{aligned}
\end{equation}

\subsection{Approximation of space-time integral of relaxed source term}
In the literature, various approaches have been explored to incorporate the stiff source term $s\left(u(x,t), \varepsilon\right)$ in the form of $\dfrac{1}{\varepsilon}g(u)$. In this case, the balance law becomes \begin{equation}\label{eq:balancemod}
u_t + f(u)_x = \dfrac{1}{\varepsilon} g(u),
\end{equation}
In the description of the scheme, to emphasise the stiff parameters, we consider the source term $s\left(u(x,t), \varepsilon\right)$ as $\dfrac{1}{\varepsilon}g(u).$ Now we define
\begin{equation*}
   I_g = \dfrac{1}{\varepsilon \Delta x }\int_{t^n}^{t^{n+1}}\!\int_{x_i}^{x_{i+1}}g\left(u(x,t)\right)\;dx\,dt,
\end{equation*}
within the numerical scheme \eqref{eq:mainnum}, using both explicit and implicit strategies.  We first propose a new way to approximate the spacetime integral $I_g$ by evaluating the spatial integration using the midpoint rule, resulting in the expression.
 \begin{equation*}
   I_g \approx  \dfrac{1}{\varepsilon} \int_{t^n}^{t^{n+1}} g(u(x_{i+1/2},t)) dt,
\end{equation*}
and later evaluate the time integral utilizing the trapezoidal rule,
\begin{equation*}
 I_g \approx  \dfrac{1}{\varepsilon} \int_{t^n}^{t^{n+1}} g(u(x_{i+1/2},t)) dt \approx \dfrac{1}{2 \varepsilon} \Delta t \left[ g(u_{i+1/2}^{n+1})+g(u_{i+1/2}^{n})\right],
\end{equation*}
where $g(u_{i+1/2}^{n})$ is approximated using a semi-implicit backward Taylor expansion,
\begin{eqnarray*}
\begin{aligned}
g(u_{i+1/2}^{n})&\approx g(u_{i+1/2}^{n+1})- \Delta t g_t (u_{i+1/2}^{n+1}),\\
&\approx g(u_{i+1/2}^{n+1})- \Delta t \left(\dfrac{\partial g}{\partial u} \dfrac{\partial u}{\partial t}\right) (u_{i+1/2}^{n+1}),\\
& \approx \ g(u_{i+1/2}^{n+1})- \Delta t \dfrac{\partial g}{\partial u}  (u_{i+1/2}^{n}) \dfrac{\partial u}{\partial t} (u_{i+1/2}^{n+1}),\\
& \approx g(u_{i+1/2}^{n+1})- \Delta t \dfrac{\partial g}{\partial u}  (u_{i+1/2}^{n}) \left(\dfrac{1}{\varepsilon} g(u_{i+1/2}^{n+1}) - f_x(u_{i+1/2}^n)\right)\\
&\approx g(u_{i+1/2}^{n+1})- \dfrac{\Delta t}{2} \left(\dfrac{\partial g}{\partial u}(u_{i}^n)+\dfrac{\partial g}{\partial u}(u_{i+1}^n)\right)\left(\dfrac{1}{\varepsilon} g(u_{i+1/2}^{n+1}) - \dfrac{f_{i+1}^n-f_{i}^n}{\Delta x}\right).
\end{aligned}
\end{eqnarray*}
Therefore, 
\begin{equation}\label{eq:source}
 I_g \approx \dfrac{\Delta t}{2 \varepsilon}  \left[ 2g(u_{i+1/2}^{n+1})- \dfrac{\Delta t}{2} \left(\dfrac{\partial g}{\partial u}(u_{i}^n)+\dfrac{\partial g}{\partial u}(u_{i+1}^n)\right)\left(\dfrac{1}{\varepsilon} g(u_{i+1/2}^{n+1}) - \dfrac{f_{i+1}^n-f_{i}^n}{\Delta x}\right)\right].
\end{equation}

In summary, by combining the piecewise‐linear reconstruction, midpoint flux integration, and semi-implicit treatment of the stiff source term, we arrive at a semi-implicit central differencing scheme, which we refer to as the \textit{second order central scheme with Euler backward Taylor numerical scheme} termed as CS-EBT2.

\medskip

\noindent\textbf{Predictor step:}\label{CSEBT2:scheme}  Solve for $u_i^{\,n+1/2}$ from
\begin{equation*}
 u_i^{\,n+1/2}= u_i^n
+ \dfrac{\Delta t}{2}\left[\dfrac{1}{\varepsilon} g(u_i^{\,n+1/2}) - \dfrac{f'_i}{\Delta x}\right].   
\end{equation*}

\noindent\textbf{Corrector step:}  Then update to $u_{i+1/2}^{n+1}$ by
\begin{equation*}
\begin{aligned}
u_{i+1/2}^{\,n+1}
&= \dfrac{1}{2}\left(u_i^n + u_{i+1}^n\right)
  + \dfrac{1}{8}\left(u'_i - u'_{i+1}\right)
  - \lambda\left[f(u_{i+1}^{\,n+1/2}) - f(u_i^{\,n+1/2})\right]\\
&
 +\dfrac{\Delta t}{2 \varepsilon}  \left[ 2g(u_{i+1/2}^{n+1})- \dfrac{\Delta t}{2} \left(\dfrac{\partial g}{\partial u}(u_{i}^n)+\dfrac{\partial g}{\partial u}(u_{i+1}^n)\right)\left(\dfrac{1}{\varepsilon} g(u_{i+1/2}^{n+1}) - \dfrac{f_{i+1}^n-f_{i}^n}{\Delta x}\right)\right],\\
\end{aligned}    
\end{equation*}
where $\lambda = \Delta t/ \Delta x$, $u'_i$ and $f'_i$ are computed using the MinMod limiter \eqref{eq:minmod}. 

\begin{remark}
The formulation presented in Equation \eqref{eq:source} follows a derivation approach analogous to that employed in \cite{Macca}. However, a key distinction lies in the underlying numerical framework: while \cite{Macca} develops the source term expression within the context of a finite-difference discretization, the current work adopts a finite-volume formulation. This difference significantly affects the treatment of spatial derivatives and fluxes, as the finite-volume method relies on integral conservation laws over control volumes and flux evaluations across cell interfaces, in contrast to the pointwise derivative approximations used in the finite-difference approach.  
\end{remark}

\begin{remark}
    Boundary conditions are enforced using ghost cells, as in standard finite-volume formulations. The reconstructed states in the ghost cells are prescribed according to the chosen boundary condition, and the predictor–corrector update is applied uniformly to both interior and boundary-adjacent cells, without introducing any additional staggering-related treatment.
\end{remark}

\begin{remark}
A relevant structural property of the proposed second-order semi-implicit scheme is that the implicit treatment of the source term is entirely local in space. In particular, no coupling between neighboring cells (or staggered control volumes) is introduced by the implicit part, and all spatial interactions arise solely from the explicitly treated flux terms. As a consequence, the nonlinear problems to be solved at each time step are fully decoupled across cells and involve only the components of the conserved vector, leading to a very efficient and easily parallelizable implementation.

In the case of systems of balance laws, the implicit treatment of the source term results in small nonlinear systems whose size is equal to the number of components of the state vector. These systems are solved independently in each control volume, under standard regularity assumptions on the source term. This locality property is nontrivial and is not always satisfied by semi-implicit discretizations available in the literature, where the implicit treatment may induce spatial coupling.

The extension of the present approach to higher-order accuracy is currently under investigation. Preliminary studies indicate that increasing the formal order while simultaneously preserving locality of the implicit solver, robustness, and stability is a delicate task, and improvements in one aspect may negatively affect others. Approaches based on Taylor or multi-derivative expansions appear as natural candidates, but their successful integration within the present framework is not straightforward. For this reason, the systematic development of higher-order variants (order greater than two) will be addressed in future research.
\end{remark}

\subsection{Extension to 2D:}
We now extend the proposed method to the two-dimensional case. Consider the hyperbolic relaxation equation
\begin{equation}
u_t + \partial_x f_1(u) + \partial_y f_2(u) = \frac{1}{\varepsilon} g(u),
\end{equation}
for $(x,y)\in\mathbb{R}^2$ and $t>0$, supplemented with the initial condition $u(x,y,0)=u_0(x,y)$. For the numerical approximation, we introduce a Cartesian grid with mesh sizes $\Delta x = x_{i+\frac12}-x_{i-\frac12}$, $\Delta y = y_{j+\frac12}-y_{j-\frac12}$ and time step $\Delta t = t^{n+1}-t^n$. The computational domain is partitioned into $N\times N$ uniform cells. The proposed CS-EBT2 scheme for the two-dimensional problem consists of the following predictor and corrector steps.\\
\noindent\textbf{Predictor step:}\label{CSEBT2:scheme2d}  Solve for $u_{i,j}^{\,n+1/2}$ from
\begin{equation*}
 u_{i,j}^{\,n+1/2}= u_{i,j}^n
+ \dfrac{\Delta t}{2}\left[\dfrac{1}{\varepsilon} g(u_{i,j}^{\,n+1/2}) - \dfrac{(f_1)^{\dagger}_{i,j}}{\Delta x}-\dfrac{(f_2)^{\ast}_{i,j}}{\Delta y}\right].   
\end{equation*}
\noindent\textbf{Corrector step:}  Then update to $u_{i+1/2,j+1/2}^{n+1}$ by
\begin{equation*}
\begin{aligned}
u_{i+1/2,j+1/2}^{\,n+1}
&= \left\langle\dfrac{1}{4}\left(u_{i,\cdot}^n + u_{i+1,\cdot}^n\right)
  + \dfrac{1}{8}\left(u^{\dagger}_{i,\cdot} - u^{\dagger}_{i+1,\cdot}\right)
  - \lambda\left[f_1(u_{i+1,\cdot}^{\,n+1/2}) - f_1(u_{i,\cdot}^{\,n+1/2})\right] \right.\\
& \qquad \left.+\dfrac{\Delta t^2}{4 \varepsilon}\left(\dfrac{\partial g}{\partial u}(u_{i,\cdot}^n)+\dfrac{\partial g}{\partial u}(u_{i+1,\cdot}^n)\right)\left(\dfrac{(f_1)_{i+1,\cdot}^n-(f_1)_{i,\cdot}^n}{\Delta x}\right)\right\rangle_{j+1/2}\\
+&\left\langle\dfrac{1}{4}\left(u_{\cdot,j}^n + u_{\cdot,j+1}^n\right)
  + \dfrac{1}{8}\left(u^{\ast}_{\cdot,j} - u^{\ast}_{\cdot,j+1}\right)
  - \mu\left[f_2(u_{\cdot,j+1}^{\,n+1/2}) - f_2(u_{\cdot,j}^{\,n+1/2})\right] \right.\\
& \qquad \left.+\dfrac{\Delta t^2}{4 \varepsilon}\left(\dfrac{\partial g}{\partial u}(u_{\cdot,j}^n)+\dfrac{\partial g}{\partial u}(u_{\cdot,j+1}^n)\right)\left(\dfrac{(f_2)_{\cdot,j+1}^n-(f_2)_{\cdot,j}^n}{\Delta y}\right)\right\rangle_{i+1/2}\\
& +\dfrac{\Delta t}{2 \varepsilon}  \left[ 2g(u_{i+1/2,j+1/2}^{n+1})-\Delta t \dfrac{\partial g}{\partial u}(u_{i+1/2,j+1/2}^n)\left(\dfrac{1}{\varepsilon} g(u_{i+1/2,j+1/2}^{n+1})\right)\right].
\end{aligned}    
\end{equation*}
Here,
\[
\langle u_{\cdot,j} \rangle_{i+\frac12} := \frac12 (u_{i,j}+u_{i+1,j}), \quad \langle u_{i,\cdot} \rangle_{j+\frac12} := \frac12 (u_{i,j}+u_{i,j+1}),
\]
denotes the staggered average in the $x$- and $y$-directions respectively. Moreover, we define $\lambda = \Delta t / \Delta x$ and $\mu = \Delta t / \Delta y$. Here, $u^{\dagger} = \partial u / \partial x$ and $u^{\ast} = \partial u / \partial y$. The quantities $u^{\dagger}_i$, $u^{\ast}_j$, $(f_1)^{\dagger}_{i,j}$, and $(f_2)^{\ast}_{i,j}$ are evaluated using the MinMod limiter defined in \eqref{eq:minmod}.

\section{Analysis of the scheme CS-EBT2}
In this section, we analyze the consistency and stability of the proposed numerical scheme. Specifically, we focus on the Jin-Xin  model with linear source to evaluate the applicability of the scheme and examine its consistency and stability characteristics within this framework.

\subsection{Consistency analysis.}
We aim to demonstrate that the proposed predictor-corrector numerical scheme \eqref{CSEBT2:scheme} in component wise to the prototypical $2 \times 2$ linear system of Jin-Xin model  
\begin{eqnarray}\label{6}
 \begin{cases}
u_t + v_x = 0,\\
v_t + u_x = -\dfrac{1}{\varepsilon}\left(v - au\right),
\end{cases}
\end{eqnarray}
with initial data
\begin{equation}\label{7}
    u(x,0) = u^0(x), \;\;\; v(x,0) = v^0(x),
\end{equation}
where $a \in \mathbb{R}$ satisfies $|a|<1$.  Note that system \eqref{6} can be written as $$U_t+F(U)_x=S(U,\varepsilon)$$ where $$U=[u, v]^T,\,\,F(U)=[v,u]^T,\,\,\,\text{and}\,\,\,\,S(U,\varepsilon)=\dfrac{1}{\varepsilon}g(U)=\dfrac{1}{\varepsilon}[0, -\left(v - au\right)]^T.$$
The proposed predictor-corrector numerical scheme becomes
\begin{equation*}
\begin{cases}
\begin{aligned}
 U_i^{\,n+1/2} &= U_i^n
+ \dfrac{\Delta t}{2}\left[\dfrac{1}{\varepsilon} g(U_i^{\,n+1/2}) - \dfrac{F'_i}{\Delta x}\right],\\
U_{i+1/2}^{\,n+1}
&= \dfrac{1}{2}\,(U_i^n + U_{i+1}^n)
  + \dfrac{1}{8}\,(U'_i - U'_{i+1})
  - \lambda\left[F(U_{i+1}^{\,n+1/2}) - F(U_i^{\,n+1/2})\right]\\
&
  +\dfrac{\Delta t}{2 \varepsilon}  \left[ 2g(U_{i+1/2}^{n+1})- \dfrac{\Delta t}{2} \left(\dfrac{\partial g}{\partial U}(U_{i}^n))+\dfrac{\partial g}{\partial U}(U_{i+1}^n))\right)\left(\dfrac{1}{\varepsilon} g(U_{i+1/2}^{n+1}) - \dfrac{F_{i+1}^n-F_{i}^n}{\Delta x}\right)\right].
\end{aligned}  
\end{cases}
\end{equation*}
Now we analyze the numerical scheme in three regimes such as fluid dynamic limit ($\varepsilon = 0$), thin regime ($\varepsilon \ll 1$) and rarefied regime ($\varepsilon = 1$).
\subsubsection*{\textbf{Fluid dynamics limit ($\varepsilon =0$).}}
In this case, from the second equation of the model \eqref{6}, we get $v=au,$ thus the first equation of the model \eqref{6} becomes $u_t+v_x=u_t+au_x=0.$ Note that in this case we do consistency analysis for the scalar equation $u_t+au_x=0,\,\,v=au$ rather than analyzing the consistency analysis for the system in component wise. The predictor-corrector numerical scheme associated with the model \eqref{6} is formulated as follows
\begin{eqnarray*}
\begin{cases}
\begin{aligned}
u_{i+1/2}^{n+1}=& \dfrac{1}{2}\left[u_i^n+u_{i+1}^n\right]+\dfrac{1}{8}\left[u'_i-u'_{i+1}\right]-\lambda \left[v_{i+1}^{n+1/2}- v_i ^{n+1/2}\right],\\
 v_{i+1/2}^{n+1}= & \dfrac{1}{2}\left[v_i^n+v_{i+1}^n\right]+\dfrac{1}{8}\left[v'_i-v'_{i+1}\right]-\lambda \left[u_{i+1}^{n+1/2}-u_i ^{n+1/2}\right] + \frac{\Delta t}{\varepsilon} \, (au_{i+1/2}^{n+1} -v_{i+1/2}^{n+1})\\
 &- \dfrac{\Delta t^2}{4} \left[\left(-\frac{1}{\varepsilon}-\frac{1}{\varepsilon}\right)\left(\frac{1}{\varepsilon}\left(au_{i+1/2}^{n+1} - v_{i+1/2}^{n+1}\right) - \dfrac{u_{i+1}^n-u_{i}^n}{\Delta x}\right)\right],
\end{aligned}
\end{cases}
\end{eqnarray*}
with the predictor step is
\begin{eqnarray*}
\begin{cases}
\begin{aligned}
u_{i}^{n+1/2} &= u_{i}^n - \dfrac{\Delta t}{2\Delta x} v'_{i}, \\
v_{i}^{n+1/2} &= v_{i}^n + \dfrac{\Delta t}{2} \left(\frac{1}{\varepsilon}\left(au_{i}^{n+1/2} - v_{i}^{n+1/2}\right)-\dfrac{u'_{i}}{\Delta x}\right).
\end{aligned}
\end{cases}
\end{eqnarray*}
To establish the consistency of the proposed scheme in the fluid regime, the relaxation terms are rescaled by multiplying $v_i^{\,n+1/2}$ in the predictor step by $\varepsilon$ and $v_{i+1/2}^{\,n+1}$ in the corrector step by $\varepsilon^2$. Passing to the limit $\varepsilon \to 0$ yields the equilibrium relations
\[
v_i^{\,n+1/2} = a\,u_i^{\,n+1/2} +\mathcal{O}(\varepsilon) \, \qquad
v_{i+1/2}^{\,n+1} = a\,u_{i+1/2}^{\,n+1} +\mathcal{O}(\varepsilon).
\]
Using the above equilibrium relations, the predictor-corrector numerical scheme takes the following reduced form in the fluid regime
\begin{eqnarray}\label{Ce0}
\begin{aligned}
\begin{cases}
 u_i^{\,n+1/2}= & u_i^n
- \dfrac{\Delta t}{2}\left(\dfrac{au'_i}{\Delta x}\right),\\
u_{i+1/2}^{\,n+1}
= & \dfrac{1}{2}\,(u_i^n + u_{i+1}^n)
  + \dfrac{1}{8}\,(u'_i - u'_{i+1})
  - \lambda\left[(au_{i+1}^{\,n+1/2}) - (au_i^{\,n+1/2})\right].
\end{cases}
\end{aligned}    
\end{eqnarray}
\begin{remark}
    We remark that the asymptotic-preserving property of the proposed scheme does not rely on the vanishing of a specific variable in the relaxation limit. More generally, the implicit treatment of the source term enforces the equilibrium condition $g(U)=0$ at the discrete level, independently of whether the relaxed variables disappear or remain physically relevant in the limiting model. The convergence of the method is therefore governed by the accuracy of the discretization of the limiting dynamics on the equilibrium manifold. Possible difficulties may only arise in cases where the equilibrium relation is ill-posed or its Jacobian is degenerate, which are issues related to the underlying physical model and not to the numerical strategy itself.
\end{remark}
Using Taylor series expansion for \eqref{Ce0}, we get
\begin{eqnarray}\label{CEN0}
\begin{aligned}
u_{i+1/2}^{\,n+1}
= & \dfrac{1}{2}\,(u_i^n + u_{i+1}^n)
  + \dfrac{1}{8}\,(u'_i - u'_{i+1})
  - \lambda a \left[(u_{i+1}^n
- \dfrac{a\Delta t}{2}\dfrac{u'_{i+1}}{\Delta x}) - (u_i^n
- \dfrac{a\Delta t}{2}\dfrac{u'_i}{\Delta x})\right]\\
=& \dfrac{1}{2}\left[u_i^n +u_{i+1}^n\right]+\dfrac{1}{8}\left[u'_i-u'_{i+1}\right]-\lambda a\left[u_{i+1}^n - u_{i}^n\right] + \dfrac{1}{2} a \lambda ^ 2 \left[v'_i-v'_{i+1}\right].\\
= & u_{i+1/2}^n - \lambda \Delta x a(u_x)_{i+1/2}^n + \dfrac{1}{2} \lambda^2 \Delta x^2 a^2 (u_{xx})_{i+1/2}^{n}+ \mathcal{O}(\textcolor{red}{\varepsilon},\Delta x^3).
\end{aligned}
\end{eqnarray}
Note that in the above equation, the last line of equation is derived by utilizing the following remark;
\begin{remark}\label{caremark}
 For a smooth function $u(x)$; the derivative approximations are of first-order accuracy and vary smoothly with respect to $x$. Thus, we have
\begin{eqnarray*}\label{17}
\begin{cases}
\begin{aligned}
&\dfrac{u'_i}{\Delta x} = \dfrac{\partial u}{\partial x} (x_i^n) + A(x_i) \Delta x +\mathcal{O}(\Delta x^2),\\
&\dfrac{v'_i}{\Delta x} = \dfrac{\partial v}{\partial x} (x_i^n) + B(x_i) \Delta x +\mathcal{O}(\Delta x^2),        
\end{aligned}
\end{cases}
\end{eqnarray*}
and hence 
\begin{eqnarray*}
\begin{cases}
\begin{aligned}
u^{n}_{i+1/2} &=\dfrac{1}{2}\left[u_i^n+u_{i+1}^n\right]+\dfrac{1}{8}\left[u'_i-u'_{i+1}\right]+ \mathcal{O}(\Delta x^3);\\
(u_x)_{i+1/2}^n &= \dfrac{1}{\Delta x} \left[u_{i+1}^n -u_{i}^n\right] + \mathcal{O}(\Delta x^2);\\
a(u_{xx})_{i+1/2}^n &= (v_{xx})_{i+1/2}^n = \dfrac{1}{\Delta x^2}\left[v'_{i+1}-v'_{i}\right] +\mathcal{O}(\Delta x);\\
\end{aligned} 
\end{cases}
\end{eqnarray*}
\end{remark}
Now, we derive the Taylor's expression for the exact solution $u(x,t)$ at $(x_{i+1/2},t^{n+1})$ 
and using the Cauchy-Kowaleski procedure to the model equation  $u_t+au_x=0,$ we have 
\begin{eqnarray}\label{CENE}
\begin{aligned}
 u(x_{i+1/2},t^{n+1})= & u(x_{i+1/2},t^{n})+\Delta t u_t(x_{i+1/2},t^{n})+\dfrac{\Delta t^2}{2} u_{tt}(x_{i+1/2},t^{n})+\mathcal{O}(\Delta t^3)\\
 =& u(x_{i+1/2},t^{n})-a\Delta t u_x(x_{i+1/2},t^{n})+\dfrac{a^2 \Delta t^2}{2} u_{xx}(x_{i+1/2},t^{n})+\mathcal{O}(\Delta t^3).
\end{aligned}
\end{eqnarray}
Upon comparing equations \eqref{CEN0} and \eqref{CENE}, we confirm that the proposed numerical scheme attains second-order accuracy.
\subsubsection*{\textbf{Thin regime ($\varepsilon \ll 1$).}}
In the system \eqref{6}, we are using Chapman-Enskog expansion for the second component $v$, i.e., $v = v^{(0)}+ \varepsilon v^{(1)} + \mathcal{O}(\varepsilon^2)$. The equilibrium states are $v^{(0)} = au +\mathcal{O}(\varepsilon)$ and $v^{(1)} = (1-a^2)u_x+\mathcal{O}(\varepsilon)$. After plugging the values of $v^{(0)}$ and $v^{(1)}$, we have $v= au -\varepsilon(1-a^2)u_x + \textcolor{red}{\mathcal{O}(\varepsilon^2)}$. Then the system \eqref{6} can be written as 
\begin{eqnarray}\label{22}
\begin{cases}
     \begin{aligned}
       &u_t +v_x = 0,\\
       &v=au-\varepsilon(1-a^2)u_{x}.
     \end{aligned}
     \end{cases}
\end{eqnarray}
Again note that the numerical scheme will be applied to scalar model rather than applying to the system in component wise. Utilizing the Remark \eqref{caremark}, the numerical solution becomes
\begin{eqnarray}\label{NCEl1}
\begin{aligned}
 u_{i+1/2}^{n+1} &= \dfrac{1}{2}\left[u_i^n +u_{i+1}^n\right]+\dfrac{1}{8}\left[u'_i-u'_{i+1}\right]-\lambda\left[v_{i+1}^{n+1/2} -v_{i}^{n+1/2}\right]\\
&= \dfrac{1}{2}\left[u_i^n +u_{i+1}^n\right]+\dfrac{1}{8}\left[u'_i-u'_{i+1}\right]\\
& \qquad -\lambda\left[a u_{i+1}^{n+1/2} -a u_{i}^{n+1/2} -\varepsilon (1-a^2)(u_x)_{i+1}^{n+1/2}+\varepsilon (1-a^2)(u_x)_{i}^{n+1/2}\right] \\
& = u_{i+1/2}^n - a \lambda \Delta x (u_x)^n_{i+1/2} +a^2 \dfrac{\lambda^2 \Delta x^2}{2}(u_{xx})^n_{i+1/2}\\ 
& \qquad  + \varepsilon (1-a^2) \lambda \Delta x ((u_{xx})^n_{i+1/2} - a\lambda \Delta x (u_{xxx})^n_{i+1/2})+\mathcal{O}(\Delta x^3, \varepsilon^2, \varepsilon\Delta x^2).
    \end{aligned}
\end{eqnarray}
Recall Taylor's expression for the exact solution $u(x_{i+1/2},t^{n+1})$ and utilizing the Cauchy-Kowaleski procedure to the model equation $u_t+v_x=0,\,\,v=au-\varepsilon(1-a^2)u_{x}$ we have 
\begin{eqnarray*}\label{CENEl1}
\begin{aligned}
  u(x_{i+1/2},t^{n+1})= & u(x_{i+1/2},t^{n})+\Delta t u_t(x_{i+1/2},t^{n})+\dfrac{\Delta t^2}{2} u_{tt}(x_{i+1/2},t^{n})+\mathcal{O}(\Delta t^3)\\
  =& u(x_{i+1/2},t^{n})-\Delta t v_x(x_{i+1/2},t^{n})+\dfrac{ \Delta t^2}{2} v_{xx}(x_{i+1/2},t^{n})+\mathcal{O}(\Delta t^3)\\
 = & u(x_{i+1/2},t^{n}) - \Delta t au_x(x_{i+1/2},t^{n}) + \dfrac{1}{2} \Delta t^2 a^2 u_{xx}(x_{i+1/2},t^{n})\\
 & + \varepsilon \Delta t (1-a^2)\left(u_{xx}(x_{i+1/2},t^{n}) - a \Delta t u_{xxx}(x_{i+1/2},t^{n}) \right)+ \mathcal{O}(\Delta t^3, \varepsilon ^2, \varepsilon\Delta x^2),
\end{aligned}
\end{eqnarray*}
which matches the numerical behavior \eqref{NCEl1} up to errors of order $\mathcal{O}(\Delta x^3,\varepsilon^2,\varepsilon\Delta x^2)$.  Thus, consistency is maintained upto second order.

\begin{remark}
The rarefied regime ($\varepsilon = 1$ or $\varepsilon=\mathcal{O}(1)$), which is commonly addressed in the literature, is not given importance in our derivation. For brevity in the presentation, we perform consistency analysis for this case also and we arrive a conclusion that the proposed numerical scheme achieves a first-order accuracy for the second term $v$ whereas second-order accuracy for the first term $u$. Even, the method remains stable, it is designed to be uniformly second-order asymptotically-accurate. Indeed, the additional dissipation induced by the implicit relaxation update may lead to an effective first-order behavior for the coupled system when the solution significantly departs from the equilibrium manifold. However, in the absence of stiffness in the source term, we advise against employing any implicit treatment, as it may introduce unnecessary complexity without providing computational benefit.  
\end{remark}

 \subsubsection*{\textbf{Rarefied regime ($\varepsilon =1$).}}
 In this case, the Jin-Xin model \eqref{6}, can be written as 
 \begin{eqnarray}
 \begin{cases}
u_t + v_x = 0,\\
v_t + u_x = au - v,
\end{cases}
\end{eqnarray}
where $U=(u,\;v)^T$, $F(U)= (v,\;u)^T$ and $S(U, \varepsilon)=\dfrac{1}{\varepsilon}g(U)= (0,\;au - v)^T$. The corrector-step of the numerical scheme for the above given system 
\begin{eqnarray*}
\begin{cases}
\begin{aligned}
u_{i+1/2}^{n+1}=& \dfrac{1}{2}\left[u_i^n+u_{i+1}^n\right]+\dfrac{1}{8}\left[u'_i-u'_{i+1}\right]-\lambda \left[v_{i+1}^{n+1/2}- v_i ^{n+1/2}\right],\\
 v_{i+1/2}^{n+1}= & \dfrac{1}{2}\left[v_i^n+v_{i+1}^n\right]+\dfrac{1}{8}\left[v'_i-v'_{i+1}\right]-\lambda \left[u_{i+1}^{n+1/2}-u_i ^{n+1/2}\right] + \Delta t \, (au_{i+1/2}^{n+1} -v_{i+1/2}^{n+1})\\
 &- \dfrac{\Delta t^2}{4} \left[\left(-1-1\right)\left(au_{i+1/2}^{n+1} - v_{i+1/2}^{n+1} - \dfrac{u_{i+1}^n-u_{i}^n}{\Delta x}\right)\right],
\end{aligned}
\end{cases}
\end{eqnarray*}
with the predictor step is
\begin{eqnarray*}
\begin{cases}
\begin{aligned}
u_{i}^{n+1/2} &= u_{i}^n - \dfrac{\Delta t}{2\Delta x} v'_{i}, \\
v_{i}^{n+1/2} &= v_{i}^n + \dfrac{\Delta t}{2} \left(au_{i}^{n+1/2} - v_{i}^{n+1/2}-\dfrac{u'_{i}}{\Delta x}\right).
\end{aligned}
\end{cases}
\end{eqnarray*}
\begin{remark}\label{remark2}
We make use of the following expressions further in the derivation for a smooth function $u(x)$ and the derivative approximations using Taylor's expansion with respect to $x$ at the point $x_{i+1/2}$;
 \begin{eqnarray}\label{18}
 \begin{cases}
\begin{aligned}
&u_{i+1}^n= u_{i+1/2}^n + \dfrac{1}{2}(u_x)_{i+1/2}^n \Delta x + \dfrac{1}{8}(u_{xx})_{i+1/2}^n \Delta x^2 +\mathcal{O}(\Delta x^3),\\
&u_{i}^n= u_{i+1/2}^n - \dfrac{1}{2}(u_x)_{i+1/2}^n\Delta x + \dfrac{1}{8}(u_{xx})_{i+1/2}^n \Delta x^2 +\mathcal{O}(\Delta x^3),\\
&\dfrac{u'_{i+1}}{\Delta x} = (u_x)_{i+1/2}^n + \left(A+\dfrac{1}{2}(u_{xx})_{i+1/2}^n\right) \Delta x +\mathcal{O}(\Delta x^2),\\
&\dfrac{u'_i}{\Delta x} = (u_x)_{i+1/2}^n + \left(A-\dfrac{1}{2}(u_{xx})_{i+1/2}^n\right) \Delta x +\mathcal{O}(\Delta x^2),
\end{aligned}
\end{cases}
\end{eqnarray}
 where $u_{i+1/2}^n=u(x_{i+1/2},t^n)$ and $A = A(x_{i+1/2},t^n)$.    
\end{remark}
Upon Taylor's expansion of predictor-corrector steps and using Remark \eqref{remark2}, we get the following
\begin{eqnarray*}
\begin{cases}
\begin{aligned}
u_{i+1/2}^{n+1} = &u_{i+1/2}^n - \lambda \Delta x (v_x)_{i+1/2}^n + \dfrac{\lambda ^2 \Delta x^2}{2} \left((v_x)_{i+1/2}^n - a(u_x)_{i+1/2}^n +(u_{xx})_{i+1/2}^n\right) +  \mathcal{O}(\Delta x^3),\\
v_{i+1/2}^{n+1} =&v_{i+1/2}^n + \lambda \Delta x \left( a u_{i+1/2}^n- v_{i+1/2}^n -(u_x)_{i+1/2}^n \right) \\ & + \dfrac{\lambda^2 \Delta x^2}{2}\left((v_{xx})_{i+1/2}^n -au_{i+1/2}^n +v_{i+1/2}^n + (u_x)_{i+1/2}^n - 2a(v_x)_{i+1/2}^n \right) + \mathcal{O}(\Delta x^3).
\end{aligned}
\end{cases}
\end{eqnarray*}
Expanding the exact solution using Taylor's expansion yields 
\begin{eqnarray*} \label{CENE1}
\begin{cases}
\begin{aligned}
u(x_{i+1/2},t^{n+1})=& u(x_{i+1/2},t^{n}) - \Delta t v_x(x_{i+1/2},t^{n}) \\
& \quad + \dfrac{\Delta t^2}{2}\left(u_{xx}(x_{i+1/2},t^{n}) + v_x(x_{i+1/2},t^{n}) -au_x(x_{i+1/2},t^{n})\right) + \mathcal{O}(\Delta t^3),\\
v(x_{i+1/2},t^{n+1})= & v(x_{i+1/2},t^{n}) - \Delta t \left(u_x(x_{i+1/2},t^{n}) +v(x_{i+1/2},t^{n}) -au(x_{i+1/2},t^{n}) \right) \\ 
& + \dfrac{\Delta t^2}{2} \left(v_{xx}(x_{i+1/2},t^{n})+ u_x(x_{i+1/2},t^{n}) - av_x(x_{i+1/2},t^{n}) \right.\\
&\left.-au(x_{i+1/2},t^{n}) + v(x_{i+1/2},t^{n}) \right) + \mathcal{O}(\Delta t^3).
\end{aligned}
\end{cases}
\end{eqnarray*}
where $u_t= - v_x$, $u_{tt}=u_{xx}-v_x+au_x$, $v_t= au-v-u_x$ and $v_{tt}=v_{xx}-av_x+u_x-au+v$.  By comparing the numerical solution  with exact solution, it is observed that the proposed scheme achieves second-order accuracy in the first component $u$, while the second component $v$,
\begin{align*}
v_{i+1/2}^{n+1}-v(x_{i+1/2}, t^{n+1})= -\dfrac{\lambda^2 \Delta x^2}{2} a(v_x)_{i+1/2}^n+\mathcal{O}(\Delta x^3), \;\; (\text{where} \;\; \Delta t = \lambda \Delta x)
\end{align*}
achieves first-order accuracy.

\subsection{Stability Analysis}
We now turn to a linear stability study of the scheme, examining its $A-$ and $L-$stability \cite{Liotta} properties when applied to systems of ordinary differential equations.\footnote{We transform the PDEs $U_t + F(U)_x = S(U,\varepsilon)$ into $y' = f(y) + g(y)$ \eqref{26}, since no  formal analysis for a general non-linear ODEs exists so far, we focus on a standard linear ODEs $y' = \lambda_1y + \lambda_2y$ \eqref{28} for which a general theory exists.  In \cite{Macca}, the analysis was much easier since the scheme is finite difference instead of finite volume as we introduced.} In this regard, let us recall the definition of $A-$ and $L-$stability. Consider a general system of ODEs of the form
 \begin{align}\label{24}
     & y'=f(y),\,\,\,\, y(0) = y_0,
 \end{align}
where $y\in \mathbb{R}^N$. In order to study the $A-$stabilily of a one-step method, let us consider the scheme applied to a test problem
\begin{align}\label{25}
    y'=\lambda y,\quad y(0)=1,
\end{align}
in which $y,\lambda\in\mathbb{C}$ with ${\rm Re}(\lambda)<0$. The exact solution of the \eqref{25} problem is  $y(t)=e^{\lambda t}$, and at time $n \Delta t$, is
\begin{align*}
    y(n \Delta t)=e^ {\lambda n \Delta t},
\end{align*}
and its absolute value decreases monotonically in time. Applying a one-step method to equation \eqref{25}, the numerical solution is
\begin{equation*}
    y_{n+1}=\Phi(z)y_n,
\end{equation*}
where $z = \lambda \Delta t$. Here $\Phi(z)$ is called absolute stability function\footnote{It is well known, that the region of the complex plane
\begin{equation*}
S_{A}=\{z \in \mathbb{C}:|\Phi(z)| \leq 1\}
\end{equation*}
is called absolute stability region.}. In particular, a one-step scheme is said to be $A-$stable if its region of absolute stability contains the complex half plane $\mathbb{C}^{-}=\{z \in \mathbb{C}: \rm Re(z) \leq 0\}$. A scheme is said to be $L-$stable if it is $A-$stable and $\lim _{z \rightarrow \infty} \Phi(z)=0$. The $L-$stability property is fundamental requirement that a scheme should have when dealing with stiff systems. 

Our scheme targets problems have the form
\begin{align}\label{26}
&y' =f(y)+g(y),\,\,\,\,y(0) =y_{0},
\end{align}
by treating $f$ explicitly and $g$ implicitly.
\par The CS-EBT2 scheme is suited for systems of the form \eqref{26}, where it is possible to identify a stiff part $g$, and a non-stiff part $f$. Our scheme, in an opportune semi-discrete way, can be written in the form
\begin{eqnarray}\label{27}
\begin{cases}
\begin{aligned}
    &y_{n+1/2} =y_{n}+\dfrac{\Delta t}{2} f\left(y_{n}\right)+\dfrac{\Delta t}{2} g\left(y_{n+1/2}\right), \\
   & y_{n+1} =y_{n}+\Delta t f\left(y_{n+1 / 2}\right)+ \Delta t \, g\left(y_{n+1}\right)- \dfrac{\Delta t^2}{4}\left(\dfrac{\partial g}{\partial y} (y_{n})\right)\left( g(y_{n+1}) + f(y_n)\right).
\end{aligned}
\end{cases}
\end{eqnarray}
\par Since the main novelty of this scheme lies in the semi-implicit treatment of the source term, we begin the $L-$stability analysis by focusing on the simpler linear system. Then, let us focus on the $A-$stability analysis, considering the linear system with $\lambda_1\ll \lambda_2$
\begin{align}\label{28}
& y' =\lambda_{1} y+\lambda_{2} y,\,\,\,\,y(0)=1,
\end{align}
thus numerical scheme (\ref{27}) predictor step is
\begin{eqnarray*}
y_{n+1 / 2} = y_{n}+\dfrac{\Delta t}{2} \lambda_1 y_{n}+ \dfrac{\Delta t}{2} \lambda_2 y_{n+1/2}=\left( \dfrac{1+ \dfrac{\Delta t}{2}\lambda_{1}}{1 - \dfrac{\Delta t}{2}\lambda_{2}}\right)  y_{n},
\end{eqnarray*}
and corrector step is
\begin{eqnarray*}
y_{n+1} =y_{n}+\Delta t \lambda_{1}\left( \dfrac{1+ \dfrac{\Delta t}{2}\lambda_{1}}{1 - \dfrac{\Delta t}{2}\lambda_{2}}\right)  y_{n}+\Delta t \lambda_2 y_{n+1} - \dfrac{\Delta t^2}{4} \lambda_2 \left( \lambda_2 y_{n+1} +\lambda_1 y_n\right),
\end{eqnarray*}
hence 
\begin{eqnarray*}
y_{n+1}=\left[\dfrac{1+ z_{1} \left(\dfrac{1+\dfrac{z_{1}}{2}}{1-\dfrac{z_{2}}{2}}\right)- \dfrac{z_1 z_2}{4}}{1- z_2 +\dfrac{z_2^2}{4}}\right]y_n.
\end{eqnarray*}
Here, the function 
\begin{equation*}
\Phi\left(z_{1}, \, z_{2}\right)=\left[1+ z_{1} \left(\dfrac{2+z_{1}}{2-z_{2}}\right)- \dfrac{z_1 z_2}{4}\right]\left(1- z_2 +\dfrac{z_2^2}{4}\right)^{-1},
\end{equation*}
is the absolute stability function with $z_{1}=\lambda_{1} \Delta t, \; z_{2}=\lambda_{2} \Delta t.$ The absolute stability region $S_{A}$ associated with scheme (\ref{27}) is defined as
\begin{equation*}
S_{A}=\left\{\left(z_{1}, \, z_{2}\right) \in \mathbb{C}^{2}:\left|\Phi\left(z_{1}, \, z_{2}\right)\right| \leq 1\right\}.    
\end{equation*} 
Since $\lambda_1\ll\lambda_2$, we can consider $z_{1}=0$ and the function reduces to
\begin{equation}\label{st_region}
\Phi\left(0,\, z_{2}\right)= \dfrac{1}{\left( 1- z_2 +\dfrac{z_2^2}{4}\right)}.
\end{equation}
The proposed scheme is $A-$stabile if the absolute stability function satisfies $|\Phi(0,\, z_2)| \leq 1$, for all $z_2 \in \mathbb{C}^-$. Furthermore, $\Phi(0,\, z_2) \to 0$ as $z_2 \to \infty$, indicating that the CS-EBT2 scheme, for the system \eqref{6} is also $L-$stable. It is evident that the region $S_A$ does not contain the set $\mathbb{C}^-
\times \mathbb{C}^-$. Our goal is to show that there exist two regions of the complex plane, $S_1 \subset \mathbb{C}, S_2 \subset \mathbb{C} $, with the following properties:
\begin{align}
    & S_A \supset S_1 \times S_2,\\
    & S_2 \supset \mathbb{C}^-,
\end{align}
The non-uniqueness of such pairs of sets, should they exist, is readily observed. Our objective is to determine, in explicit form, the largest set $S_1$ such that the corresponding set $S_2$ encompasses the entire left half of the complex plane, i.e., $S_2 \supset \mathbb{C}^-$. This admissible domain is formally defined by
\begin{equation*}
S_{1}=\left\{z_{1} \in \mathbb{C}: \max_{z_2 \in \mathbb{C}^-} \left|\Phi\left(z_{1}, \, z_{2}\right)\right| \leq 1\right\}.    
\end{equation*}
 \begin{lemma} 
     Following \cite{Liotta}, for any fixed value of $z_{1} \in \mathbb{C}$, the modulus of the function $\Phi\left(z_{1}, \, z_{2}\right)$ assumes its maximum value in the complex half plane $\mathbb{C}^-$ for some $z_2$ belonging to the imaginary axis.
 \end{lemma}
\begin{proof}
     The stability function $\Phi(z_1,z_2)$ admits a rational representation of the form
    \[
     \Phi(z_1,\, z_2)=\frac{P(z_1,\,z_2)}{Q(z_2)},
    \]
    where $P$ and $Q$ are polynomials. For a fixed parameter $z_1\in\mathbb{C}$, $\Phi$ is regarded as a function of the complex variable $z_2$. The A-stability of the numerical scheme ensures that its region of absolute stability contains the closed left half-plane. Accordingly, the stability function satisfies the uniform bound
    \[
    |\Phi(z_1,z_2)|\leq 1 \quad \text{for all} \; \;z_2\in\mathbb{C}^-
    \]
    Since $\Phi(z_1,z_2)$ is rational in $z_2$, this boundedness excludes the presence of poles in $\mathbb{C}^-$; otherwise, the modulus of $\Phi$ would diverge in any neighborhood of such a singularity. Hence, for each fixed $z_1$, $\Phi(z_1,z_2)$ is holomorphic in $z_2$ throughout $\mathbb{C}^-$.
    As a consequence of holomorphicity, the maximum modulus principle applies and implies that the maximum of $|\Phi(z_1,z_2)|$ over $\mathbb{C}^-$ cannot be attained in the interior of the domain. Therefore, any point at which the maximum is achieved must lie on the boundary $\partial\mathbb{C}^-$, with the maximum at infinity excluded by the decay at infinity.
    This is guaranteed by the L-stability of the scheme, which yields the decay property
    \[
     \lim_{z_2\to\infty}\Phi(z_1,z_2)=0.
    \]
    Thus, $|\Phi(z_1,z_2)|$ vanishes as $|z_2|\to\infty$ within $\mathbb{C}^-$, excluding infinity as a possible maximizer. It follows that the maximum of $|\Phi(z_1,z_2)|$ over $\mathbb{C}^-$ is attained on the boundary, which coincides with the imaginary axis. This completes the proof.\end{proof}

\noindent  Based on the preceding lemma, the set $S_1$ can be characterized as
\begin{equation}\label{S1}
S_{1}=\left\{z_{1} \in \mathbb{C}: \max_{y \in \mathbb{R}} \left|\Phi\left(z_{1}, \, i y\right)\right| \leq 1\right\}.   
\end{equation}
The corresponding boundary of this region is defined by the condition
\begin{equation*}
\partial S_{1}=\left\{z_{1} \in \mathbb{C}: \max_{y \in \mathbb{R}} \left|\Phi\left(z_{1}, \, i y\right)\right| = 1 \right\}.    
\end{equation*}

\begin{figure}[!ht]
     \centering
     \begin{subfigure}[b]{0.49\textwidth}
         \centering   \includegraphics[trim=13.3cm 0.21cm 13.4cm 0.19cm, clip=true,width=\linewidth]{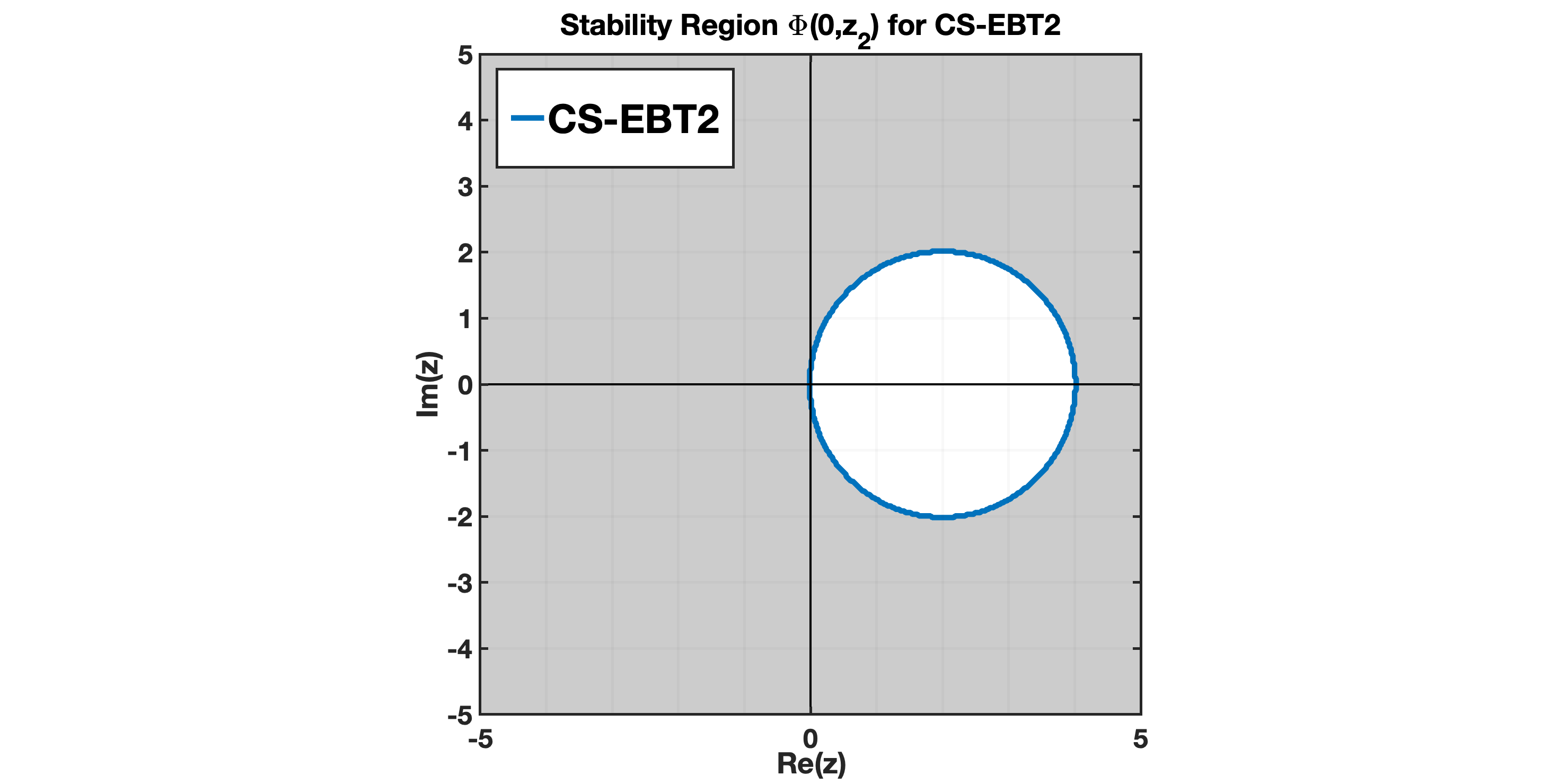}
    \caption{Stability regions (grey shaded areas) for the CS-EBT2 scheme.}
         \label{(fig:s1)}
     \end{subfigure}
     \hfill
     \begin{subfigure}[b]{0.49\textwidth}
         \centering   \includegraphics[trim=12.5cm 0.27cm 13.0cm 0.05cm, clip=true,width=\linewidth]{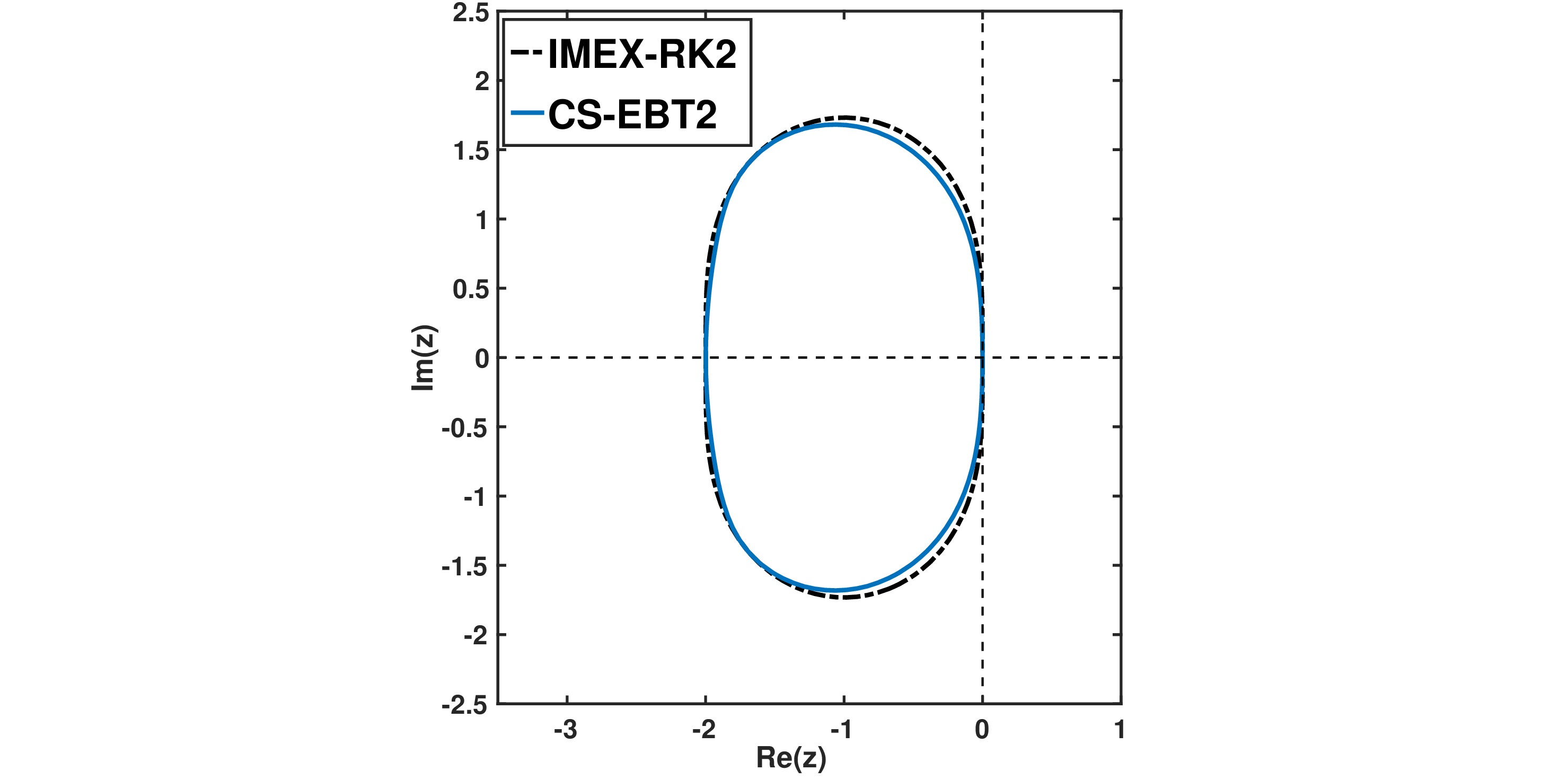} \caption{Stability regions $S_1$ for IMEX-RK2 and CS-EBT2 scheme.}
         \label{(fig:s2)}
     \end{subfigure}
   \caption{Stability region of the proposed scheme for $\Phi(0,z_2)$ (left) and $S_1$ (right). }
     \label{fig:stability_region}
\end{figure}

\hyperref[(fig:s1)]{Figure~\ref*{(fig:s1)}} shows the stability region of the proposed CS–EBT2 scheme, where the analysis is restricted to the source term $\Phi(0,\, z_{2})$ by setting $z_1 = 0$. In this case, the shaded area highlights the corresponding stability domain of the scheme. For further details on the definition of $S_2$ and the computation of $S_1$, refer to \cite{SB}. \hyperref[(fig:s2)]{Figure~\ref*{(fig:s2)}} presents the stability regions $S_1$ for both IMEX–RK schemes and the proposed method, providing a comparative analysis. The results demonstrate that the stability region of the proposed scheme is very close to that of IMEX–RK2 \cite{Boscarino-Filbet,MaccaRussoBumi,Macca2024}. 

As IMEX-RK2, we consider the scheme defined by the following double Butcher tableau :
\begin{equation}
\label{tableau}
\begin{array}{c|cc}
     & 0 &  \\
    c & c & 0\\ \hline
     & 1 - \gamma & \gamma
\end{array}
\hspace{3cm}
\begin{array}{c|cc}
    \gamma & \gamma &  \\
    1 & 1 - \gamma & \gamma\\ \hline
     & 1 - \gamma & \gamma
\end{array}
\end{equation}
where $\gamma = 1 - \frac{1}{\sqrt{2}}$ and $c = \frac{1}{2\gamma}$.

\section{Numerical Examples}\label{sec:Numerical_Experiments}
We assess the proposed Central Scheme with Modified Trapezoidal Rule (CS-EBT2) on a set of standard benchmark problems for hyperbolic balance laws with stiff relaxation, in one and two space dimensions. The test suite is selected to cover: (i) smooth well-prepared data (accuracy and convergence), (ii) unprepared data and discontinuities (robustness in the presence of sharp gradients), and (iii) very stiff regimes ($\varepsilon \ll 1$) (performance in the asymptotic limit) \textit{asymptotic-preserving} and \textit{asymptotic-accurate} schemes \cite{IMEX_book}.

In all experiments, the time step is chosen according to
\begin{equation}
\Delta t^n = \mathrm{CFL}\,\frac{\Delta x}{\Lambda_{\max}^n},
\end{equation}
where $\Lambda_{\max}^n$ denotes the maximum spectral radius of the Jacobian $\partial F/\partial U$ at time $t^n$; in two dimensions, $\Delta x$ denotes the minimum mesh spacing. Unless otherwise stated, we use uniform Cartesian grids and specify within each subsection the mesh size, boundary conditions, stiffness parameter $\varepsilon$, and final time.

The numerical campaign is organized as follows:
\begin{itemize}
\item Jin-Xin relaxation model: smooth (well-prepared/unprepared) data and Riemann-type profiles.
\item Shallow water model with relaxation: smooth and non-smooth configurations.
\item Broadwell model: convergence tests and a canonical discontinuous setup in the conserved variables.
\item Euler equations with heat transfer: discontinuous unprepared data in the stiff heat-transfer regime.
\item Euler system with stiff friction: discontinuous profiles for large friction coefficients.
\item Two-dimensional relaxation benchmarks: a 2D Jin-Xin test and 2D Euler Jin-Xin-type relaxation Riemann problems.
\end{itemize}

Whenever available, we use exact solutions; otherwise, we compare against suitable reference solutions computed on refined meshes.

Although the scheme is formally second-order accurate, its main strength lies in the design of a robust semi-implicit central discretization in which the stiff source term is treated through only two implicit stages. This represents a notable simplification compared to several second-order asymptotic-accurate schemes available in the literature, which typically rely on three or more implicit steps.

\subsection{Jin-Xin Model} We begin our analysis with the Jin–Xin model \eqref{6}, for which a smooth initial condition is considered over the spatial domain $[0,1]$, subject to periodic boundary conditions,
\begin{equation*}
\begin{cases}
u_t + v_x = 0,\\
v_t + u_x = -\dfrac{1}{\varepsilon}\,(v-a u).
\end{cases}
\end{equation*}

\noindent{\textbf{Smooth case (Well-Prepared):}}  The initial data is given by
\begin{eqnarray} \label{smoothdata:XinJin}
\begin{cases}
\begin{aligned}
 &u(x,0) = \sin(2\pi x),\\
&v(x,0) = au(x,0),
\end{aligned}
\end{cases}
\end{eqnarray}
where the parameter $a$ is fixed at $0.7$. This configuration ensures a smooth evolution of the solution, providing a suitable setting to assess the accuracy and stability of the proposed numerical scheme in the absence of discontinuities. The computations are performed using $N=320$ uniformly distributed grid points\footnote{In order to enhance the clarity of the plots, only about 40–60\% of the data points were retained for display.} he CFL number was set to $1/3$ in this test in order to reproduce the reference results available in the literature \cite{Liotta}. However, we remark that the proposed scheme remains stable and accurate also for significantly larger CFL numbers (up to CFL $=0.9$), as confirmed by additional numerical experiments. The solutions in this test are computed up to the final time $T=0.35$.
\begin{figure}[!ht]
     \centering
    \centering
     \begin{subfigure}[b]{0.5\textwidth}
         \centering
         \includegraphics[width=\textwidth]{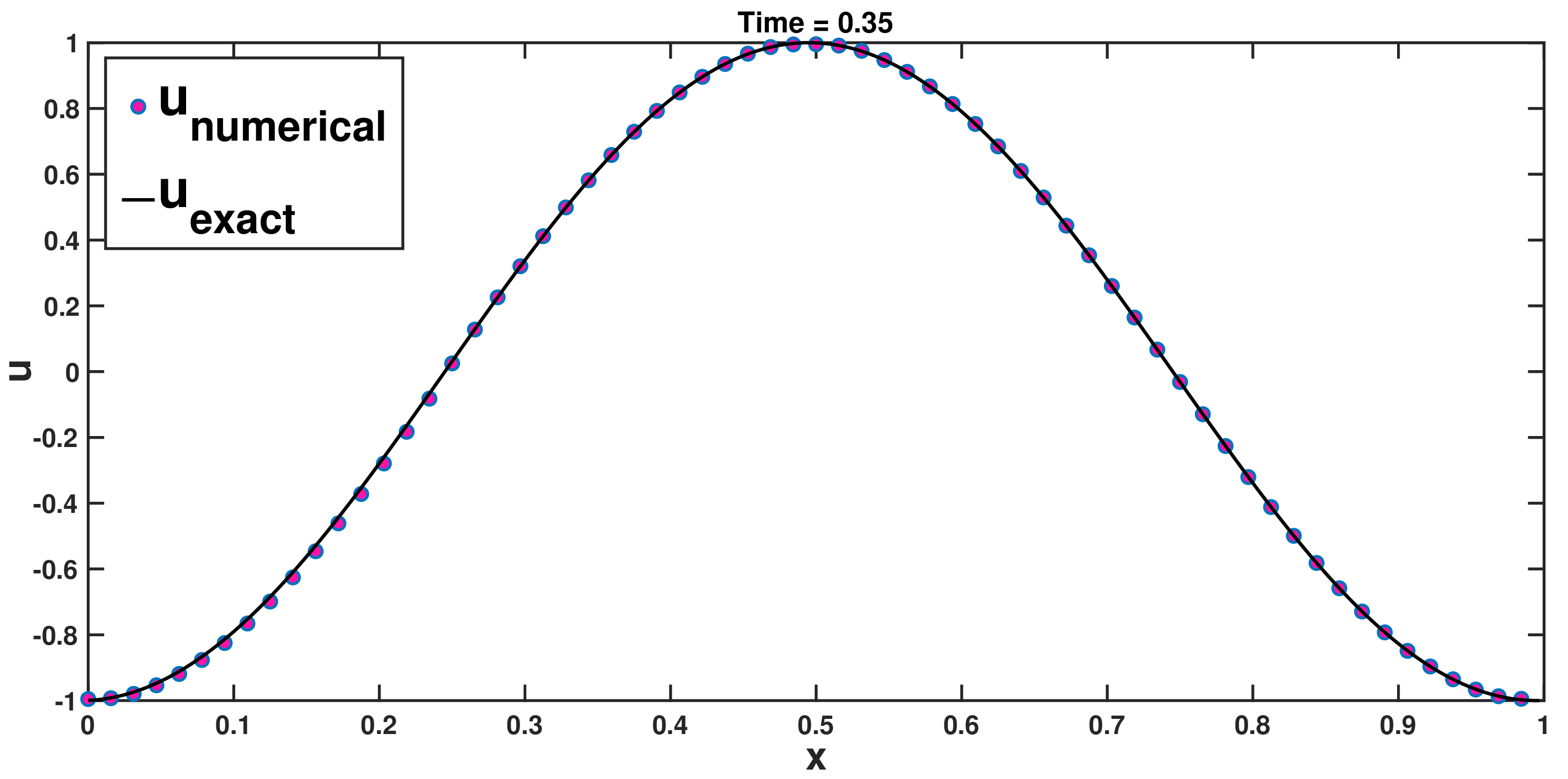}
         \caption{$u$ and $\varepsilon = 10^{-7}$}
         \label{(12a))}
     \end{subfigure}
     \hfill
     \begin{subfigure}[b]{0.48\textwidth}
         \centering
         \includegraphics[width=\textwidth]{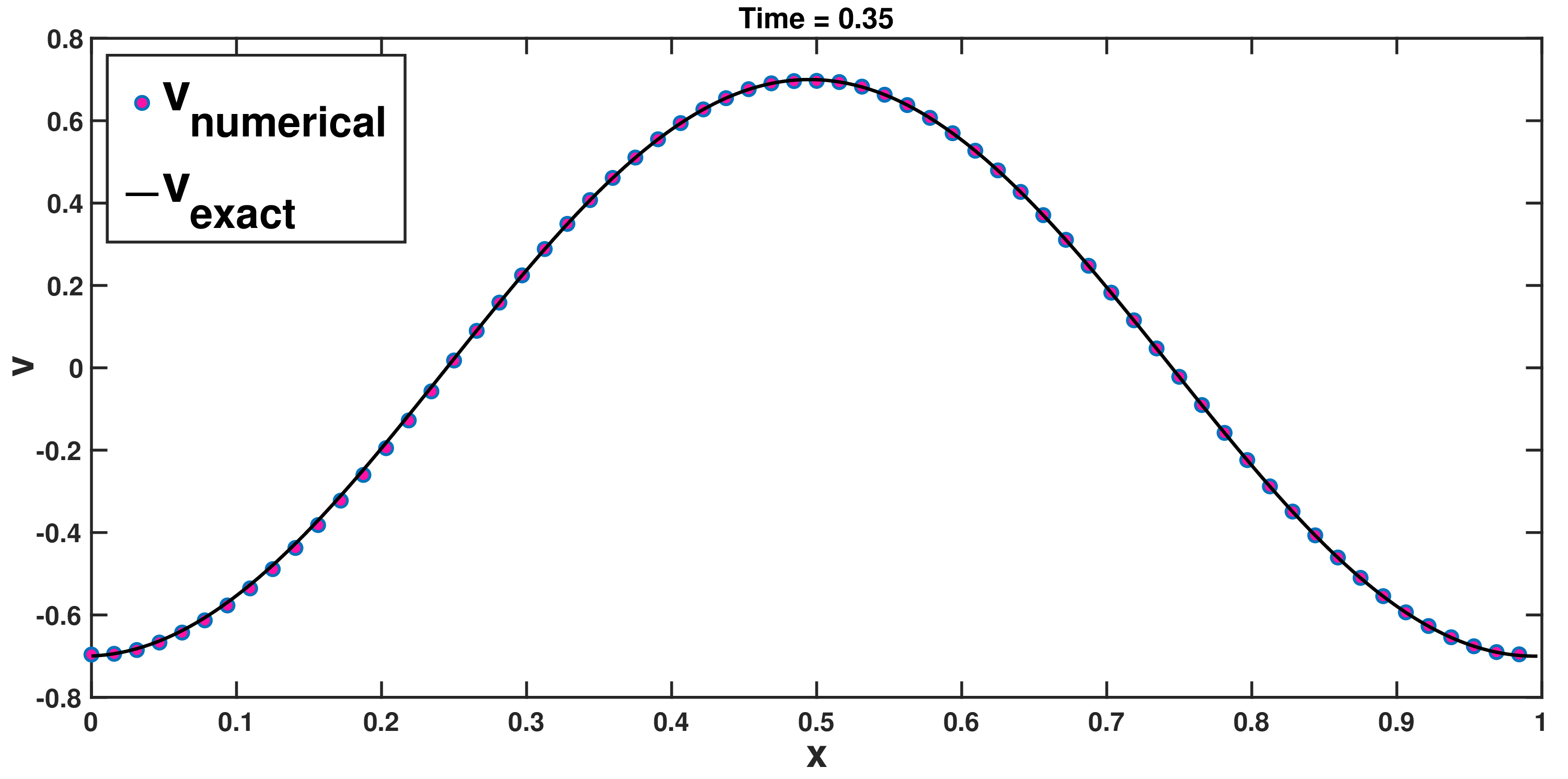}
         \caption{$v$ and $\varepsilon = 10^{-7}$}
         \label{(12b))}
     \end{subfigure}
     \medskip
    \centering
     \begin{subfigure}[b]{0.5\textwidth}
         \centering
         \includegraphics[width=\textwidth]{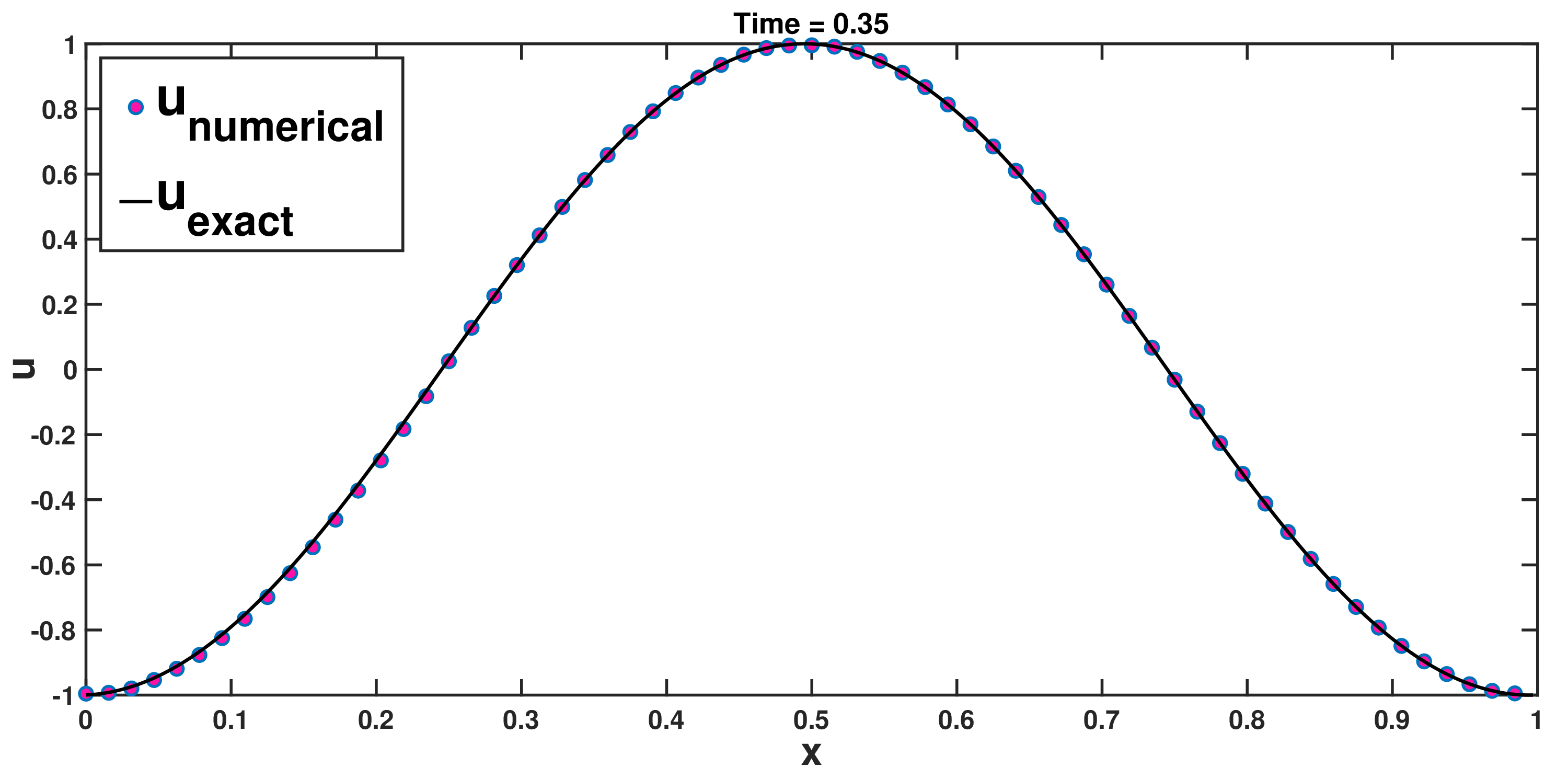}
         \caption{$u$ and $\varepsilon = 10^{-8}$}
         \label{(13a)}
     \end{subfigure}
     \hfill
     \begin{subfigure}[b]{0.48\textwidth}
         \centering
         \includegraphics[width=\textwidth]{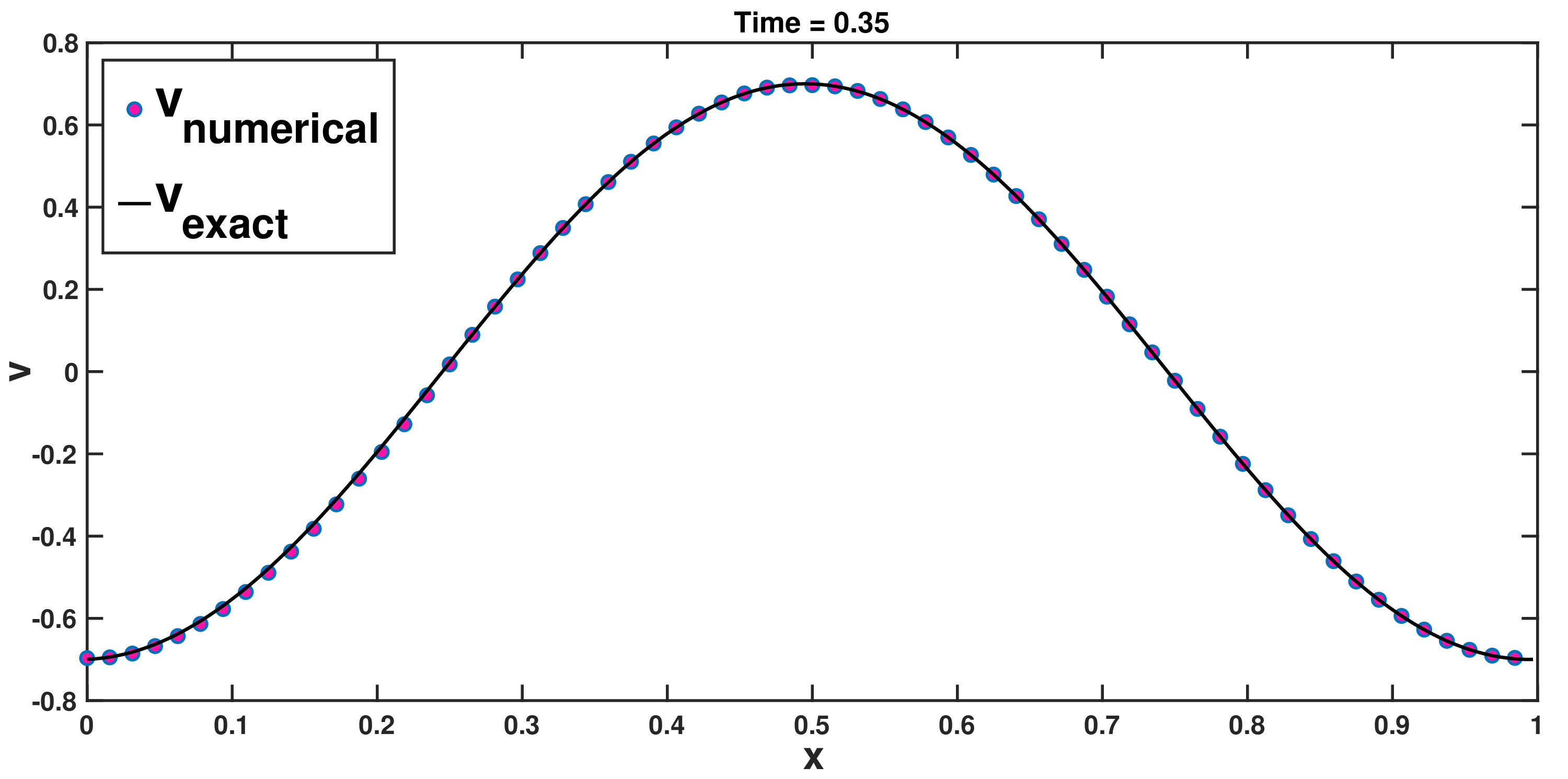}
         \caption{$v$ and $\varepsilon = 10^{-8}$}
         \label{(13b)}
     \end{subfigure}
    \medskip
    \centering
     \begin{subfigure}[b]{0.5\textwidth}
         \centering
         \includegraphics[width=\textwidth]{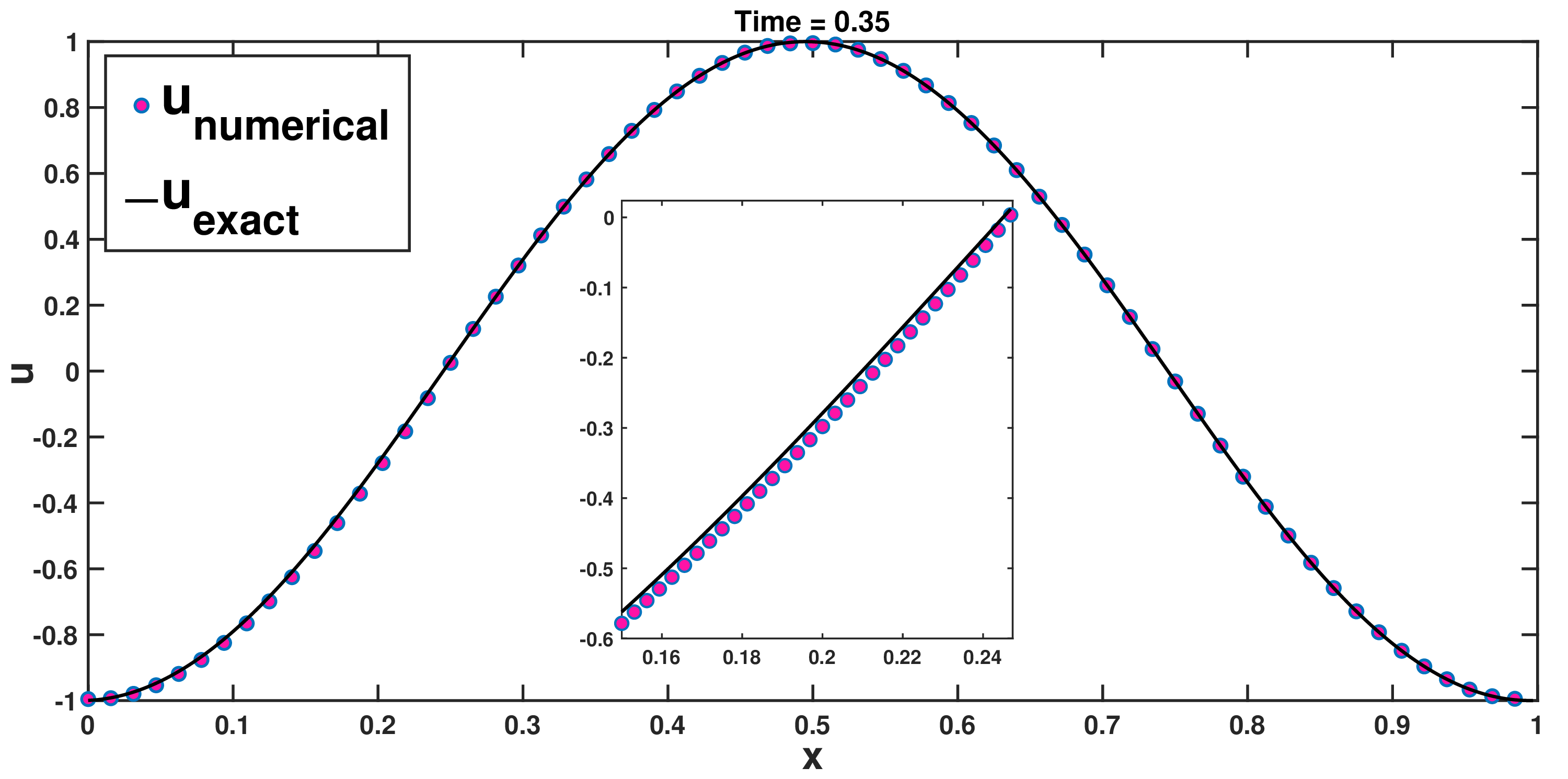}
         \caption{$u$ and $\varepsilon = 10^{-10}$}
         \label{(14a)}
     \end{subfigure}
     \hfill
     \begin{subfigure}[b]{0.48\textwidth}
         \centering
         \includegraphics[width=\textwidth]{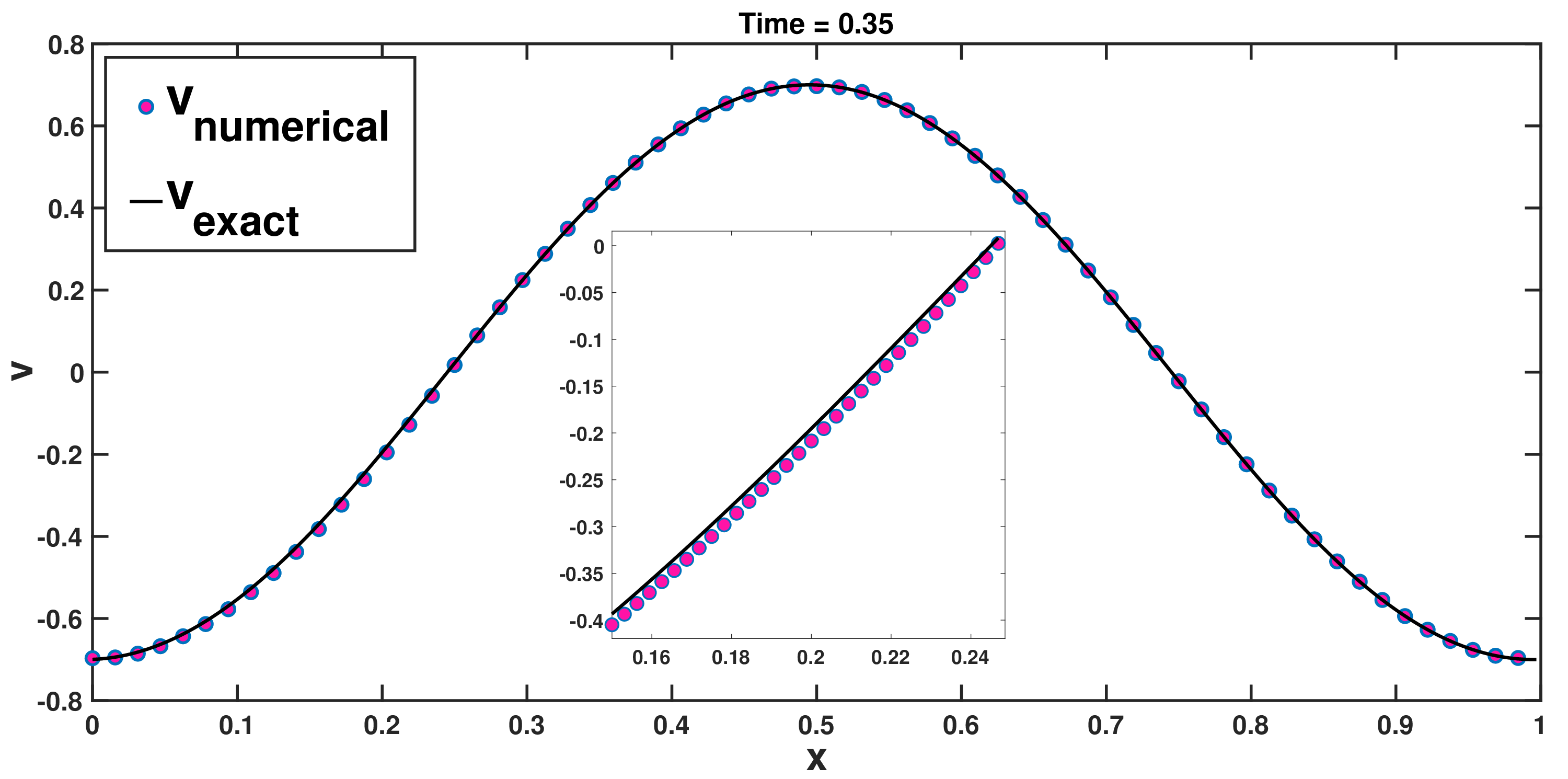}
         \caption{$v$ and $\varepsilon = 10^{-10}$}
         \label{(14b)}
     \end{subfigure}
     \caption{Jin-Xin model with smooth case: comparison between Numerical solution $u$(left) and $v$(right) and exact solution with CFL $1/3$ and $N=320$.}
        \label{fig:example1}
\end{figure}
\par Figure~\ref{fig:example1} illustrates the numerical solutions for the variables $u$ and $v$, shown on the left and right, respectively, for various values of the relaxation parameter $\varepsilon$, specifically $\varepsilon = 10^{-7},\, 10^{-8}$ and $10^{-10}$. To enhance clarity in the visual representation, a subsampling strategy is employed, plotting data at every fifth grid point\footnote{Indices $1:5:N$} instead of displaying all $N=320$ points. These results are compared against the exact solution, which is computed analytically using the Fourier transform method. The comparison offers a clear view of the accuracy of the numerical scheme under various relaxation regimes.
\begin{figure}[!ht]
     \centering
     \begin{subfigure}[b]{0.5\textwidth}
         \centering
         \includegraphics[width=\textwidth]{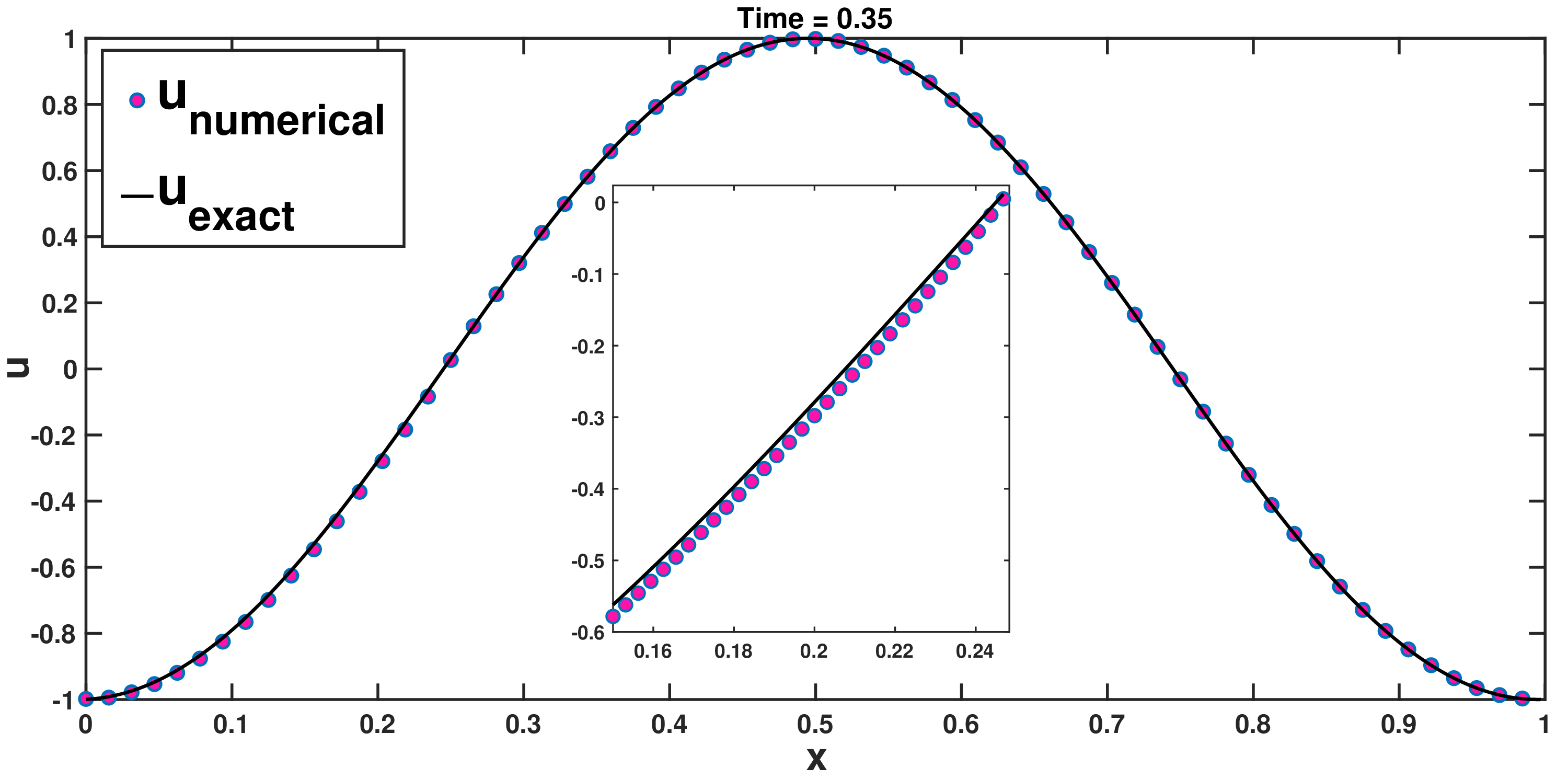}
         \caption{$u$ and $\varepsilon = 10^{-10}$}
         \label{(15a))}
     \end{subfigure}
     \hfill
     \begin{subfigure}[b]{0.48\textwidth}
         \centering
         \includegraphics[width=\textwidth]{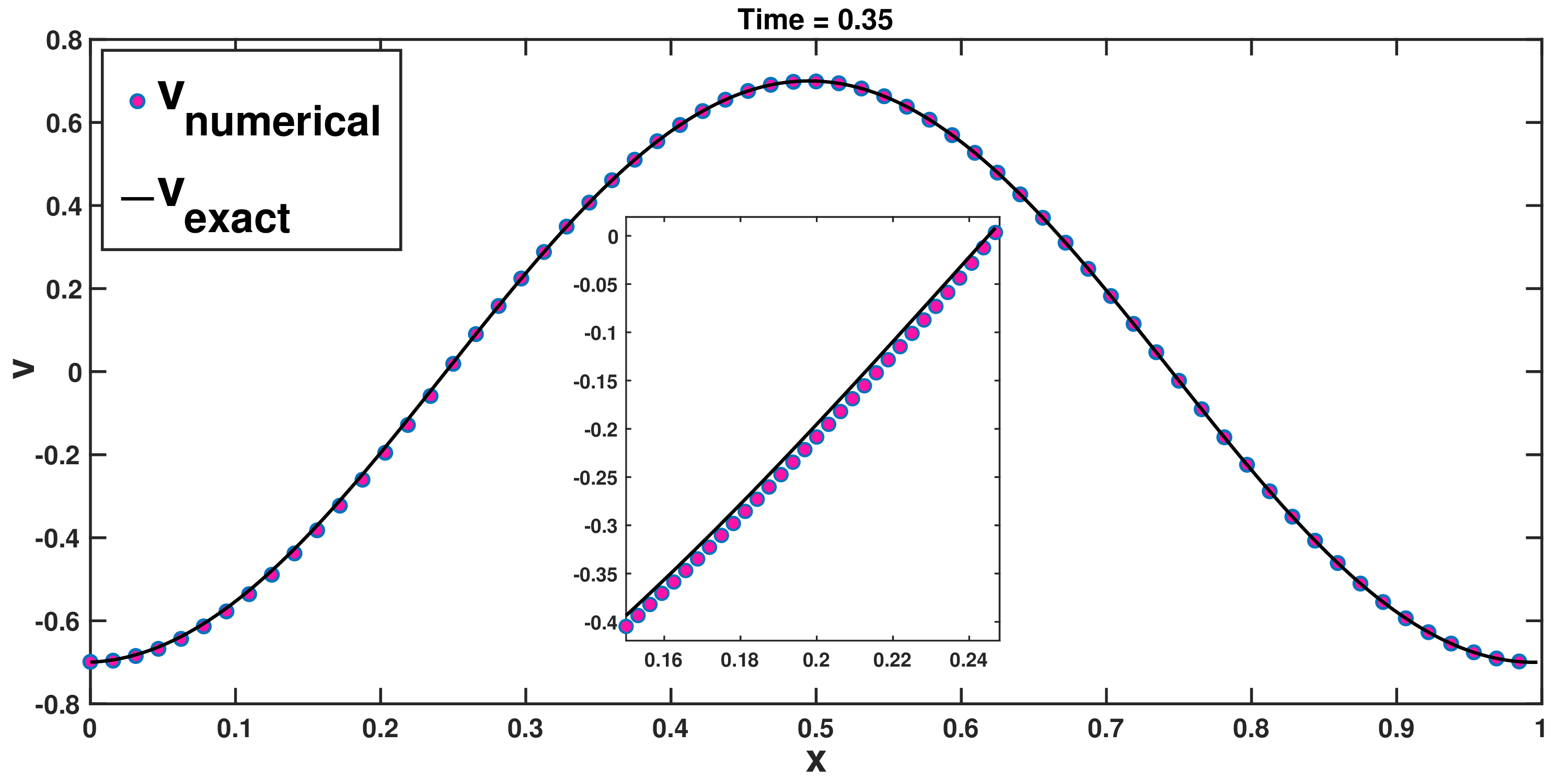}
         \caption{$v$ and $\varepsilon = 10^{-10}$}
         \label{(15b))}
     \end{subfigure}
     \caption{Jin-Xin model with smooth case: comparison between Numerical solution $u$(left) and $v$(right) and exact solution with CFL $0.9$ and $N=320$.}
        \label{fig:example1a}
\end{figure}
\par Figure~\ref{fig:example1a} displays the numerical solutions for the variables $u$ and $v$, shown on the left and right, respectively, for the relaxation parameter $\varepsilon = 10^{-10}$. The numerical results are compared with the corresponding exact solution. The computations are carried out using $N=320$ uniformly spaced grid points, with a CFL number of $0.9$ imposed to ensure numerical stability. The solution is advanced in time up to the final time $T=0.35$.
\renewcommand{\arraystretch}{1.2}{
\begin{table}[ht!]
 \caption{Error and order of convergence for Jin-Xin model with smooth initial condition \eqref{smoothdata:XinJin} for $u$}
 \centering
 \begin{tabular}{*{9}{c}}
 \toprule
\multirow{2}{*}{N} & 
\multicolumn{2}{c}{$\varepsilon=10^{-10}$} &
\multicolumn{2}{c}{$\varepsilon=10^{-8}$} & 
\multicolumn{2}{c}{$\varepsilon=10^{-7}$} \\
\cmidrule(lr){2-3}
\cmidrule(lr){4-5}
\cmidrule(lr){6-7}

& $L^{1}$ - Error & Order 
& $L^{1}$ - Error & Order 
& $L^{1}$ - Error & Order \\
\midrule
20  & 2.5785e-03 & -& 2.5785e-03 & - & 2.5785e-03& - \\

40 & 6.6550e-04 & 1.95 & 6.6550e-04 & 1.95 & 6.6550e-04 & 1.95 \\

 80 &1.6381e-04 & 2.02 & 1.6381e-04 & 2.02 &  1.6381e-04 & 2.02 \\

 160 & 4.0213e-05 &2.03 & 4.021e-05 & 2.03 & 4.021e-05 & 2.03 \\

320 & 9.9344e-06 &2.02& 9.9341e-06  & 2.02 &  9.9345e-06 & 2.02 \\

640 & 2.4653e-06 &2.01& 2.4652e-06 & 2.01 & 2.4654e-06 & 2.01\\

\bottomrule
\end{tabular}
\label{Tab3}
\end{table}
\begin{table}[ht!]
 \caption{Error and order of convergence for Jin-Xin model with smooth initial condition \eqref{smoothdata:XinJin} for $v$}
 \centering
 \begin{tabular}{*{9}{c}}
 \toprule
\multirow{2}{*}{N} & 
\multicolumn{2}{c}{$\varepsilon=10^{-10}$} &
\multicolumn{2}{c}{$\varepsilon=10^{-8}$} & 
\multicolumn{2}{c}{$\varepsilon=10^{-7}$} \\
\cmidrule(lr){2-3}
\cmidrule(lr){4-5}
\cmidrule(lr){6-7}

& $L^{1}$ - Error & Order 
& $L^{1}$ - Error & Order 
& $L^{1}$ - Error & Order \\
\midrule
20  & 1.8050e-03 & -& 1.8050e-03 & - & 1.8049e-03& -\\

40 &4.6585e-04 & 1.95 & 4.6583e-04 & 1.95 & 4.6568e-04 & 1.95  \\

 80 &1.1466e-04 & 2.02&1.1453e-04& 2.02& 1.1341e-04 &2.04 \\

 160 &2.8148e-05  &2.03 &2.8078e-05 &  2.03& 2.7846e-05 & 2.03 \\

320 & 6.9537e-06  &2.02 & 6.9218e-06& 2.02 & 7.0057e-06 & 1.99 \\

640 & 1.7255e-06 & 2.01& 1.7229e-06& 2.01 & 1.7815e-06 & 1.98 \\

\bottomrule
\end{tabular}
\label{Tab4}
\end{table}
}

\begin{table}[ht!]
 \caption{Error and order of convergence for Jin-Xin model with smooth initial condition \eqref{smoothdata:XinJin} for $u$ and $v$}
 \centering
 \begin{tabular}{*{9}{c}}
 \toprule
\multirow{3}{*}{N} & 
\multicolumn{4}{c}{$\varepsilon=1$}\\
\cmidrule(lr){2-5}
\multicolumn{2}{c}{$u$} &
\multicolumn{2}{c}{$v$} \\
\cmidrule(lr){2-3}
\cmidrule(lr){4-5}
& $L^{1}$ - Error & Order 
& $L^{1}$ - Error & Order\\
\midrule
10  & 1.0615e-02 & -& 1.3595e-02 & - \\

20  & 2.9940e-03 & 1.83 & 4.9350e-03 & 1.46 \\

40 & 1.0307e-03 & 1.54 & 2.2526e-03 & 1.13 \\

 80 &4.3360e-04 & 1.25 & 1.1945e-03 & 0.91  \\

 160 & 2.1277e-04 & 1.02 & 6.1986e-04 & 0.95 \\

320 & 1.0999e-04 & 0.95 & 3.1625e-04  & 0.97 \\

640 & 5.6706e-05 &0.96& 1.5978e-04 & 0.98 \\

1280 & 2.8850e-05 &0.97& 8.0317e-05 & 0.99 \\
\bottomrule
\end{tabular}
\label{Tabep1}
\end{table}
To investigate the effect of stiffness, we consider three different values of $\varepsilon$, namely $\varepsilon = 10^{-7}$, $\varepsilon = 10^{-8}$, and $\varepsilon = 10^{-10}$. For the stiffest case, $\varepsilon = 10^{-10}$, zoomed-in portions of the numerical profiles for both $u$ and $v$ are included to examine the fine-scale behavior of the scheme in regions where small relaxation time significantly influences the solution. As observed from the plots, the scheme performs exceptionally well in the stiff regime, closely approximating the exact solution when $\varepsilon$ tends toward zero. However, as $\varepsilon$ increases and approaches $1$, we observed the numerical solution exhibits noticeable dissipation, which is consistent with the inherent dissipative characteristics of the scheme in such regimes. This trend highlights the method's accuracy and robustness in capturing the correct solution structure under stiffness-dominated scenarios. 


In addition to the graphical comparisons, the accuracy of the scheme is further validated through quantitative error analysis. Table~\ref{Tab3} demonstrates that as $\varepsilon$ tends toward zero, the numerical method achieves second-order accuracy in the variable $u$, independently of the stiff parameter. The table report the computed $L^1$ errors alongside the corresponding observed orders of convergence for a range of grid resolutions, confirming the consistency of the scheme. Similarly, Table \ref{Tab4} presents the error metrics and convergence rates for the variable $v$, showing comparable second-order accuracy. Together, these results substantiate the scheme's effectiveness in the smooth setting and verify its designed order of accuracy in the stiff relaxation limit. Table~\ref{Tabep1} reports the error and the experimental order of convergence (EOC) of the proposed scheme in the rarefied regime, i.e., $\varepsilon = 1$, for the Jin-Xin test case with smooth initial data. In this regime, the relaxation effect is weak and the defect term becomes dominant, which leads to a reduction in the observed EOC to values close to one. Therefore, for the rarefied regime, it is more appropriate to employ an explicit central scheme rather than an implicit formulation. 
\begin{figure}[!ht]
    \centering
    \begin{subfigure}[b]{0.49\textwidth}
         \includegraphics[width=\textwidth]{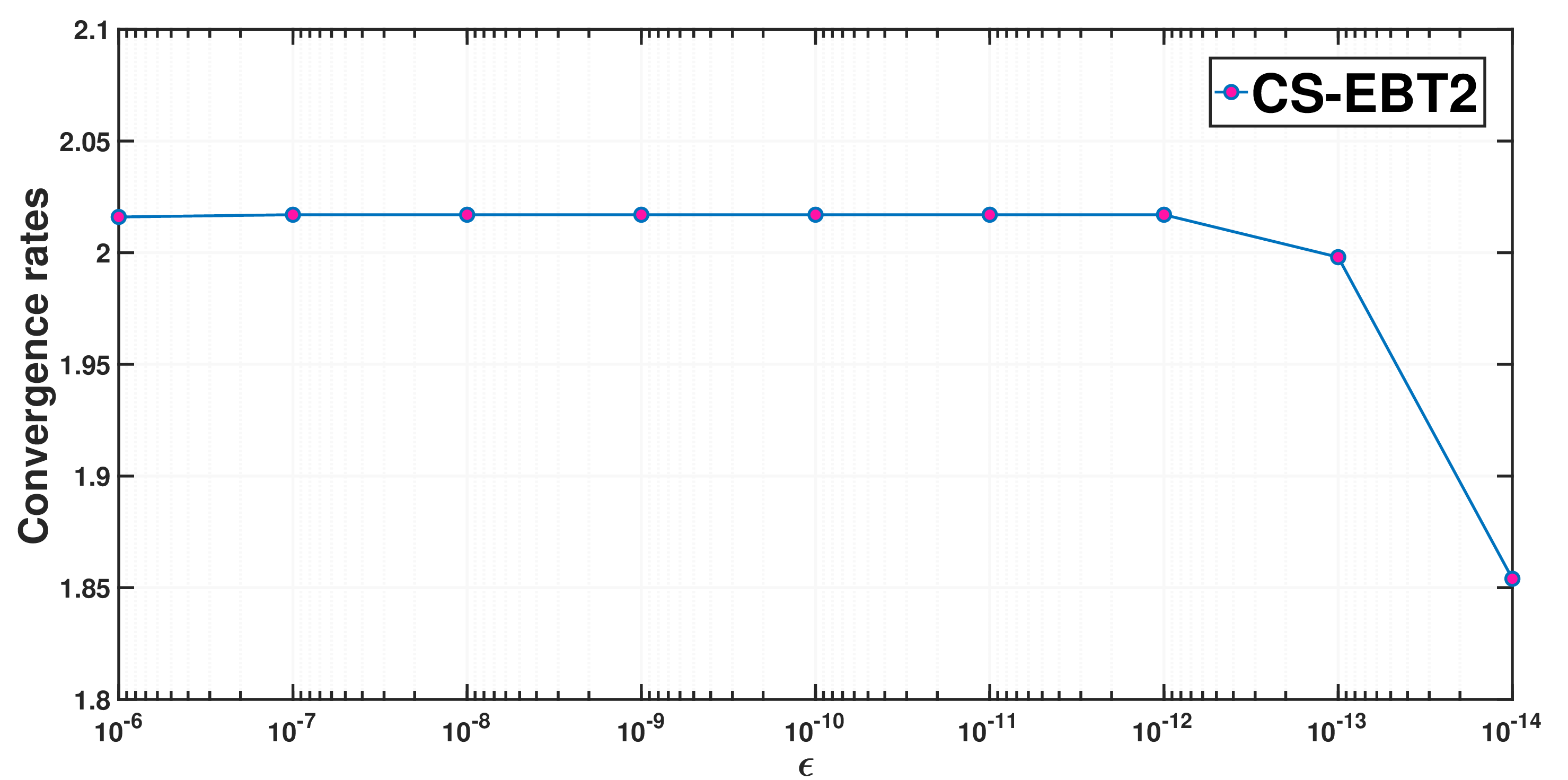} 
        \caption{Convergence rate for $u$}
        \label{fig:top}
     \end{subfigure}
     \hfill
     \begin{subfigure}[b]{0.49\textwidth}
        \centering
        \includegraphics[width=\textwidth]{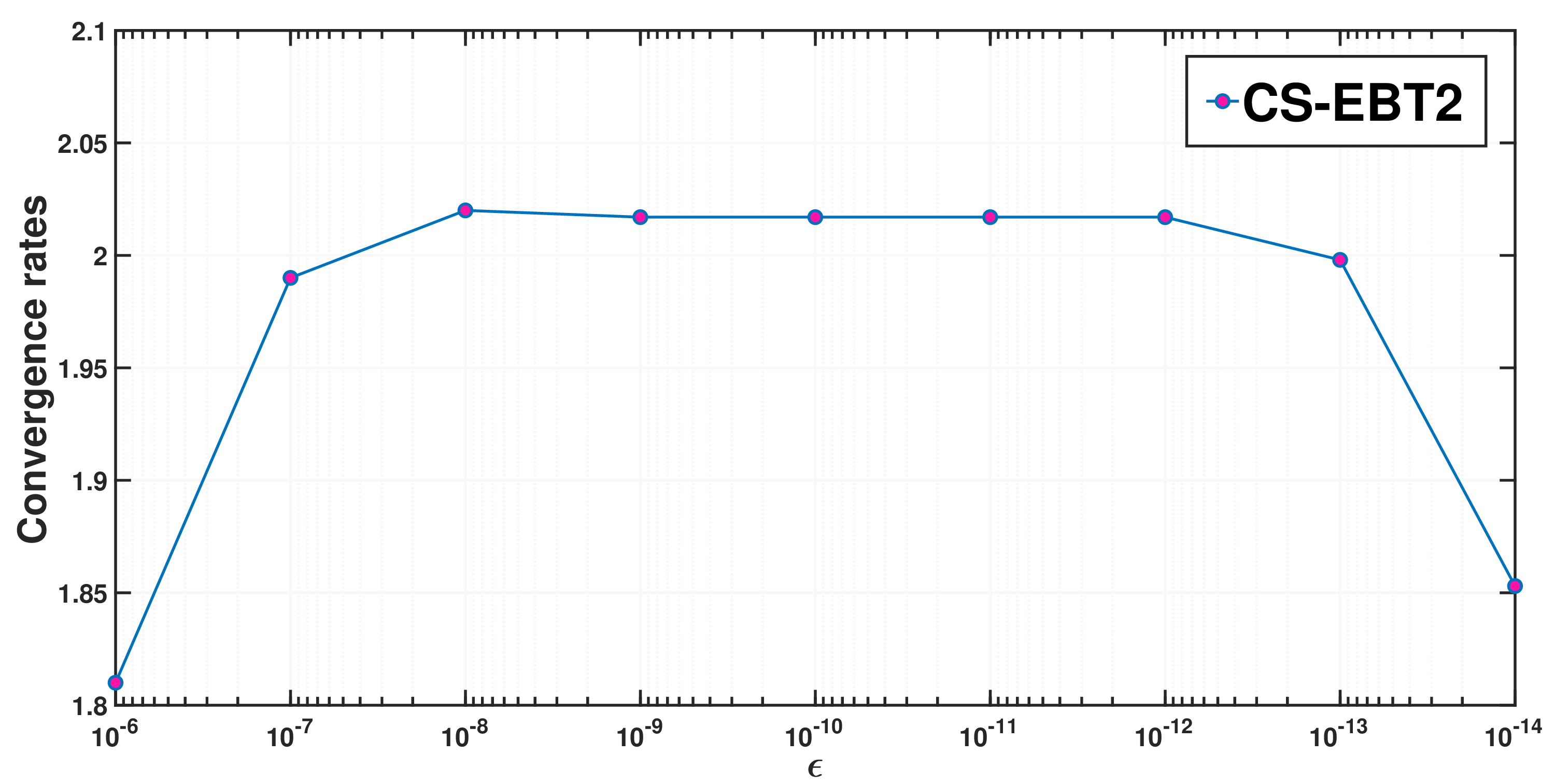} 
        \caption{Convergence rate for $v$}
        \label{fig:bottom}
     \end{subfigure}
    \caption{Convergence rates for the numerical solutions $u$ and $v$ in the smooth case of the Jin–Xin model, computed on the spatial domain $[0,1]$ up to final time $T=0.35$, using a CFL number of $0.9$ and relaxation parameters $\varepsilon$ ranging from $10^{-6}$ to $10^{-14}$.}
    \label{fig:orderConvergence}
\end{figure}
Figure \ref{fig:orderConvergence} depicts the observed order of convergence of the proposed scheme for the Jin-Xin model with respect to variables $u$ and $v$, respectively, across a range of relaxation parameters $\varepsilon = 10^{-6}$ to $10^{-14}$. As $\varepsilon$ approaches zero, corresponding to the stiff regime, the scheme consistently demonstrates second-order accuracy for both variables. These graphical findings align with the corresponding error and convergence order values reported in the table, thereby providing visual confirmation of the scheme’s second order performance.

\noindent{\textbf{Smooth case (Unprepared):}
We consider an additional test case for the Jin-Xin model with $a=0.7$ involving a smooth initial condition that is unprepared. The simulation is performed on the spatial domain $[0,1]$, with periodic boundary conditions. The initial data is explicitly defined as
\begin{eqnarray}
 \begin{cases}    
 \begin{aligned}
 &u(x,0) = \sin(2\pi x),\\
&v(x,0) = 0.1 \; u(x,0).
\end{aligned}
\end{cases}
\end{eqnarray}
 This test case demonstrates that, even when the initial condition lacks the well-preparedness property, the proposed numerical scheme performs effectively, due to its $L$-stablity property. The computations are carried out using $N=320$ uniformly spaced grid points, and a CFL number of $0.9$ is enforced for $\varepsilon = 10^{-10}$. The numerical solution is advanced in time up to a final time of  $T=0.35$.
 \begin{figure}[!ht]
     \centering
     \begin{subfigure}[b]{0.5\textwidth}
         \centering
         \includegraphics[width=\textwidth]{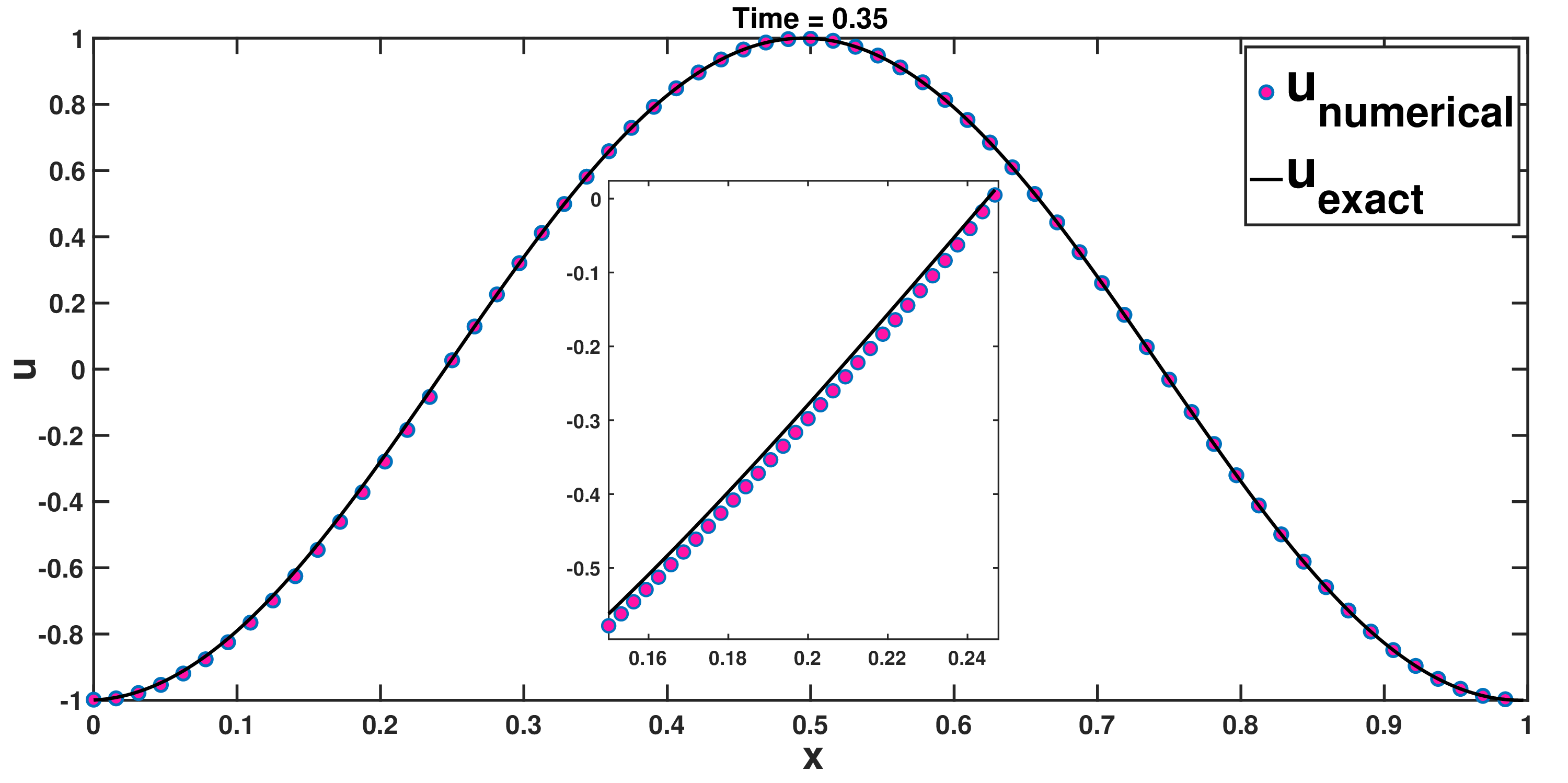}
         \caption{$u$ and $\varepsilon = 10^{-10}$}
         \label{(16a))}
     \end{subfigure}
     \hfill
     \begin{subfigure}[b]{0.48\textwidth}
         \centering
         \includegraphics[width=\textwidth]{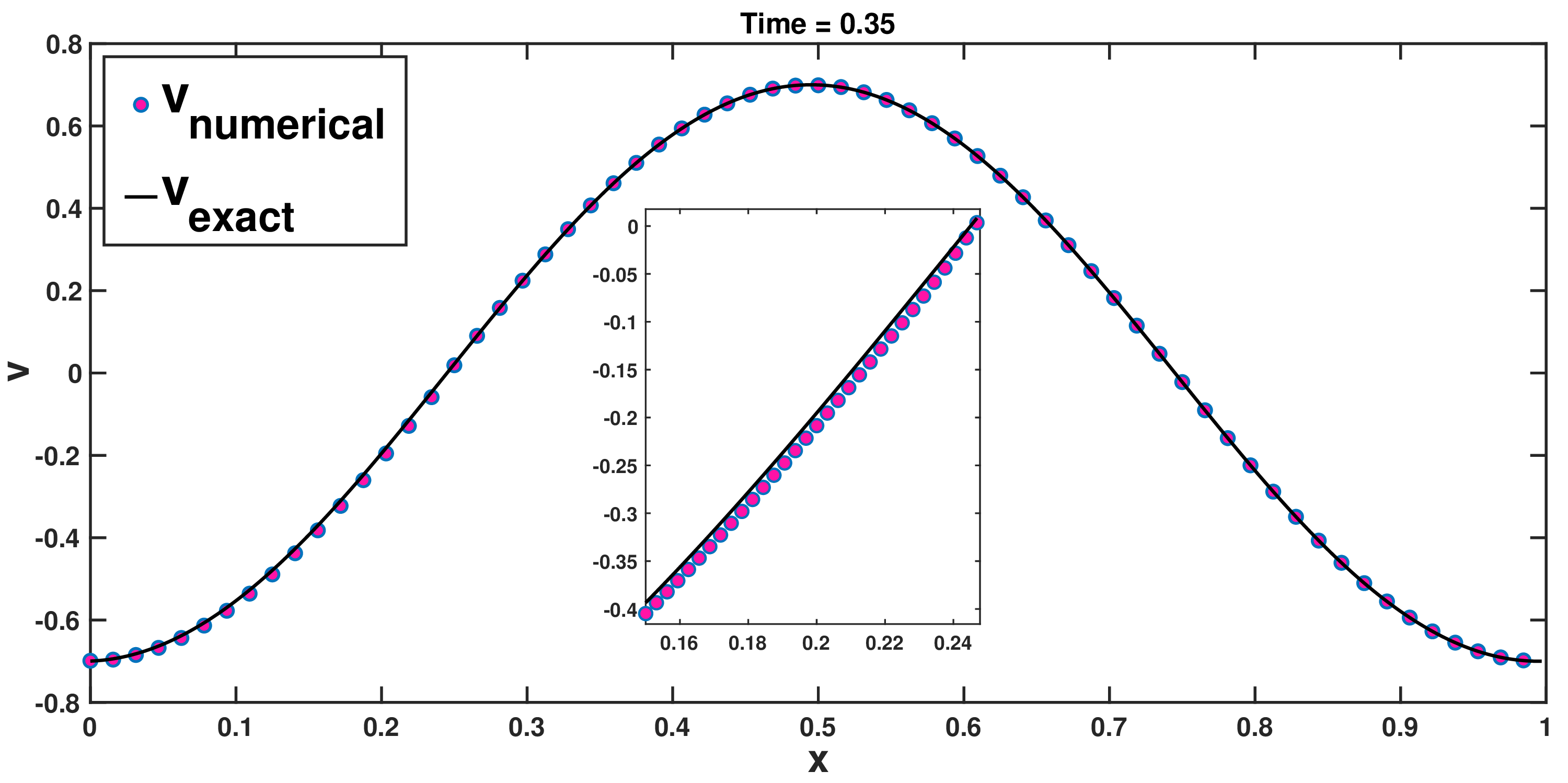}
         \caption{$v$ and $\varepsilon = 10^{-10}$}
         \label{(16b))}
     \end{subfigure}
     \caption{Jin-Xin model with an unprepared  smooth case: comparison between Numerical solution $u$(left) and $v$(right) and exact solution with CFL $0.9$ and $N=320$.}
        \label{fig:example1b}
\end{figure}
As illustrated in Figure \ref{fig:example1b}, the scheme exhibits satisfactory performance. For improved clarity, a zoomed-in view of a selected region is also provided. In the figure, the left panel corresponds to the variable $u$, while the right panel shows the variable $v$, both for $\varepsilon = 10^{-10}$.\\



\begin{figure}[!ht]
    \centering
     \begin{subfigure}[b]{0.48\textwidth}
         \centering
         \includegraphics[width=\textwidth]{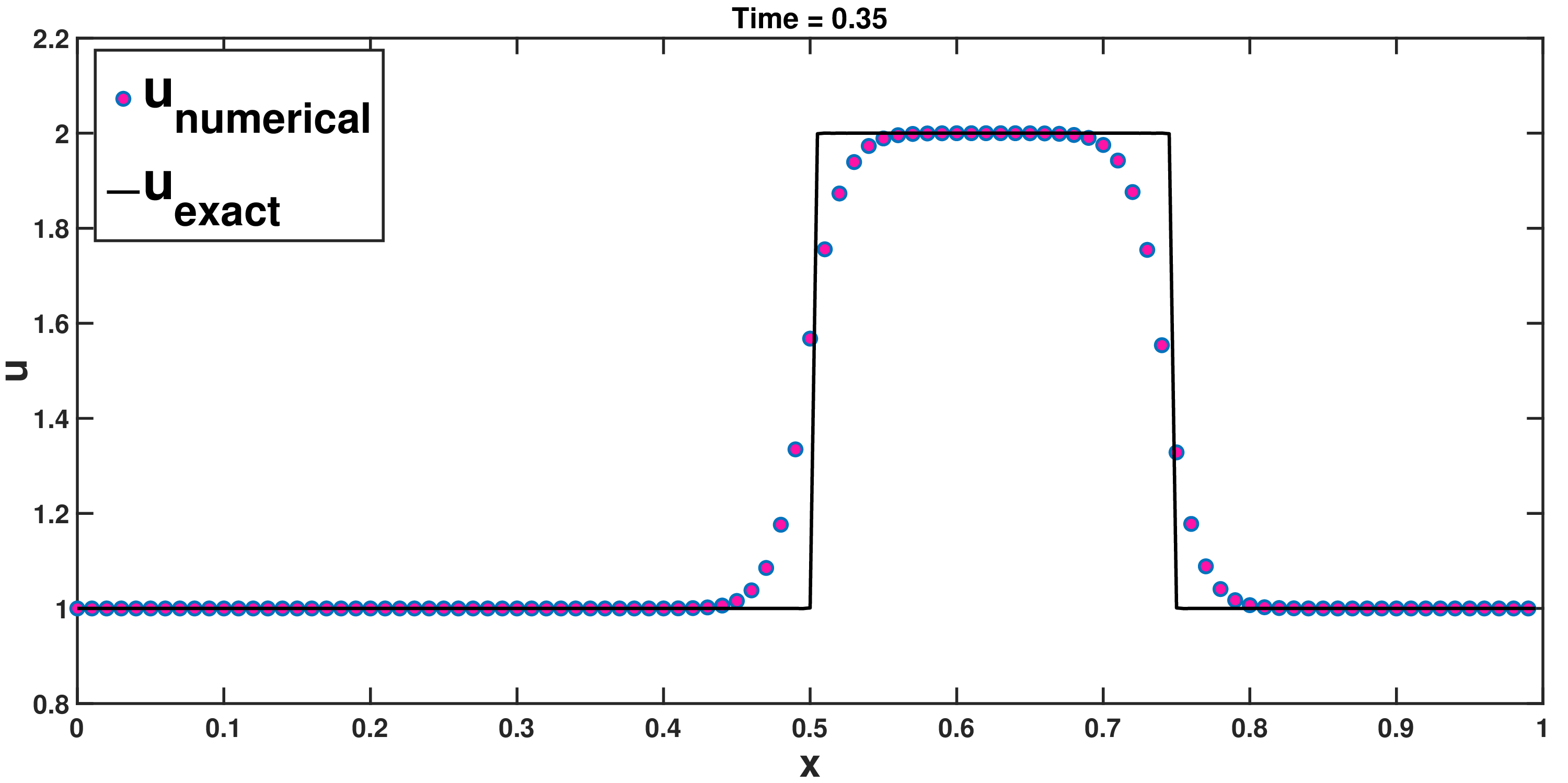}
         \caption{$u$ and $\varepsilon = 10^{-7}$}
         \label{(22a))}
     \end{subfigure}
     \hfill
     \begin{subfigure}[b]{0.48\textwidth}
         \centering
         \includegraphics[width=\textwidth]{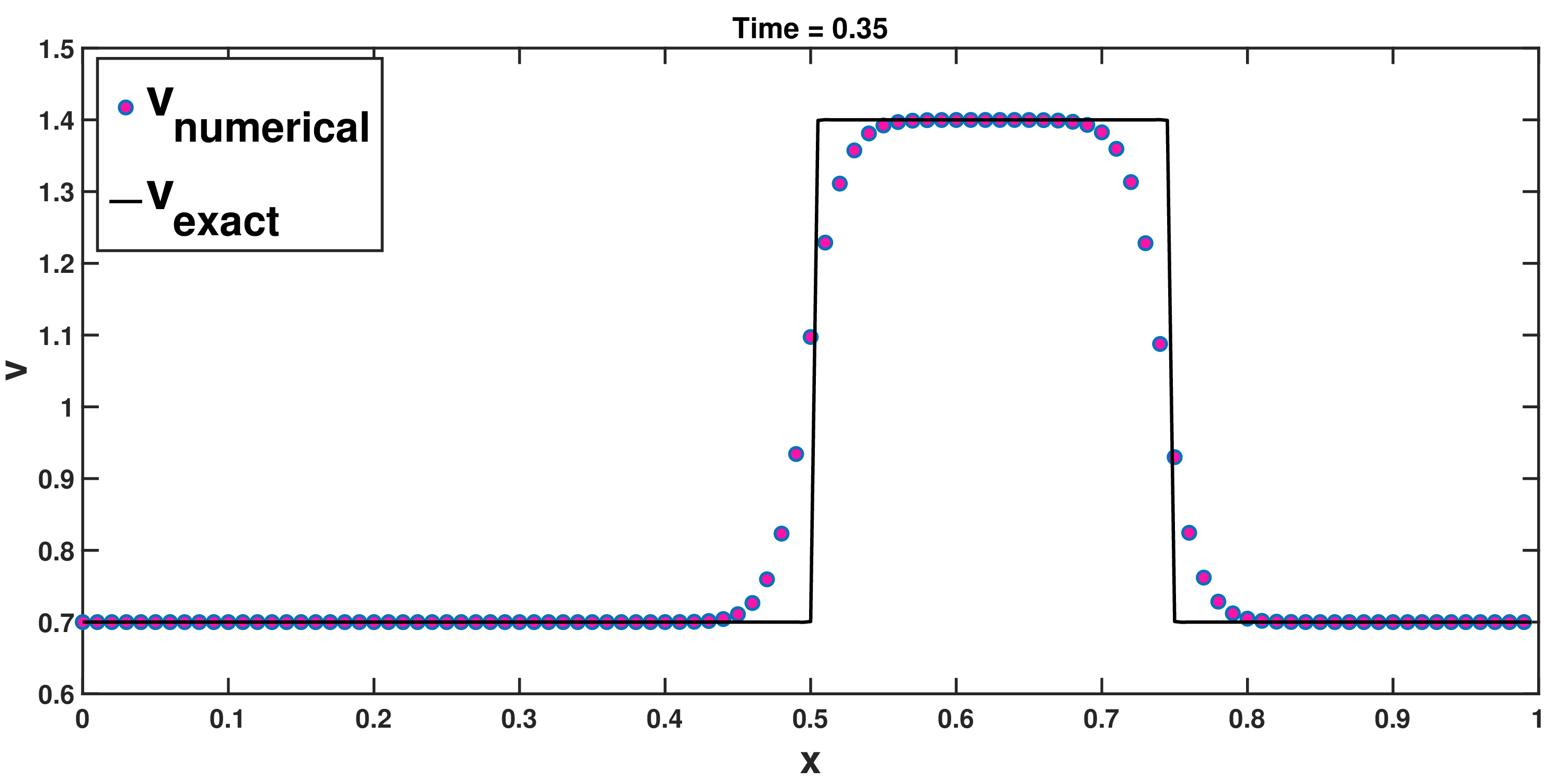}
         \caption{$v$ and $\varepsilon = 10^{-7}$}
         \label{(22b))}
     \end{subfigure}
     \medskip
    \centering
     \begin{subfigure}[b]{0.48\textwidth}
         \centering
         \includegraphics[width=\textwidth]{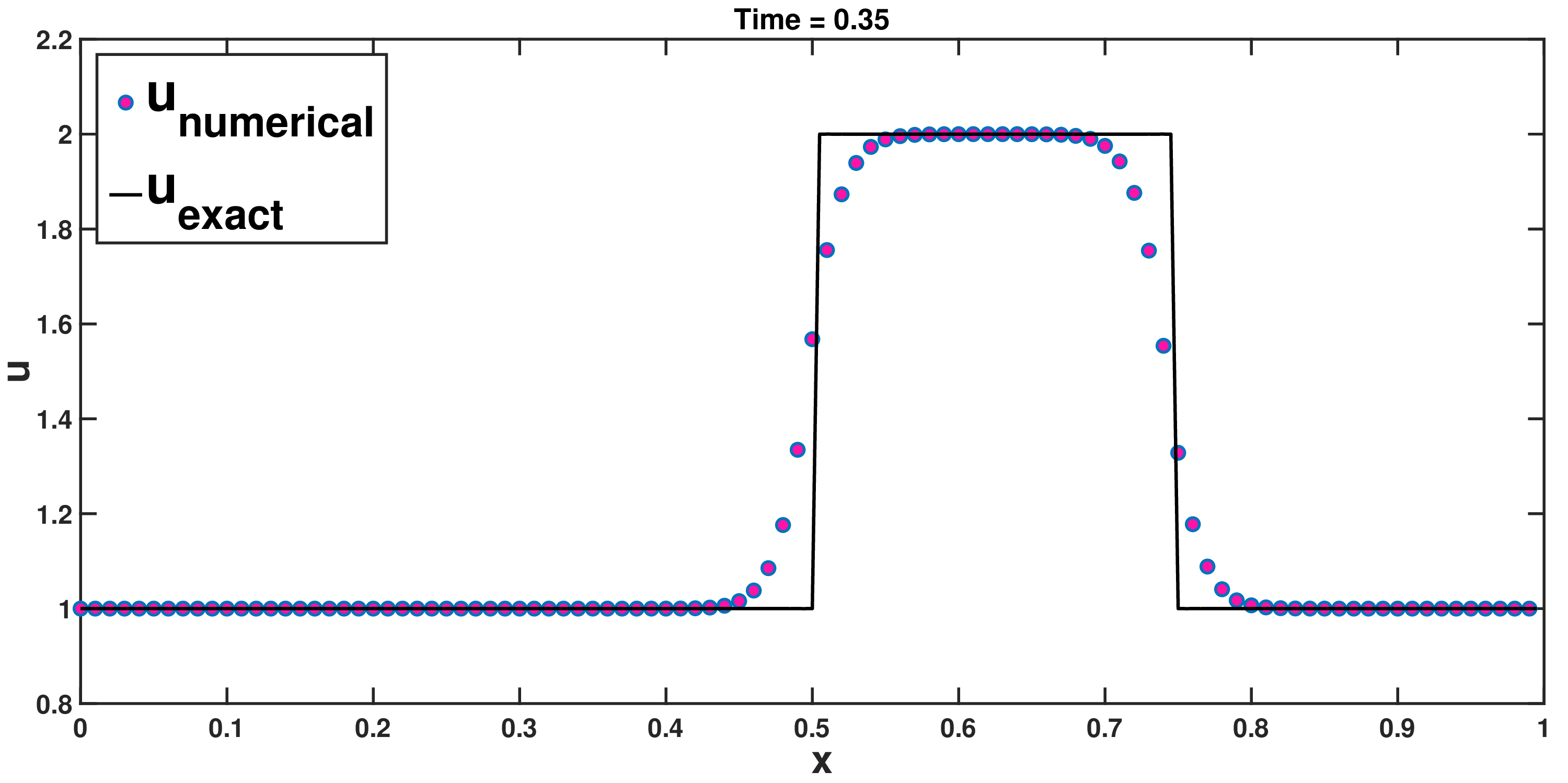}
         \caption{$u$ and $\varepsilon = 10^{-8}$}
         \label{(23a)}
     \end{subfigure}
     \hfill
     \begin{subfigure}[b]{0.48\textwidth}
         \centering
         \includegraphics[width=\textwidth]{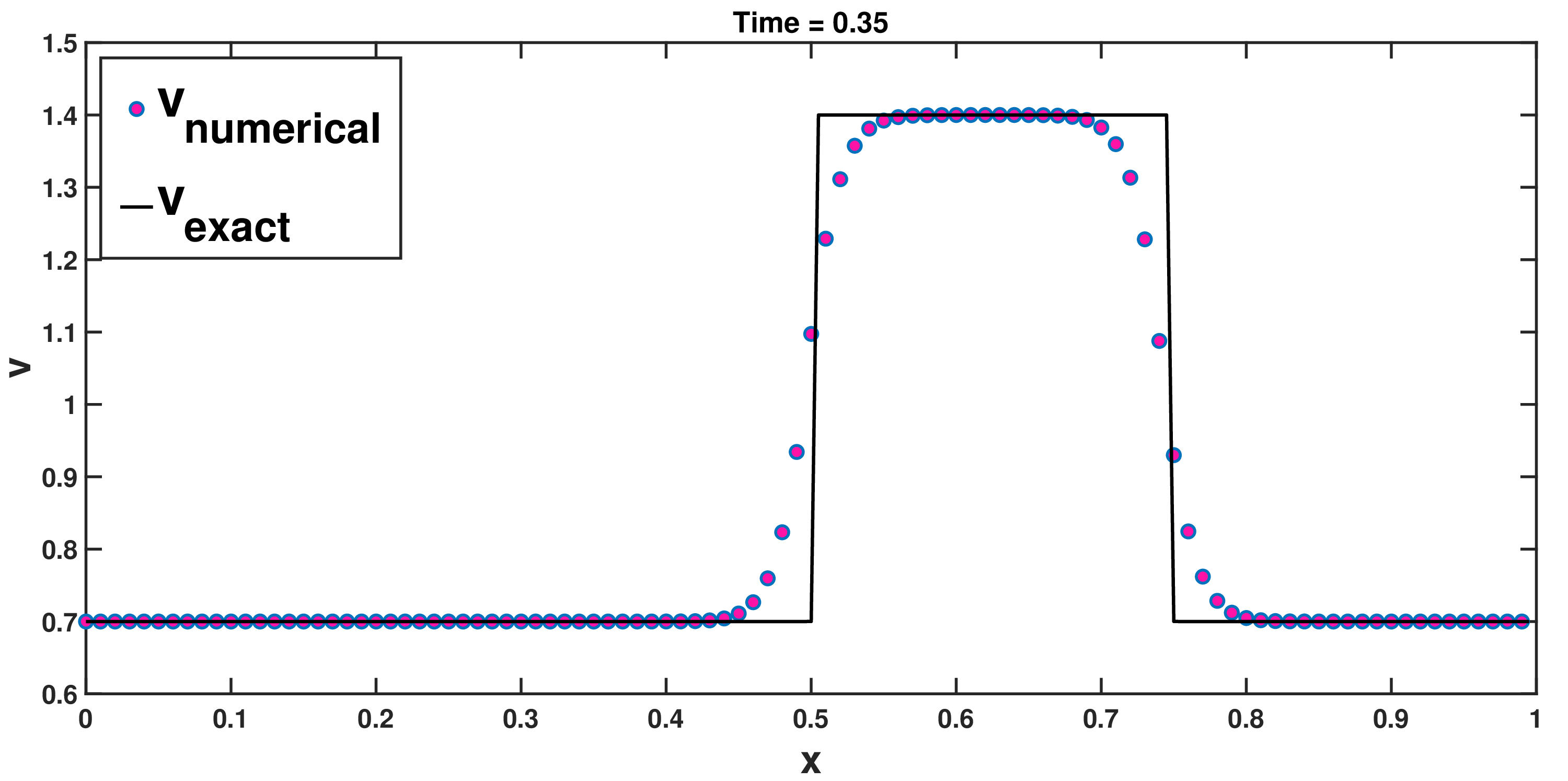}
         \caption{$v$ and $\varepsilon = 10^{-8}$}
         \label{(23b)}
     \end{subfigure}
     \medskip
    \centering
     \begin{subfigure}[b]{0.48\textwidth}
         \centering
         \includegraphics[width=\textwidth]{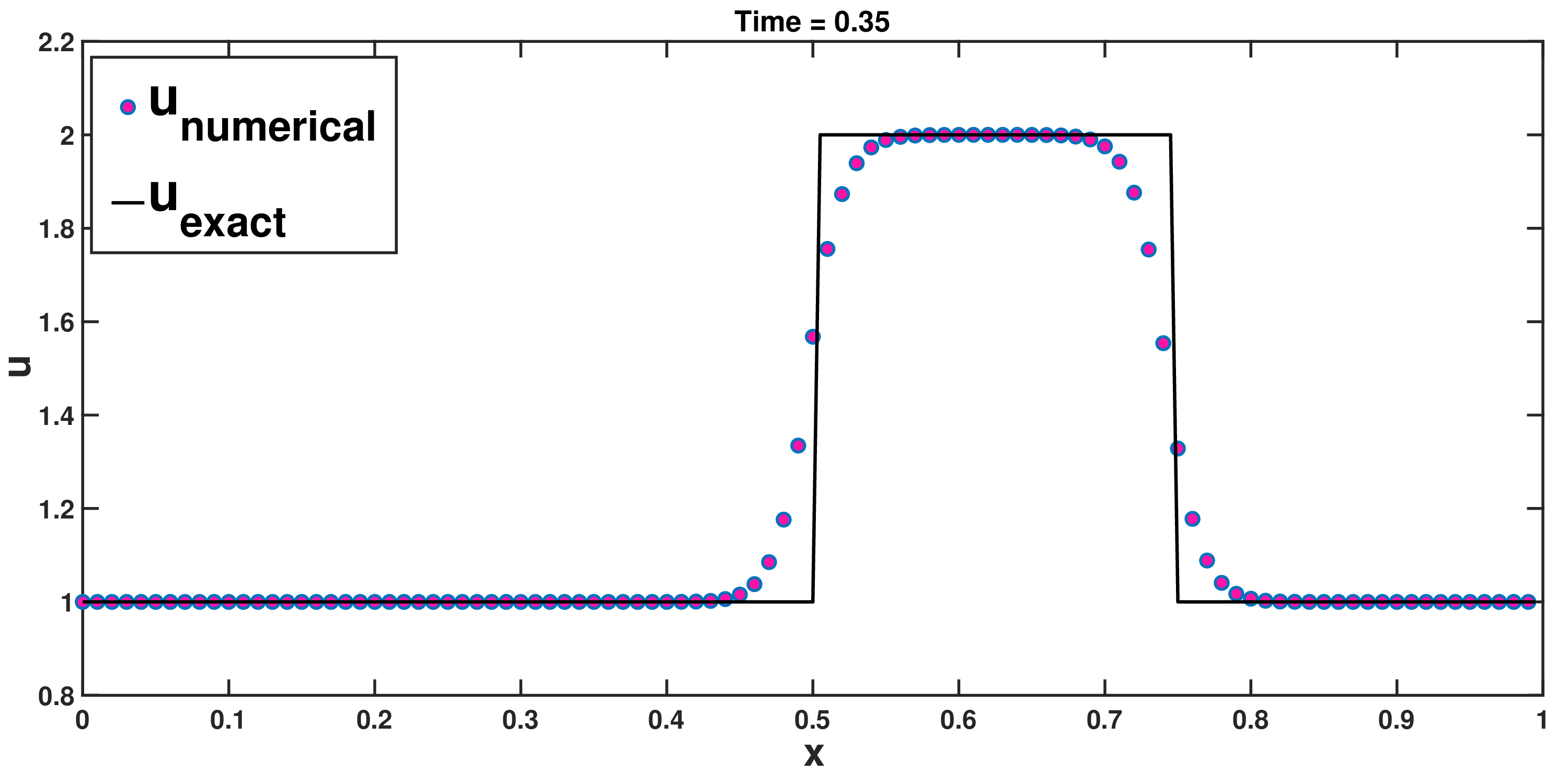}
         \caption{$u$ and $\varepsilon = 10^{-10}$}
         \label{(24a)}
     \end{subfigure}
     \hfill
     \begin{subfigure}[b]{0.48\textwidth}
         \centering
         \includegraphics[width=\textwidth]{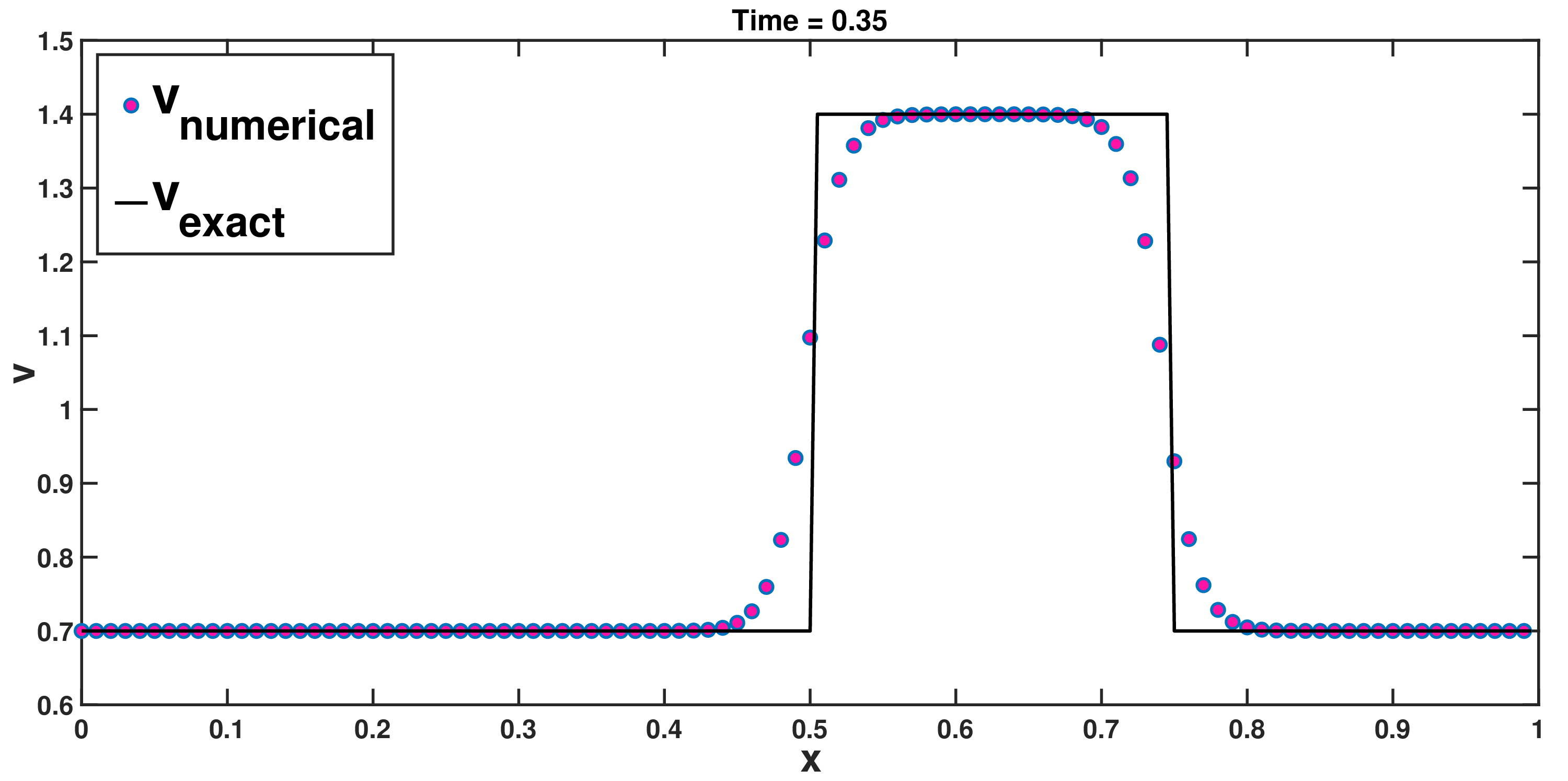}
         \caption{$v$ and $\varepsilon = 10^{-10}$}
         \label{(24b)}
    \end{subfigure}
     \caption{Jin-Xin model with non-smooth case: comparison between Numerical solution $u$(left) and $v$(right) and exact solution with CFL $1/3$ and $N=200$.}
        \label{fig:example2}
\end{figure}
\begin{figure}[!ht]
     \centering
     \begin{subfigure}[b]{0.5\textwidth}
         \centering
         \includegraphics[width=\textwidth]{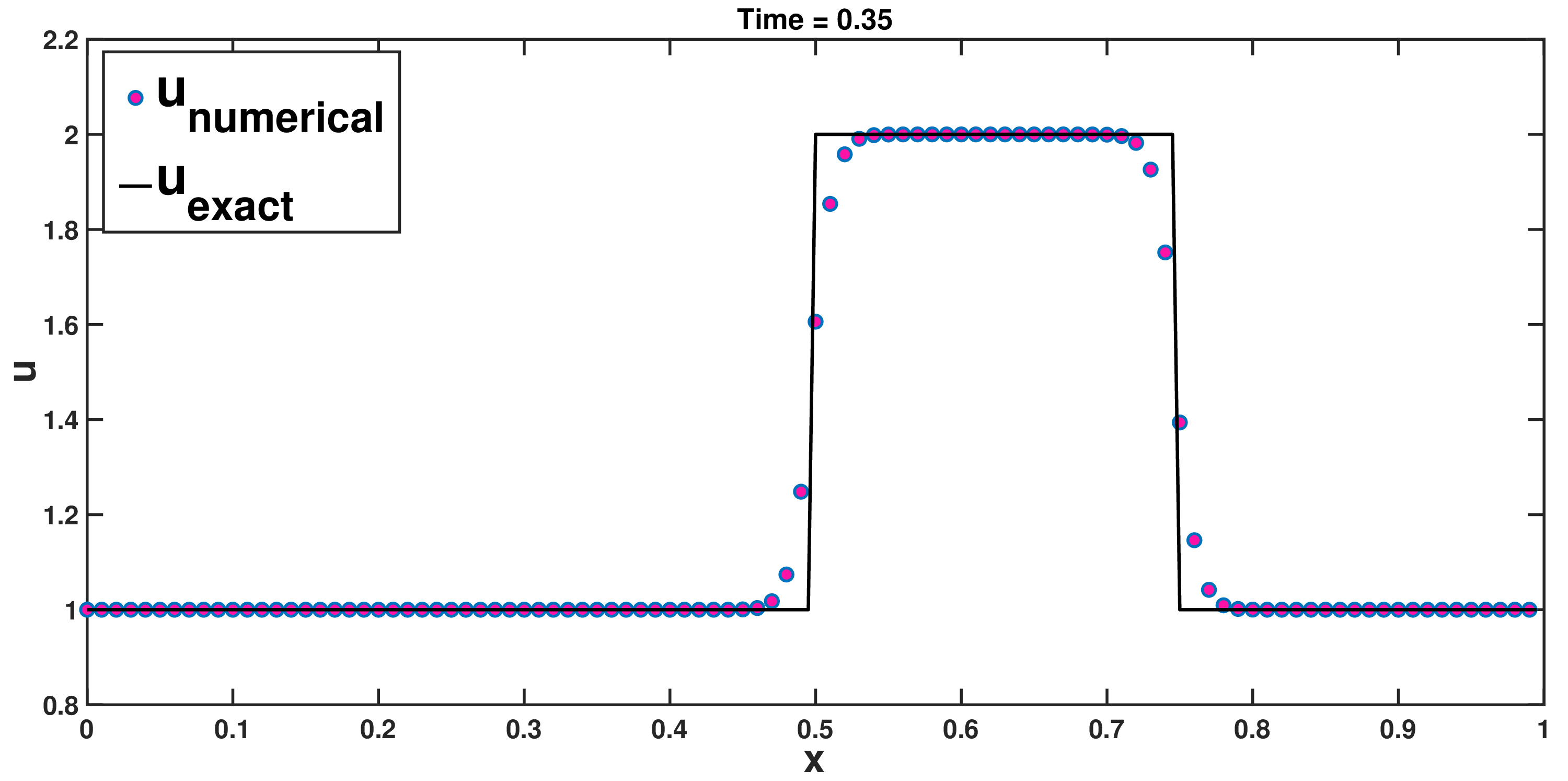}
         \caption{$u$ and $\varepsilon = 10^{-10}$}
         \label{(25a))}
     \end{subfigure}
     \hfill
     \begin{subfigure}[b]{0.48\textwidth}
         \centering
         \includegraphics[width=\textwidth]{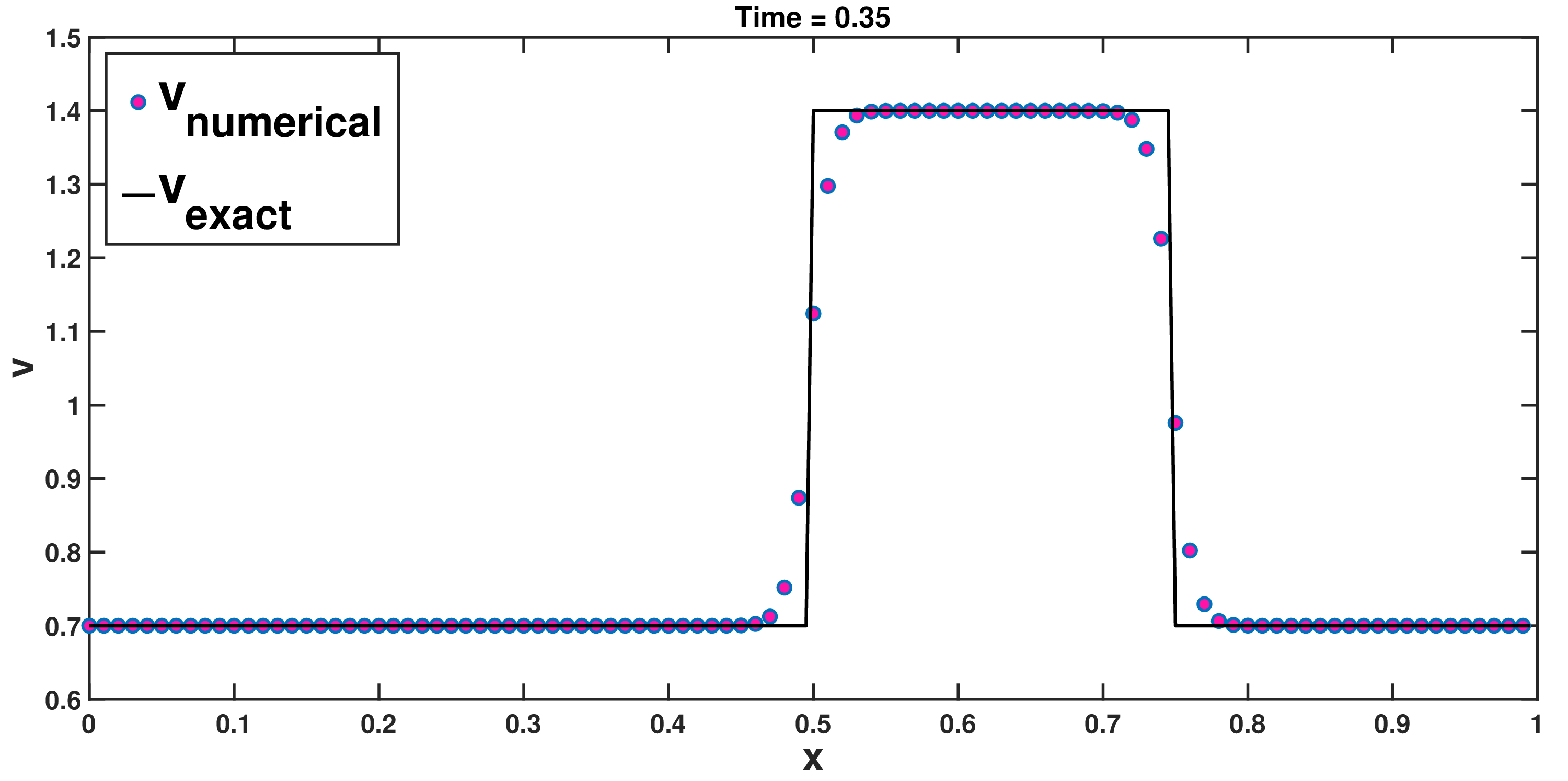}
         \caption{$v$ and $\varepsilon = 10^{-10}$}
         \label{(25b))}
     \end{subfigure}
     \caption{Jin-Xin model with non-smooth case: comparison between Numerical solution $u$(left) and $v$(right) and exact solution with CFL $0.9$ and $N=200$.}
        \label{fig:example2a}
\end{figure}
\noindent{\textbf{Non-Smooth case:}}
To further examine how the proposed CS-EBT2 scheme performs in the presence of discontinuities within the Jin–Xin model, we investigate a non-smooth unprepared initial condition. This case helps us to understand how the scheme behaves when faced with sharp gradients or discontinuities, which are common in many physical phenomena, especially in stiff regimes. Specifically, the non-smooth initial profile is defined as follows:
\begin{eqnarray}
\begin{aligned}
&u(x,0)=
\begin{cases}
 2, \, \text{if}\,\,\,\, 0.25<x<0.5,\\
1, \, \text{otherwise},
\end{cases} \\
&v(x,0) = au(x,0),
\end{aligned}
\end{eqnarray}
with $a=0.7$.  This setup creates a sharp discontinuity in the variable $u$ between the intervals  $0.25 < x < 0.5$, and the variable $v$ is defined in terms of $u$, scaled by the factor $a$. To investigate the impact of different stiffness levels on the accuracy and stability of the CS-EBT2 scheme, we consider three distinct values of the relaxation parameter $\varepsilon$, specifically $\varepsilon = 10^{-7}$, $\varepsilon = 10^{-8}$, and $\varepsilon = 10^{-10}$.
The numerical simulations are performed on the spatial domain $[0,1]$ with periodic boundary conditions, using $N=200$ uniformly spaced grid points. The final time is set to $T=0.35$, and the initial condition considered in each cases is non-smooth. The numerical solutions for the variables $u$ and $v$ are compared with the corresponding exact solutions to evaluate the performance of the proposed scheme.
Figure~\ref{fig:example2} illustrates the shock-capturing capability of the CS-EBT2 scheme under a CFL condition of $1/3$. The numerical results for $u$ and $v$, shown on the left and right panels respectively, highlight the scheme’s effectiveness in resolving discontinuities with high accuracy.

Figure~\ref{fig:example2a} presents the numerical solutions for $u$ and $v$(left and right panels, respectively) for the relaxation parameter $\varepsilon = 10^{-10}$, using a CFL number of $0.9$. The close agreement between the numerical and exact solutions further demonstrates the stability and robustness of the method, even in the presence of sharp gradients and non-smooth initial data.
As the stiffness parameter $\varepsilon$ approaches zero, the scheme continues to resolve discontinuities sharply and remain stable, faithfully preserving the solution’s steep gradients. This behavior underscores its capacity to accommodate stiff source terms without resorting to excessive numerical dissipation. Overall, the results attest to the method’s robustness in the face of both discontinuous initial data and stiff dynamics, guaranteeing that the computed solution stays physically consistent and stable throughout the simulation.

\subsection{Shallow Water Model}
We begin our analysis with the shallow water model \eqref{ShallowWater_model}, for which a smooth initial condition is considered over the spatial domain $[0,1]$, subject to periodic boundary conditions for smooth case and the spatial domain $[-1, 1]$ with transmissive boundary condition for non-smooth case,
\begin{equation*}
\begin{cases}
h_t + (hu)_x = 0,\\
(hu)_t + \left(hu^2 + \dfrac{1}{2}h^2\right)_x =
-\dfrac{1}{\varepsilon}\left(hu - \dfrac{1}{2}h^2\right).
\end{cases}
\end{equation*}

\noindent{\textbf{Smooth case (Well-prepared):}}
The initial conditions are defined as follows:
\begin{eqnarray}
 \begin{cases}
\begin{aligned}
 &h(x,0) =1+0.2 \sin(8\pi x),\\
&hu(x,0) = \dfrac{1}{2}\;h^2(x,0),
\end{aligned}
\end{cases}
\end{eqnarray}
\begin{figure}[!ht]
     \centering
     \begin{subfigure}[b]{0.5\textwidth}
         \centering
         \includegraphics[width=\textwidth]{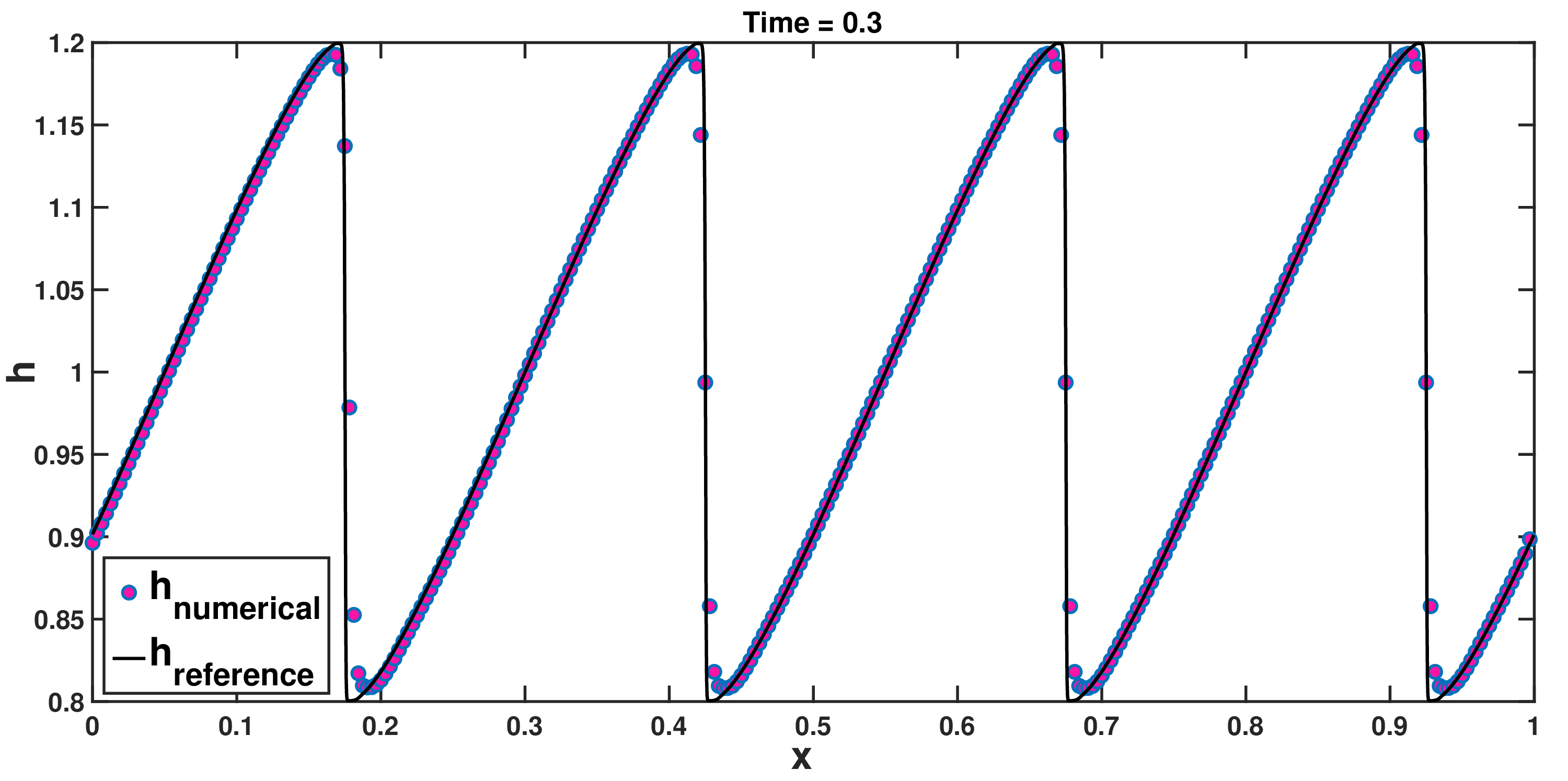}
         \caption{$h$ and $\varepsilon = 10^{-8}$}
         \label{shallowWater:1a}
     \end{subfigure}
     \hfill
     \begin{subfigure}[b]{0.48\textwidth}
         \centering
         \includegraphics[width=\textwidth]{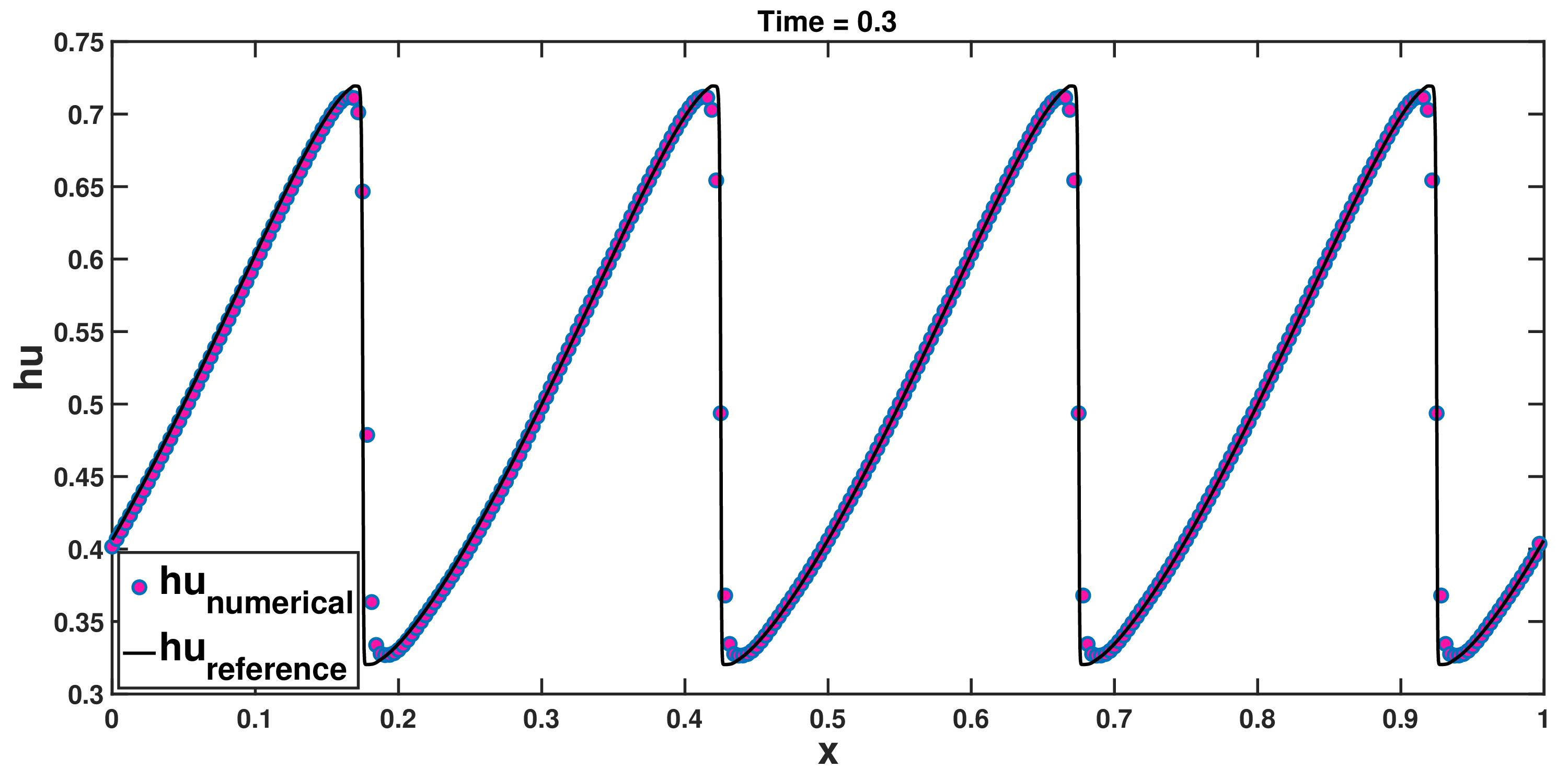}
         \caption{$hu$ and $\varepsilon = 10^{-8}$}
         \label{shallowWater:2a}
     \end{subfigure}
     \caption{Shallow water model with a smooth case: comparison between Numerical solution $h$(left) and $hu$(right) and reference solution with CFL $0.9$ and $N=320$.}
        \label{ShallowWaterfig:example1}
\end{figure}
The numerical experiment is carried out with $N=320$ grid points under periodic boundary conditions. The CFL number is set to $0.9$, and the computation is advanced up to the final time $T=0.3$. 
Figure \ref{shallowWater:1a} reports the numerical solution for the water height $h$, while Figure \ref{shallowWater:2a} displays the corresponding momentum $hu$. The obtained solution is compared with a high-resolution reference solution computed with second order implicit-explicit Runge-Kutta scheme on a much finer mesh with $N=3200$.  
The numerical results clearly show that the proposed scheme is able to capture the early stages of shock formation with high fidelity. In particular, the steepening of the wave profile is well reproduced without the introduction of spurious oscillations, and the position of the discontinuity is accurately tracked. Moreover, the scheme preserves the overall shape of the solution and provides a very good agreement with the reference even on the relatively coarse grid used in this test. 
These observations highlight both the robustness of the method in the presence of stiff relaxation and its ability to correctly resolve nonlinear wave interactions in the shallow water regime.


\noindent{\textbf{Non-smooth case (Unprepared):}}
The initial conditions are defined as follows:
\begin{eqnarray}
\begin{cases}
\begin{aligned}
&h(x,0)=
\begin{cases}
 1, \, \text{if}\,\,\,\, 0<x<0.2,\\
0.2, \, \text{otherwise},
\end{cases} \\
&hu(x,0) = -\dfrac{1}{2}\;h^2(x,0),
\end{aligned}
\end{cases}
\end{eqnarray}
where $hu$ be the non-local-equilibrium initial data.
For this case the computation is carried out with $N=320$ grid points under transmissive boundary conditions. The CFL number is 
is fixed at $0.9$ with final time $T=0.5.$ Figure \ref{shallowWater:1b} illustrates the results for the height component $h$, while Figure \ref{shallowWater:2b} shows the corresponding momentum $hu$. In both figures, the numerical solution is compared with a high-resolution reference solution obtained on a finer mesh with $N=3200.$
\begin{figure}[!ht]
     \centering
     \begin{subfigure}[b]{0.5\textwidth}
         \centering
         \includegraphics[width=\textwidth]{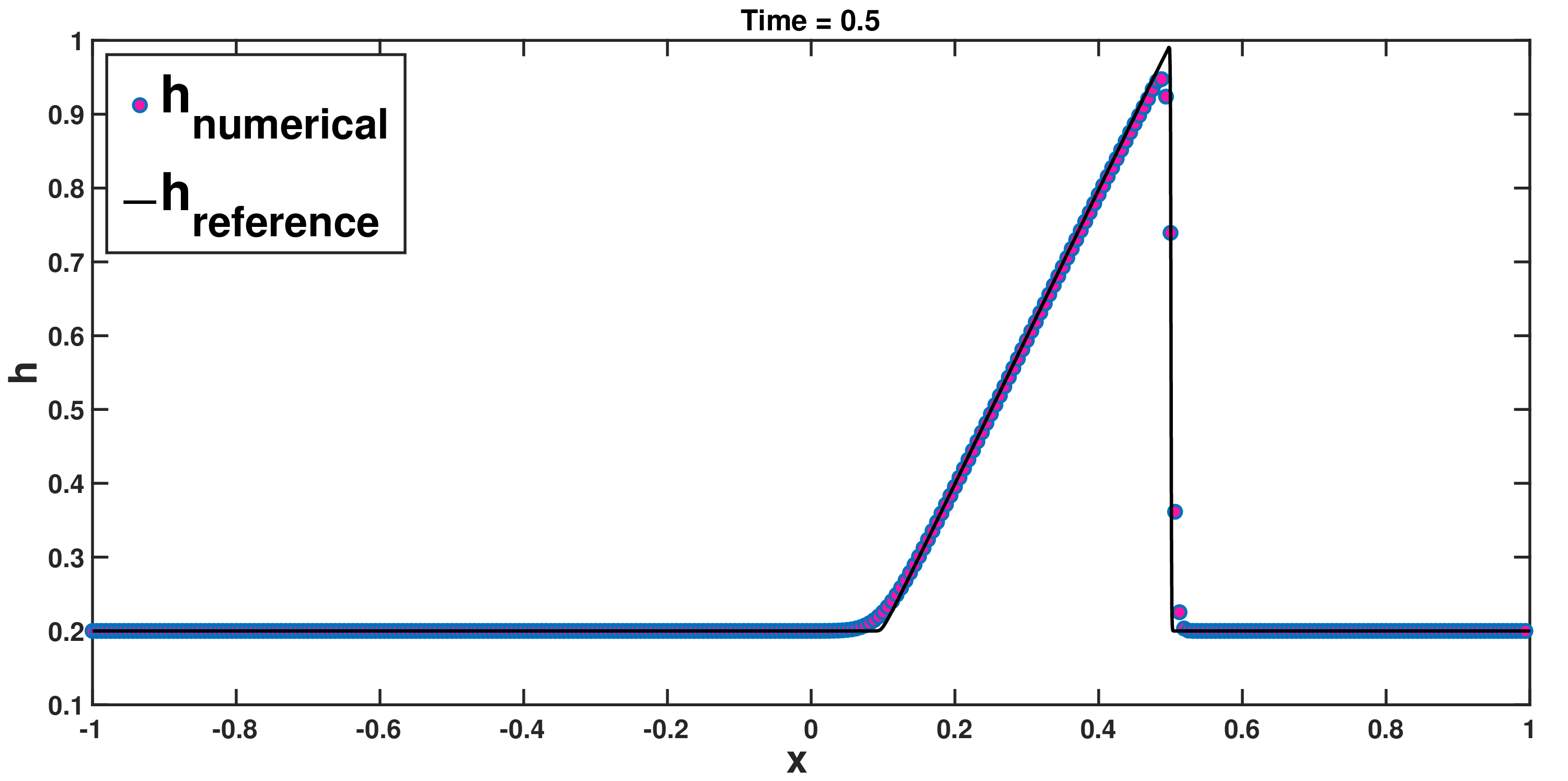}
         \caption{$h$ and $\varepsilon = 10^{-8}$}
         \label{shallowWater:1b}
     \end{subfigure}
     \hfill
     \begin{subfigure}[b]{0.48\textwidth}
         \centering
         \includegraphics[width=\textwidth]{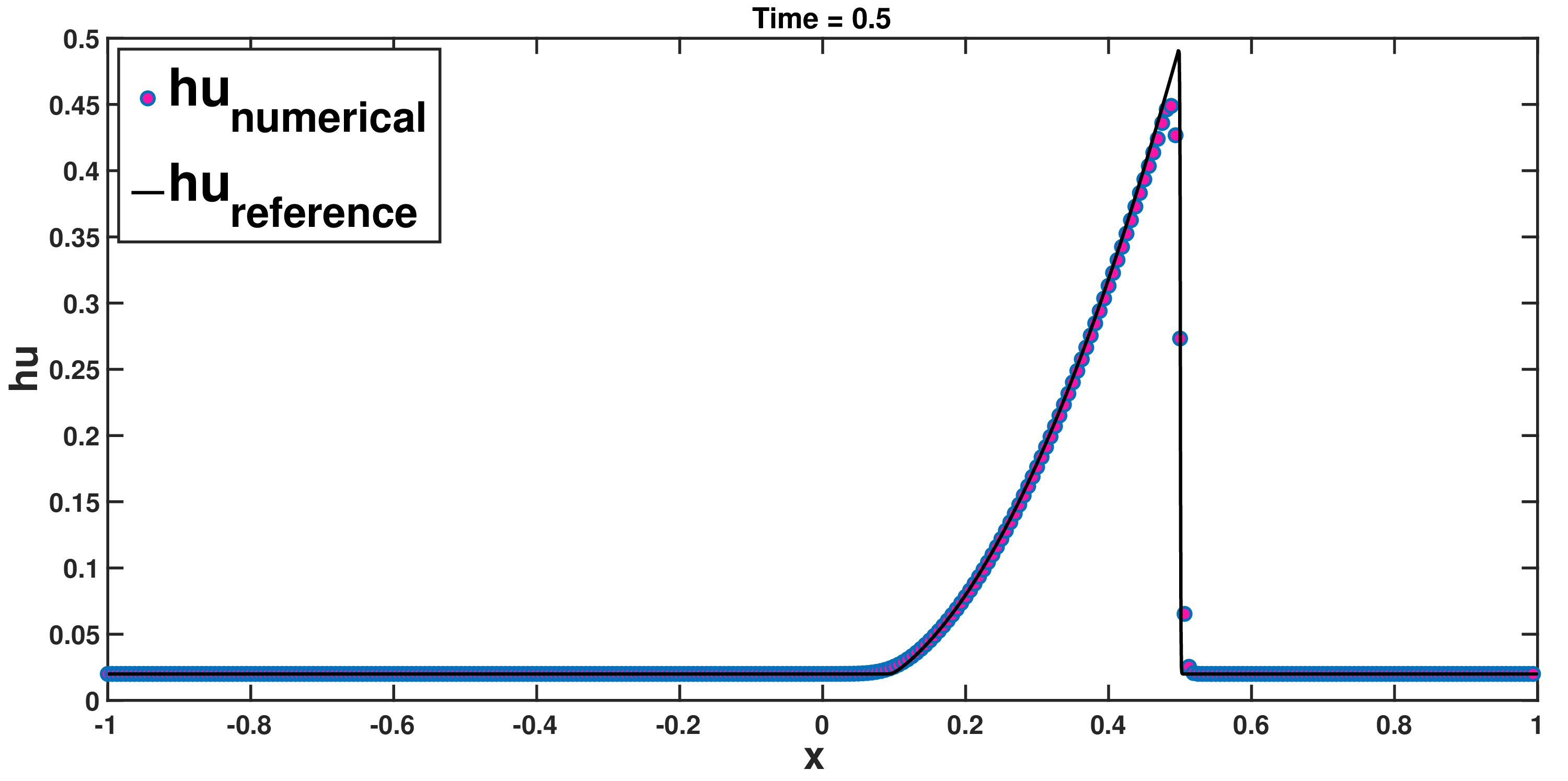}
         \caption{$hu$ and $\varepsilon = 10^{-8}$}
         \label{shallowWater:2b}
     \end{subfigure}
     \caption{Shallow water model with an non-smooth case: comparison between Numerical solution $h$(left) and $hu$(right) and reference solution with CFL $0.9$ and $N=320$.}
        \label{ShallowWaterfig:example2}
\end{figure}
We now consider the case where $hu$ represents a non-equilibrium initial state. 
The computation is performed with $N=320$ grid points under transmissive boundary conditions. 
The CFL number is set to $0.9$, and the final time is chosen as $T=0.5$. 
Figure \ref{shallowWater:1b} displays the numerical solution for the water height $h$, while Figure \ref{shallowWater:2b} reports the corresponding momentum $hu$. 
In both cases, the computed solution is compared against a reference solution obtained on a much finer mesh with $N=3200$.  

The results confirm that the scheme is able to correctly capture the relaxation dynamics induced by the non-equilibrium initial condition. 
The transition toward the equilibrium state is well resolved, with the steepening of the wave profile clearly visible in the height variable and consistently reflected in the momentum component. 
The comparison with the reference solution highlights the good accuracy of the method, even on the coarse grid, and shows that the scheme is robust in handling the interplay between nonlinear wave propagation and relaxation effects. 
Overall, the test demonstrates the effectiveness of the proposed approach in reliably reproducing the correct asymptotic behavior of the system, even in the presence of an unprepared initial condition. 
This robustness is a direct consequence of the \textit{AP-property} of the scheme, which ensures the correct coupling between the relaxation dynamics and the underlying hyperbolic structure across different regimes.

\subsection{Broadwell Model}
Following the Jin-Xin model, we extend our investigation to the Broadwell model \eqref{Broadwell_model}. For smooth initial profiles, periodic boundary conditions are imposed, whereas transmissive boundary conditions are employed for non-smooth initial data,
\begin{equation*}
\begin{cases}
\rho_t + m_x = 0,\\
m_t + z_x = 0,\\
z_t + m_x = \dfrac{1}{2\varepsilon}\left(\rho^2 + m^2 - 2\rho z\right).
\end{cases}
\end{equation*}

\noindent{\textbf{Smooth case:}}
The initial conditions are defined as follows:
\begin{eqnarray}\label{smoothdata:Broadwell}
 \begin{cases}
\begin{aligned}
 &\rho(x,0) =1+0.3 \sin(2\pi x),\\
&u(x,0) = 0.5+0.1 \sin(2\pi x), \;\; m(x,0) = \rho(x,0) u(x,0),\\
 &z(x,0) = 0.5\rho(1+u^2).
\end{aligned}
\end{cases}
\end{eqnarray}

\begin{table}[ht]
 \caption{Error and order of convergence for Broadwell model with \textcolor{red}{initial} condition \eqref{smoothdata:Broadwell} for density $\rho$}
 \centering
 \begin{tabular}{*{9}{c}}
 \toprule
\multirow{2}{*}{N} & 
\multicolumn{2}{c}{$\varepsilon=10^{-8}$} &
\multicolumn{2}{c}{$\varepsilon=0.02$} & 
\multicolumn{2}{c}{$\varepsilon=1$} \\
\cmidrule(lr){2-3}
\cmidrule(lr){4-5}
\cmidrule(lr){6-7}

& $L^{1}$ - Error & Order 
& $L^{1}$ - Error & Order 
& $L^{1}$ - Error & Order\\
\midrule

20  & 6.0363e-04 &-& 5.9336e-04 & - & 4.7706e-04 & - \\

 40 & 1.5790e-04  &1.93 & 1.6464e-04 & 1.85 & 1.2922e-04 & 1.88 \\

 80 & 4.0953e-05 &1.95 & 4.4846e-05 & 1.88 & 3.5338e-05  & 1.87\\

 160 & 1.0422e-05 &1.97 & 1.1890e-05  & 1.92 & 9.1849e-06 & 1.94\\

320 & 2.6288e-06 &1.99 & 3.0807e-06  & 1.95 & 2.3315e-06  & 1.98  \\
 
640 & 6.5991e-07 &1.99 & 7.8600e-07 & 1.97 & 5.8835e-07 & 1.99 \\

 \bottomrule
\end{tabular}
\label{Tab1b}
\end{table}

\begin{table}[ht]
 \caption{Error and order of convergence for Broadwell model with initial condition \eqref{smoothdata:Broadwell} for momentum $m$}
 \centering
 \begin{tabular}{*{9}{c}}
 \toprule
\multirow{2}{*}{N} & 
\multicolumn{2}{c}{$\varepsilon=10^{-8}$} &
\multicolumn{2}{c}{$\varepsilon=0.02$} & 
\multicolumn{2}{c}{$\varepsilon=1$} \\
\cmidrule(lr){2-3}
\cmidrule(lr){4-5}
\cmidrule(lr){6-7}

& $L^{1}$ - Error & Order 
& $L^{1}$ - Error & Order 
& $L^{1}$ - Error & Order\\
\midrule

20  & 5.6300e-04 &-& 6.8380e-04 & - & 5.8149e-04 & - \\

 40 & 1.4642e-04  &1.94 & 2.0535e-04 & 1.74 & 1.5486e-04 & 1.91 \\

 80 & 3.6789e-05 &1.99 & 5.8150e-05 & 1.82 & 4.0122e-05  & 1.95\\

 160 & 9.2637e-06 &1.99& 1.5673e-05  & 1.89 & 1.0245e-05 & 1.97\\

320 & 2.3266e-06 &1.99 & 4.0796e-06  & 1.94 & 2.5881e-06  & 1.99  \\
 
640 & 5.8357e-07 &2.00 & 1.0417e-06 & 1.97 & 6.5033e-07 & 1.99 \\

 \bottomrule
\end{tabular}
\label{Tab2b}
\end{table}

\begin{table}[ht!]
 \caption{Error and order of convergence for Broadwell model with initial condition \eqref{smoothdata:Broadwell} for flux of momentum $z$}
 \centering
 \begin{tabular}{*{9}{c}}
 \toprule
\multirow{2}{*}{N} & 
\multicolumn{2}{c}{$\varepsilon=10^{-8}$} &
\multicolumn{2}{c}{$\varepsilon=0.02$} & 
\multicolumn{2}{c}{$\varepsilon=1$} \\
\cmidrule(lr){2-3}
\cmidrule(lr){4-5}
\cmidrule(lr){6-7}

& $L^{1}$ - Error & Order 
& $L^{1}$ - Error & Order 
& $L^{1}$ - Error & Order\\
\midrule

20  & 4.9267e-04 &-& 7.1371e-04 & - & 5.5222e-04 & - \\

 40 & 1.2781e-04  &1.95 & 2.0313e-04 & 1.81 & 1.4363e-04 & 1.94 \\

 80 & 3.2266e-05 &1.99 & 5.8082e-05 & 1.81 & 3.7115e-05  & 1.95\\

 160 & 8.1615e-06 &1.98 & 1.5600e-05  & 1.90 & 9.4782e-06 & 1.97\\

320 & 2.0571e-06 &1.99 & 4.0697e-06  & 1.94 & 2.3930e-06  & 1.99  \\
 
640 & 5.1617e-07 &1.99 & 1.0397e-06 & 1.97 & 6.0141e-07 & 1.99 \\

 \bottomrule
\end{tabular}
\label{Tab3b}
\end{table}
\begin{figure}[!ht]
     \centering
     \begin{subfigure}[b]{0.3\textwidth}
         \centering  
         \includegraphics[width=\linewidth, height=3.5cm]{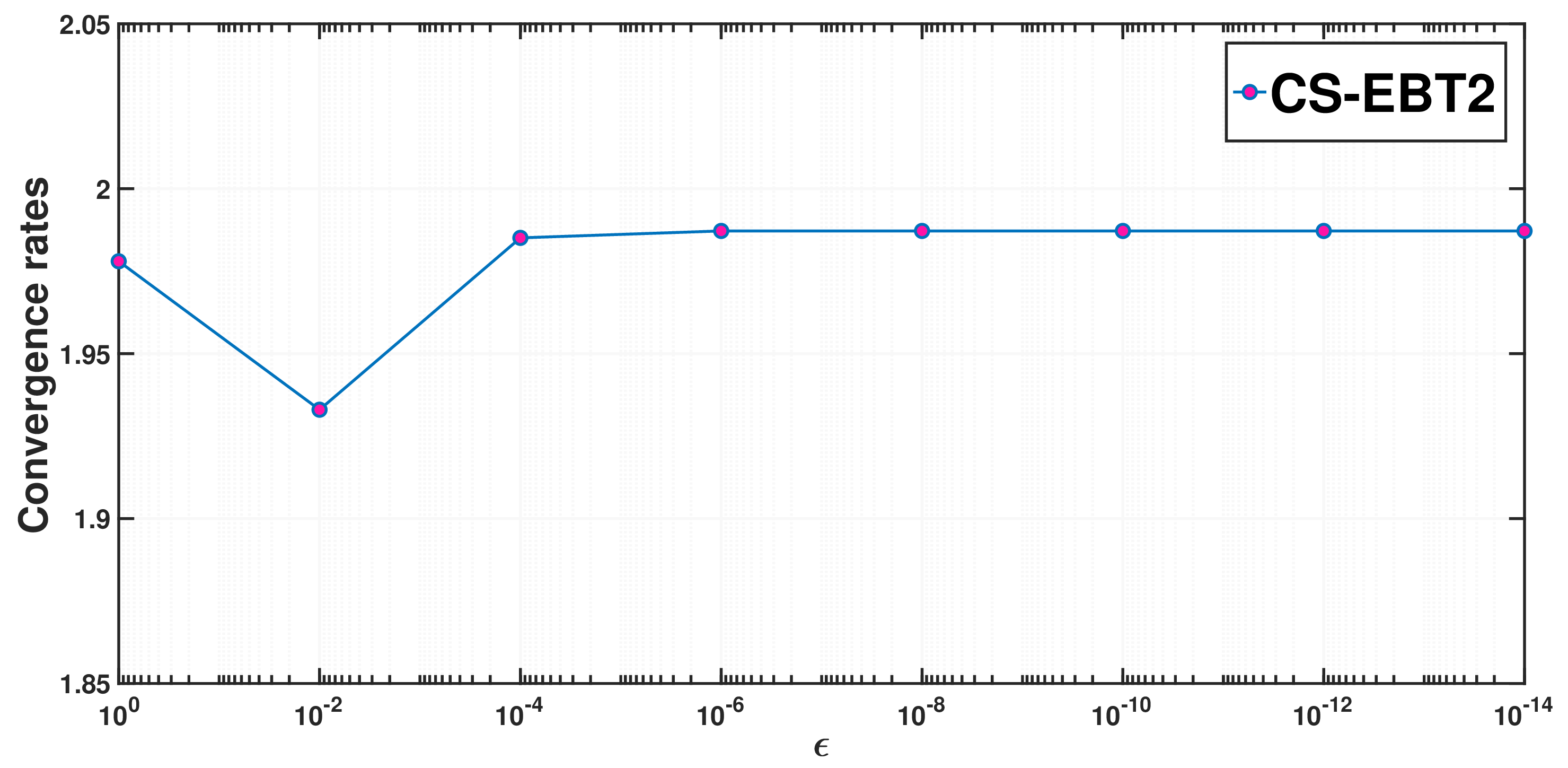}
      \caption{Convergence rate for $\rho$}
         \label{(fig:Broadwell_s1)}
     \end{subfigure}
     \hfill
     \begin{subfigure}[b]{0.3\textwidth}
         \centering   
         \includegraphics[width=\linewidth, height=3.5cm]{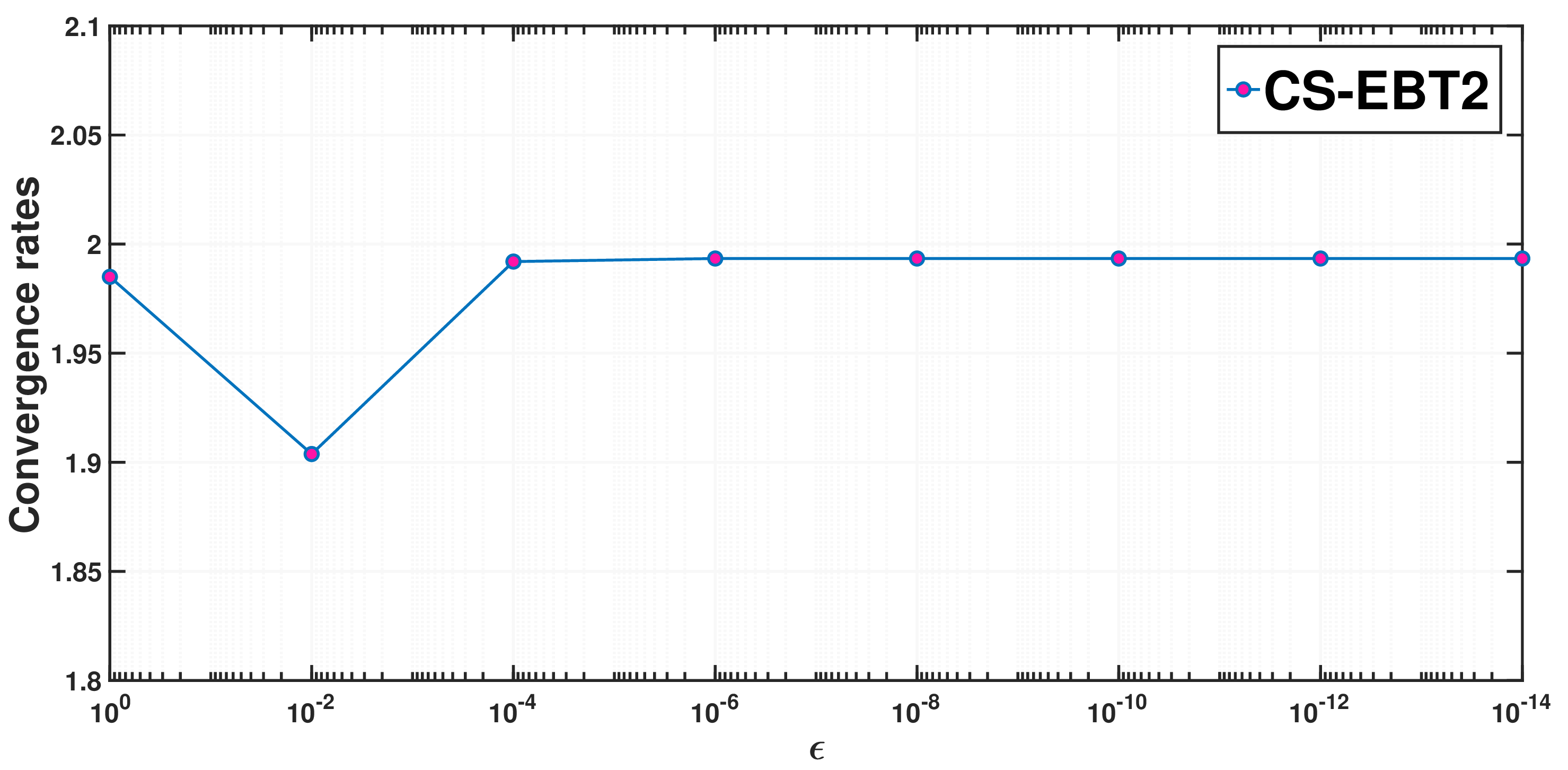}
         \caption{Convergence rate for $m$}
         \label{(fig:Broadwell_s2)}
     \end{subfigure}
     \hfill
    \begin{subfigure}[b]{0.3\textwidth}
         \centering
      \includegraphics[width=\linewidth, height=3.5cm]{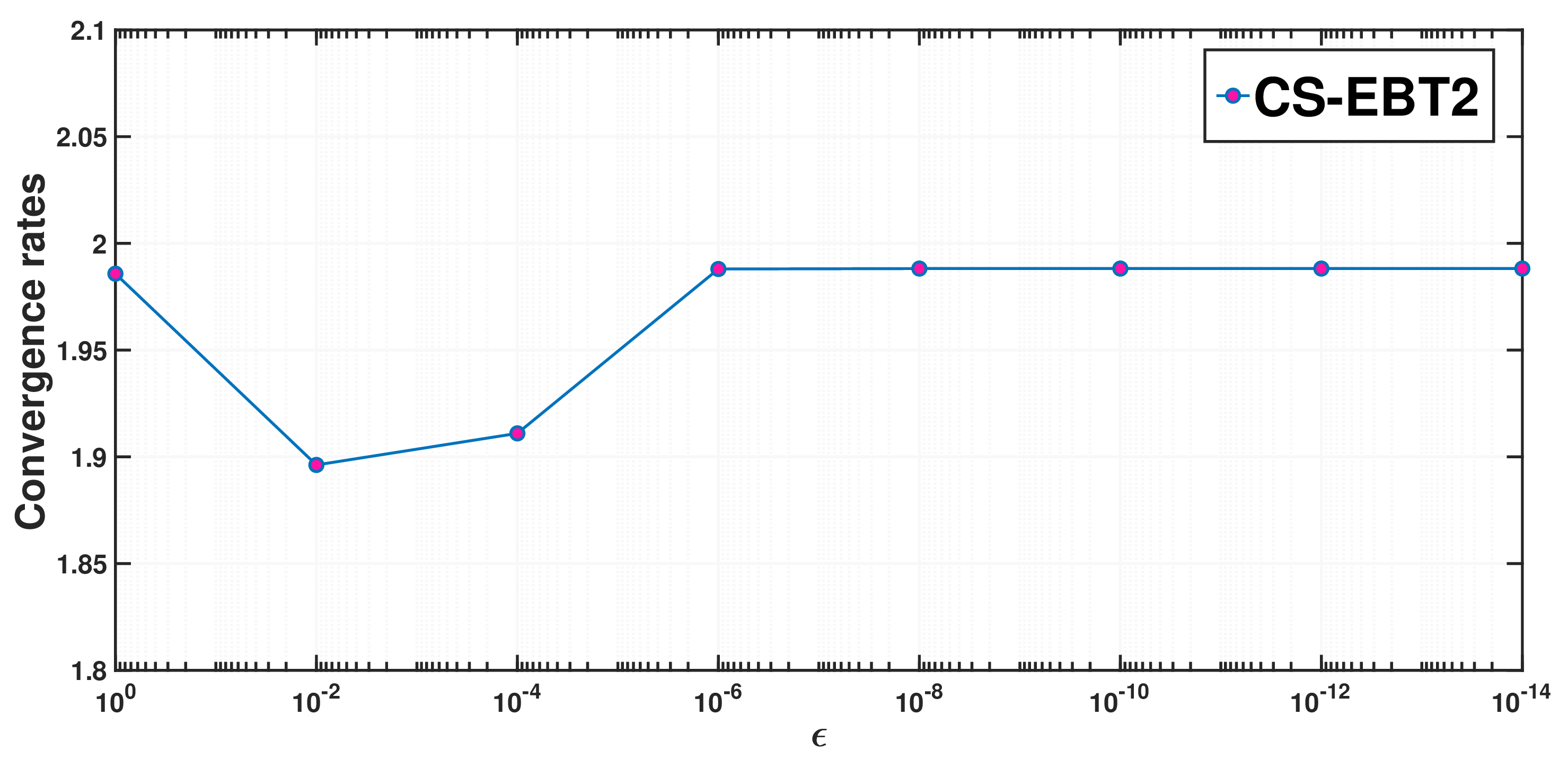}
     \caption{Convergence rate for z}
         \label{(fig:Broadwell_s3)}
     \end{subfigure}
         \caption{Convergence rates for the numerical solutions $\rho$, $m$ ans $z$ corresponding to the smooth case of the Broadwell model, computed on the domain $[0,1]$ up to the final time $T=0.3$, with a CFL number of $0.9$, using $\varepsilon \in \{10^0, 10^{-2}, 10^{-4}, 10^{-6},10^{-8},10^{-10}, 10^{-12}, 10^{-14}\}$.}
    \label{(broadwell:order)}
\end{figure}
This setup guarantees a smooth solution throughout the simulation, making it well-suited for evaluating the accuracy and stability of the numerical scheme in a regime free of discontinuities. The computations are performed using a uniform grid, and a CFL of $0.9$ upto the final time $T=0.3$, is applied to ensure numerical stability. 
To evaluate the scheme’s performance, we report results for various values of the relaxation parameter $\varepsilon$, specifically $10^{-8}, 0.02$ and $1$. The corresponding errors and convergence orders have been provided for the variables $\rho, m$ and $z$ in Tables \ref{Tab1b}, \ref{Tab2b} and \ref{Tab3b} which contain the $L^{1}$-errors, respectively. These tables demonstrate that the scheme achieves second-order accuracy.

To illustrate the uniform second-order accuracy of the scheme, we consider a broader range of relaxation parameters $\varepsilon \in \{ 10^{-2\ell}, \ell = 0,\ldots,7 \}$. The convergence behavior for this set of values is depicted in Figure \ref{(broadwell:order)}, which displays the order of accuracy for the smooth Broadwell case, computed on the interval $[0,1]$, using CFL$=0.9$ and final time T$= 0.3$. 
In this smooth Broadwell setup the solution stays closer to equilibrium also for $\varepsilon=1$, hence the kinetic-regime dissipative defect is less visible and the EOC remains close to two.
\begin{figure}[!ht]
    \begin{minipage}[b]{0.48\linewidth}
        \includegraphics[width=\linewidth]{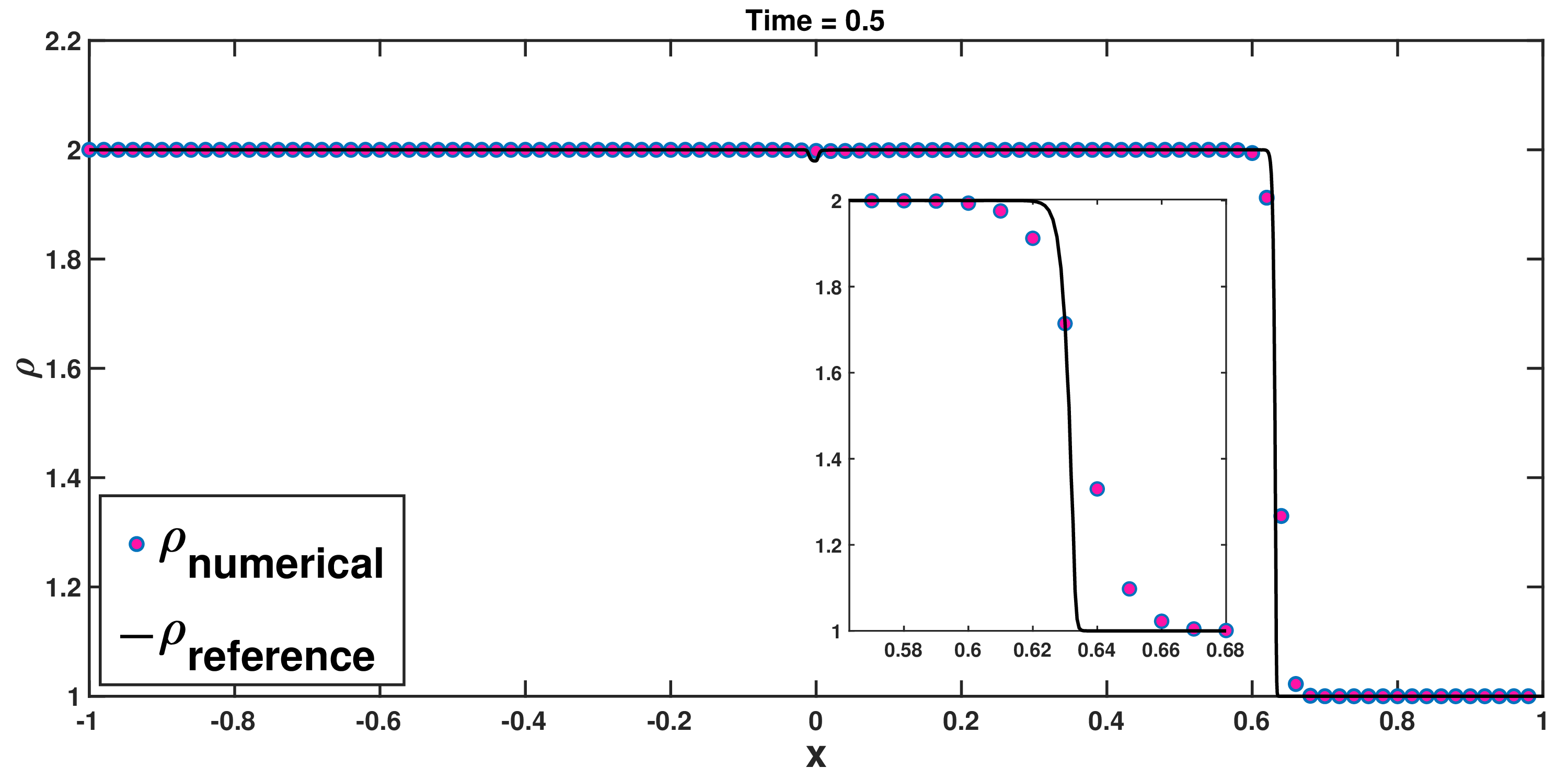}
    \end{minipage}
    \hfill
    \begin{minipage}[b]{0.48\linewidth}
        \includegraphics[width=\linewidth]{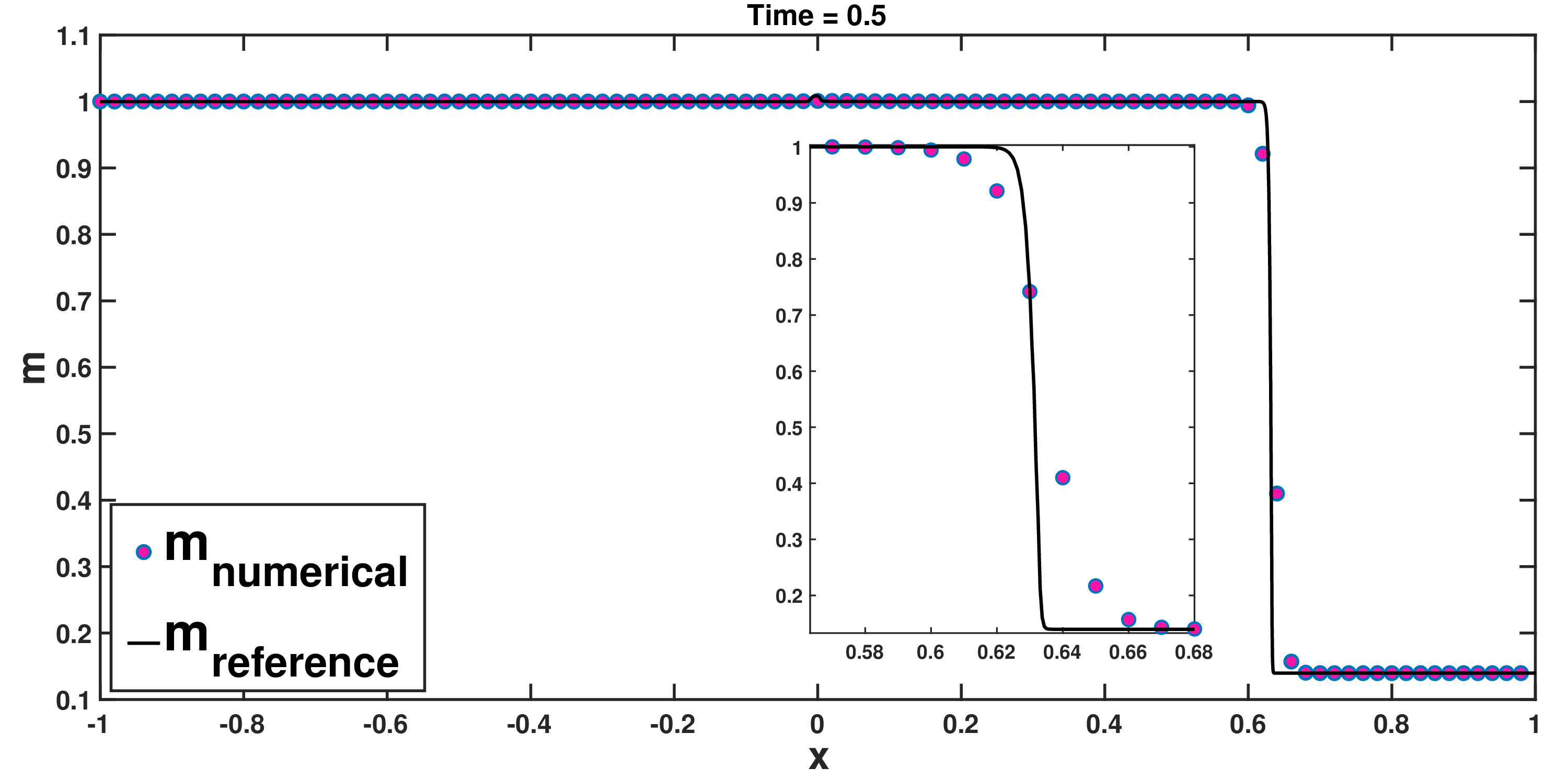}
    \end{minipage}

    \vspace{0.5cm} 

    \centering
    \includegraphics[width=0.5\linewidth]{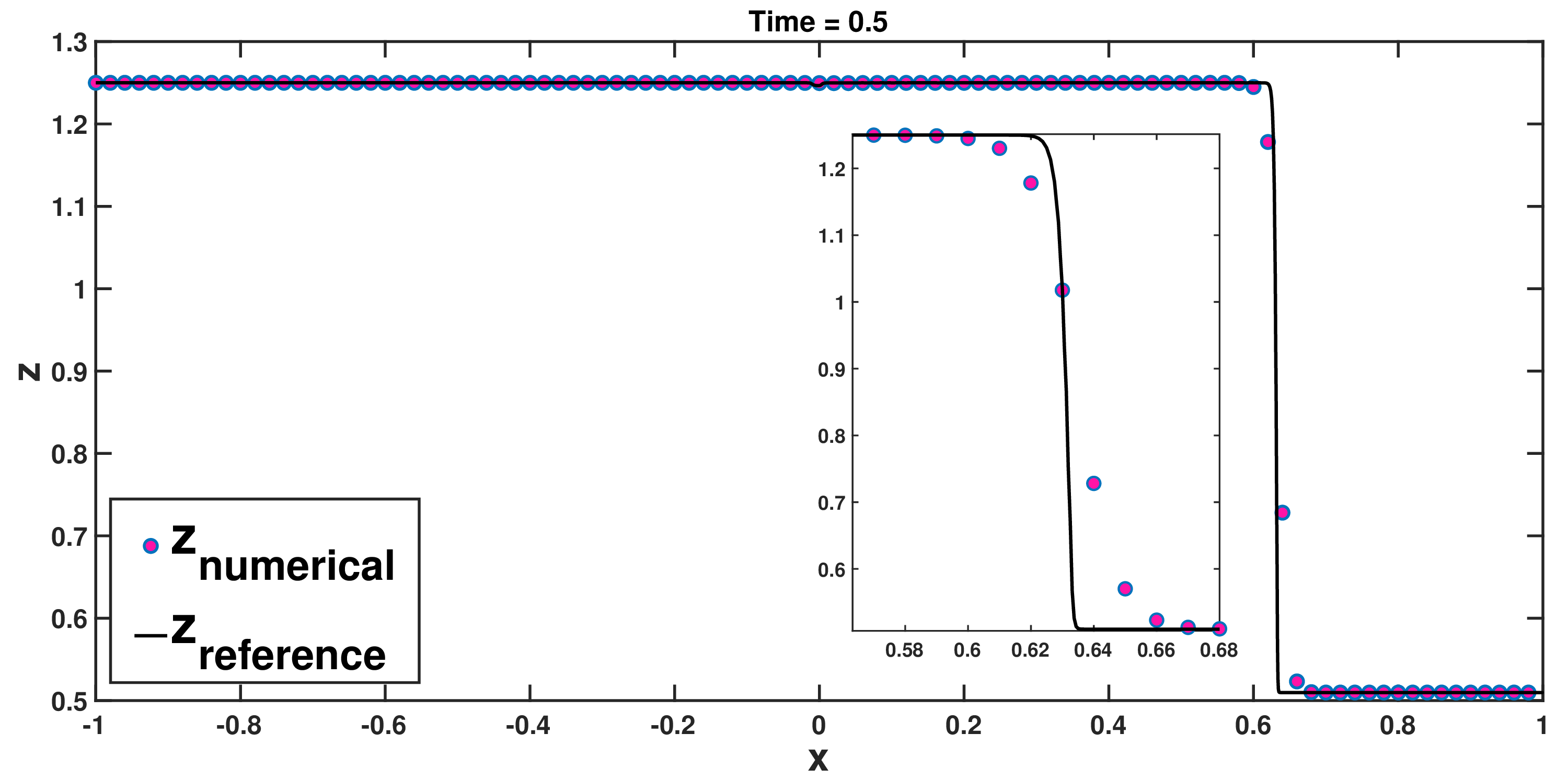}
    
    \caption{Broadwell model with non-smooth case: comparison between Numerical solution $\rho$(left), $m$(right) and $z$(center) and the reference solution (IMEX-RK2) with $\varepsilon = 10^{-8}$, CFL $0.9$ and $N=200$.}
    \label{(3a)}
\end{figure}
\begin{figure}[!ht]
    \begin{minipage}[b]{0.49\linewidth}
        \includegraphics[width=\linewidth]{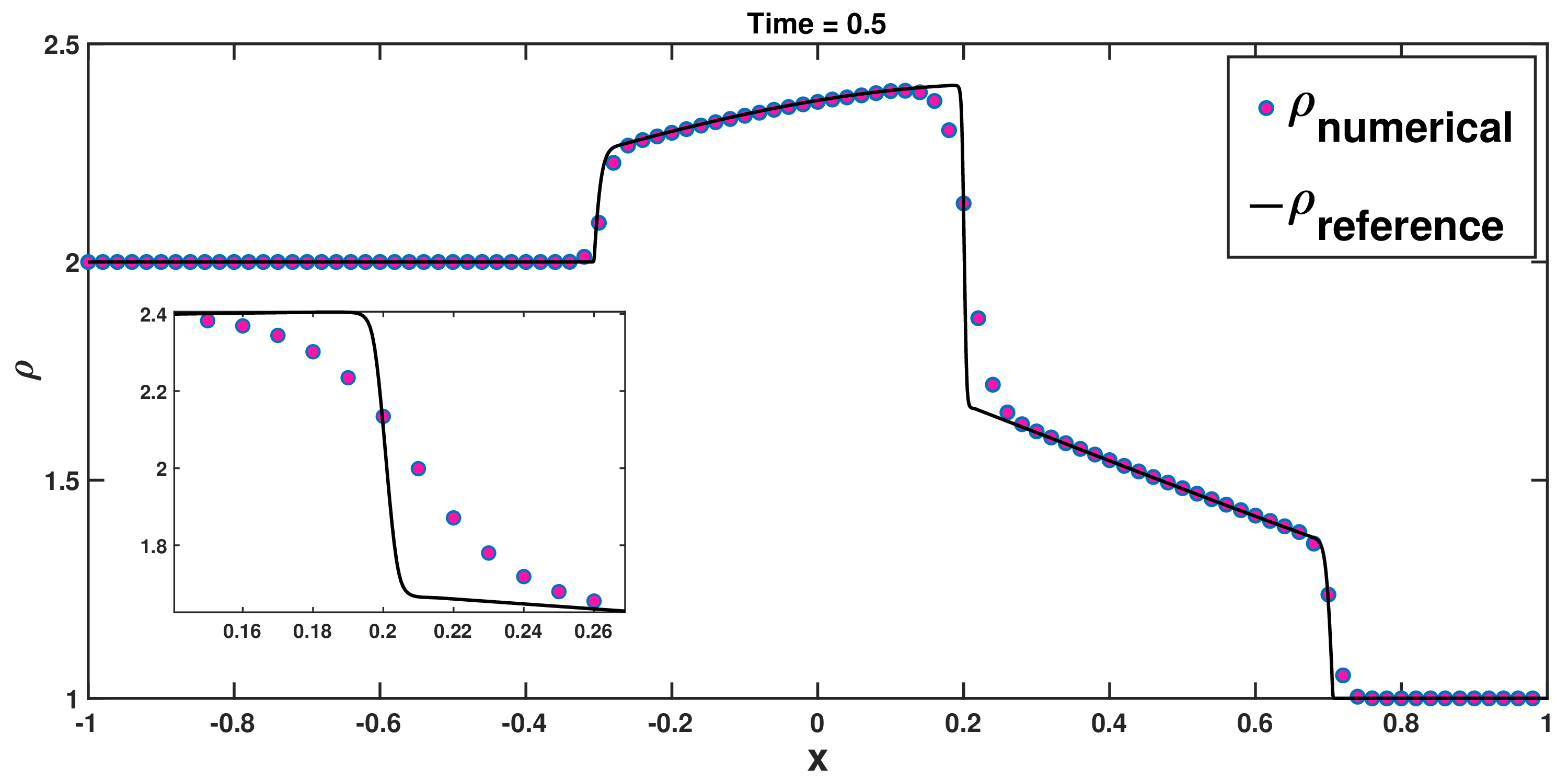}
    \end{minipage}
    \hfill
    \begin{minipage}[b]{0.49\linewidth}
        \includegraphics[width=\linewidth]{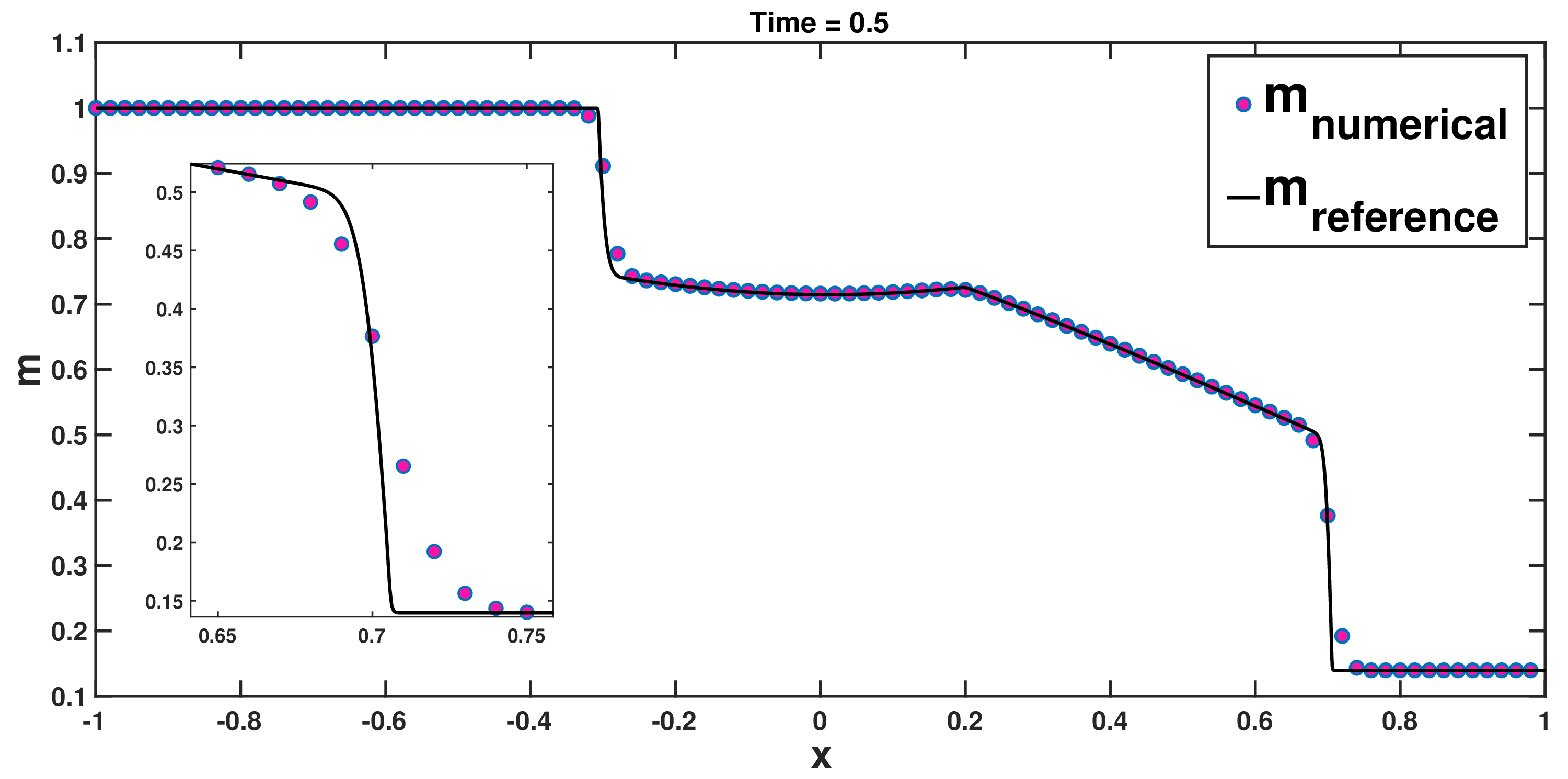}
    \end{minipage}

    \vspace{0.5cm} 

    \centering
    \includegraphics[width=0.5\linewidth]{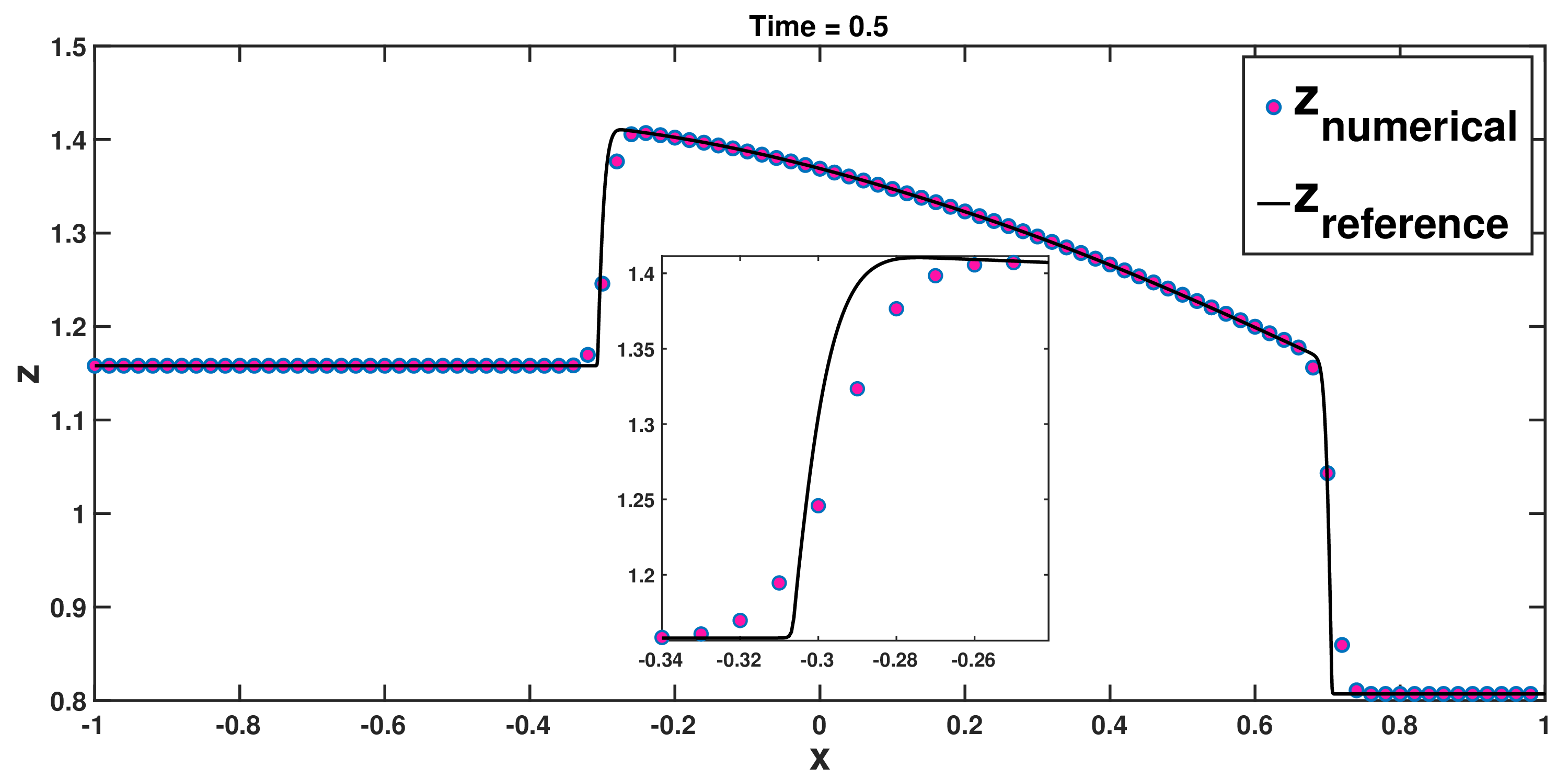}
    
    \caption{Broadwell model with non-smooth case: comparison between Numerical solution $\rho$(left), $m$(right) and $z$(center) and the reference solution (IMEX-RK2) with $\varepsilon = 1$, CFL $0.9$ and $N=200$.}
    \label{(3b)}
\end{figure}

\noindent{\textbf{Non-Smooth case:}}
In the context of the Broadwell model, we consider a non-smooth initial profile to evaluate the robustness and accuracy of the proposed numerical scheme. The initial conditions for the conservative variables $\rho$, $m$ and $z$ are defined as
\begin{eqnarray}
    \begin{aligned}
    (\rho(x,0), m(x,0), z(x,0))=
        \begin{cases}
            (2,1,1), \, \text{if}\,\,\,\,x\leq0.2,\\
            (1,0.13962,1), \, \text{if}\,\,\,\,x>0.2,
            \end{cases} \\
    \end{aligned}
\end{eqnarray}
which introduces a discontinuity at $x = 0.2$. This sharp interface serves as a useful test case for examining how well the numerical method handles strong gradients, particularly in the stiff relaxation regime. Transmissive boundary conditions are imposed at the endpoints of the domain $[-1,1]$, and the solution is evolved up to the final time $T = 0.5$.

\begin{figure}[!ht]
    \begin{minipage}[b]{0.48\linewidth}
        \includegraphics[width=\linewidth]{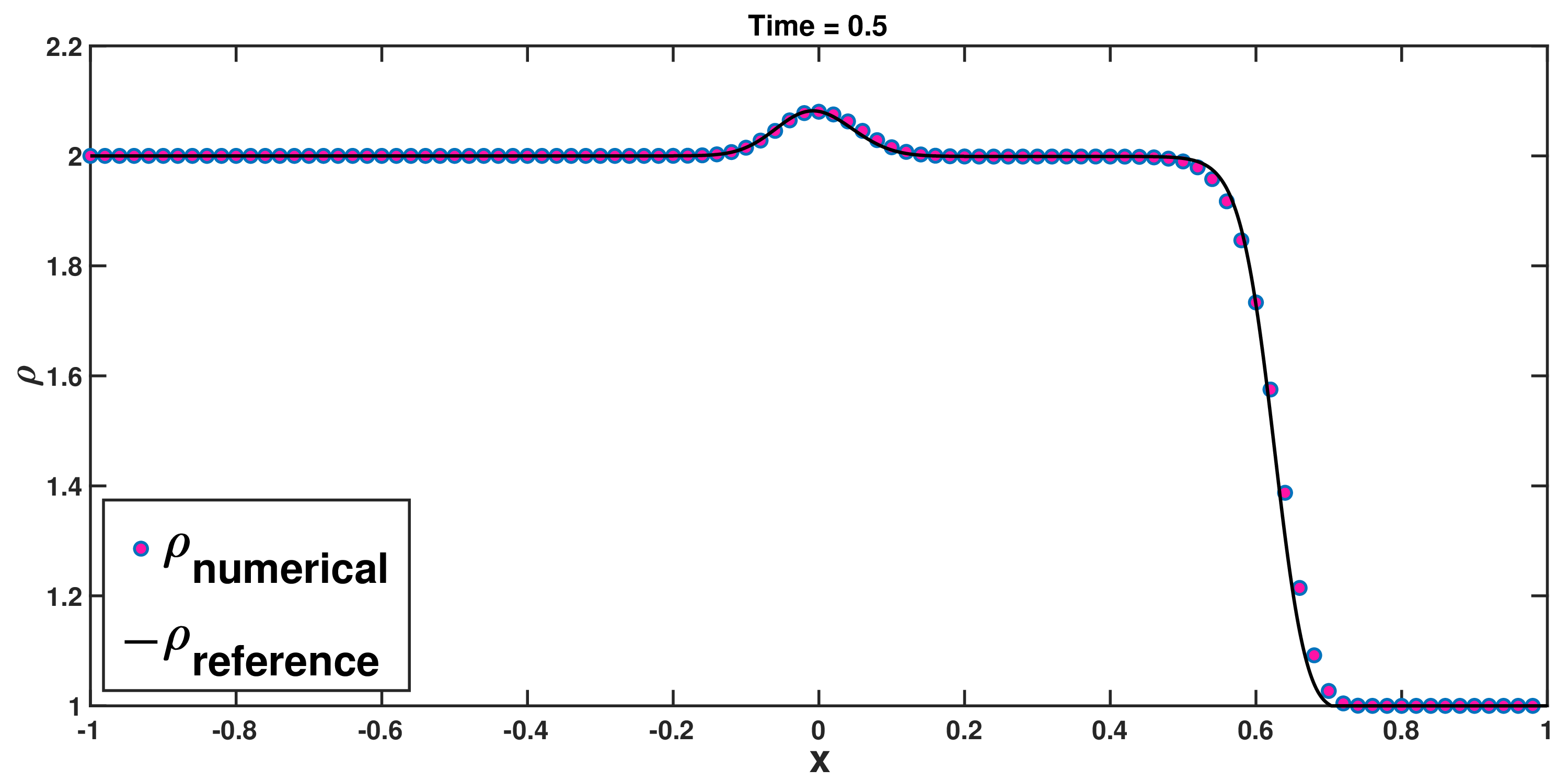}
    \end{minipage}
    \hfill
    \begin{minipage}[b]{0.48\linewidth}
        \includegraphics[width=\linewidth]{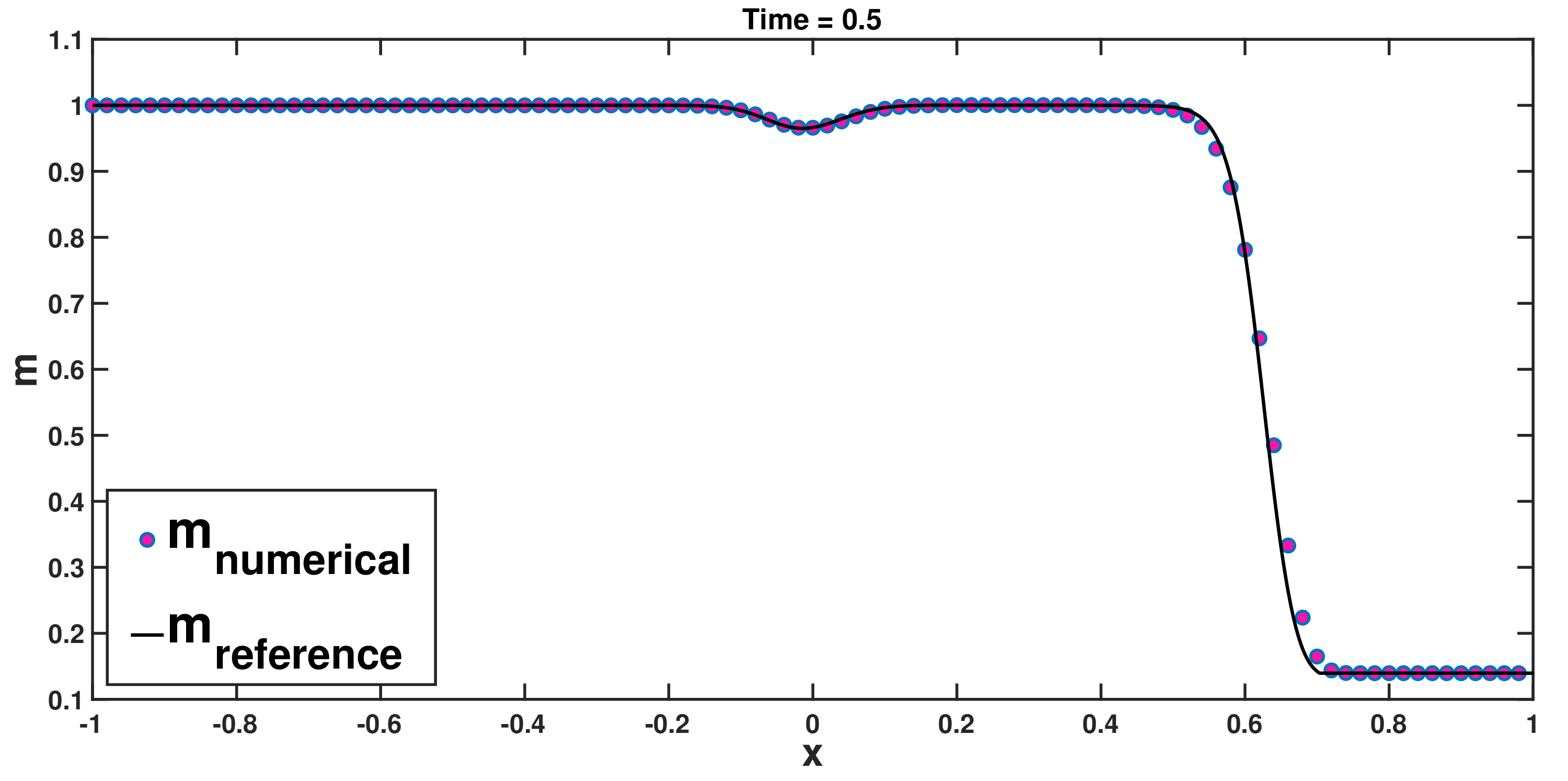}
    \end{minipage}

    \vspace{0.5cm} 

    \centering
    \includegraphics[width=0.5\linewidth]{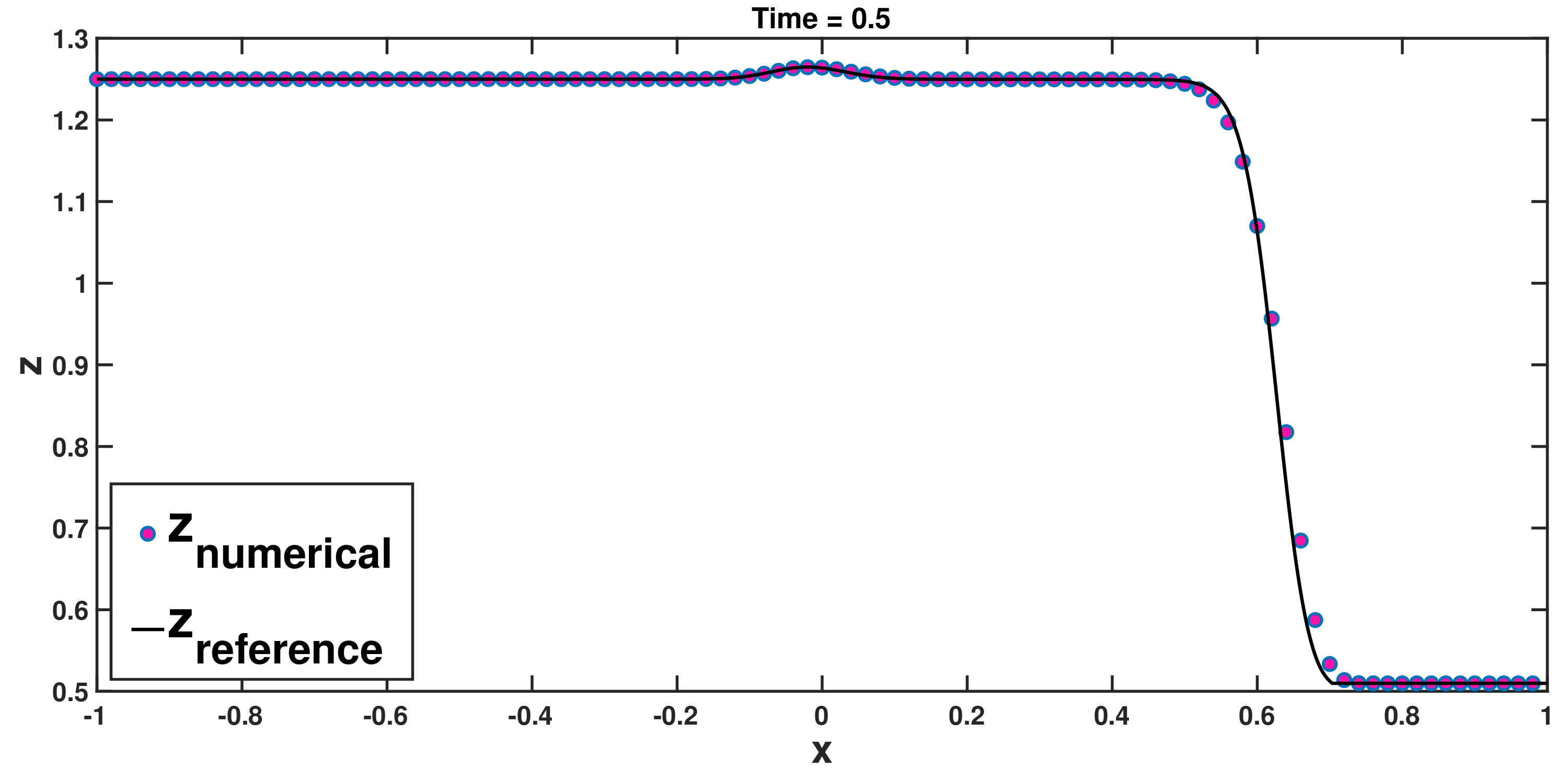}
    
    \caption{Broadwell model with non-smooth case: comparison between Numerical solution $\rho$(left), $m$(right) and $z$(center) and the reference solution (IMEX-RK2) with $\varepsilon = 0.02$, CFL $0.9$ and $N=200$.}
    \label{(3c)}
\end{figure}
\par
Figure~\ref{(3a)} illustrates the numerical solutions for the density ($\rho$), momentum ($m$), and flux momentum ($z$) represented from left to right, for the case $\varepsilon = 10^{-8}$. These results are compared with a reference solution obtained using the second-order IMEX Runge–Kutta (IMEX-RK2) method computed on a highly refined mesh with $N=3200$ grid points. The proposed scheme is implemented using $N=200$ uniformly distributed grid points and a CFL number of $0.9$. As observed, the scheme performs well in the stiff regime, accurately capturing the solution structure and maintaining consistency with the reference profile. To explore the scheme’s adaptability under varying stiffness, additional simulations are carried out for two larger values of the relaxation parameter, namely $\varepsilon=1$ and $\varepsilon=0.02$. The results corresponding to these cases are displayed in Figures~\ref{(3b)} and~\ref{(3c)}, respectively. To highlight the scheme’s resolution properties near discontinuities, zoomed-in portions around the interface are included for both the highly stiff case $\varepsilon = 10^{-8}$ and the fully kinetic regime $\varepsilon= 1$. These views provide detailed insight into the scheme’s capability to resolve sharp transitions without introducing spurious oscillations or numerical diffusion.

Moreover, although no zoomed plot is shown explicitly for the intermediate case $\varepsilon=0.02$, the combined observations from the extreme cases help to qualitatively interpret the behavior at this intermediate stiffness level. Across all tested regimes, the proposed method consistently demonstrates strong agreement with the reference solution. In particular, it exhibits the desired robustness and accuracy near discontinuities, confirming its reliability for simulating relaxation systems across both stiff and non-stiff regimes. 
\subsection{Euler equation with heat transfer Model}
The model Euler equation with heat transfer \eqref{Euler_heat} is supplemented with transmissive boundary conditions at both boundaries in the case of non-smooth profiles,
\begin{equation*}
\begin{cases}
\rho_t + (\rho u)_x = 0,\\
(\rho u)_t + (\rho u^2 + p)_x = 0,\\
(\rho E)_t + (\rho u E + pu)_x = -K\rho (T-T_0),\qquad K=1/\varepsilon.
\end{cases}
\end{equation*}
\noindent{\textbf{Non-Smooth case:}}
To evaluate the performance of the numerical scheme in capturing sharp discontinuities and non-equilibrium dynamics, an unprepared non-smooth initial condition has been imposed within the spatial domain. The initial data for the conserved variables are density $(\rho)$, velocity ($u$) and total energy ($E$) is given by
\begin{eqnarray}
    \begin{aligned}
    (\rho(x,0), u(x,0), E(x,0))=
        \begin{cases}
            (1,0,1), \, \text{if}\,\,\,x\leq0.5\\
            (1/8,0,1) \, \text{if}\,\,\,x>0.5.
            \end{cases} \\
    \end{aligned}
\end{eqnarray}
This discontinuous profile introduces a jump at the midpoint of the domain and is specifically chosen to test the robustness of the numerical method under challenging initial configurations. Transmissive boundary conditions are applied at both ends of the domain. The governing model is the Euler system with heat transfer, where the small parameter $\varepsilon = 10^{-8}$ represents a regime close to thermal equilibrium.

Figure~\ref{(4a)} presents the numerical results for the conservative variables: density $\rho$ (left), momentum
$\rho u$ (right), and total energy $\rho E $(center). These results are compared with a reference solution obtained using a second-order IMEX Runge–Kutta (IMEX-RK2) scheme with a highly refined grid of  $N=3200$ points. The proposed scheme is applied on a coarser mesh with $N=200$ uniformly spaced grid points over the spatial interval $[0,1]$, using CFL$=0.9$ and a final simulation time of $T=0.3$.
\begin{figure}[!ht]
    \begin{minipage}[b]{0.48\linewidth}
        \includegraphics[width=\linewidth]{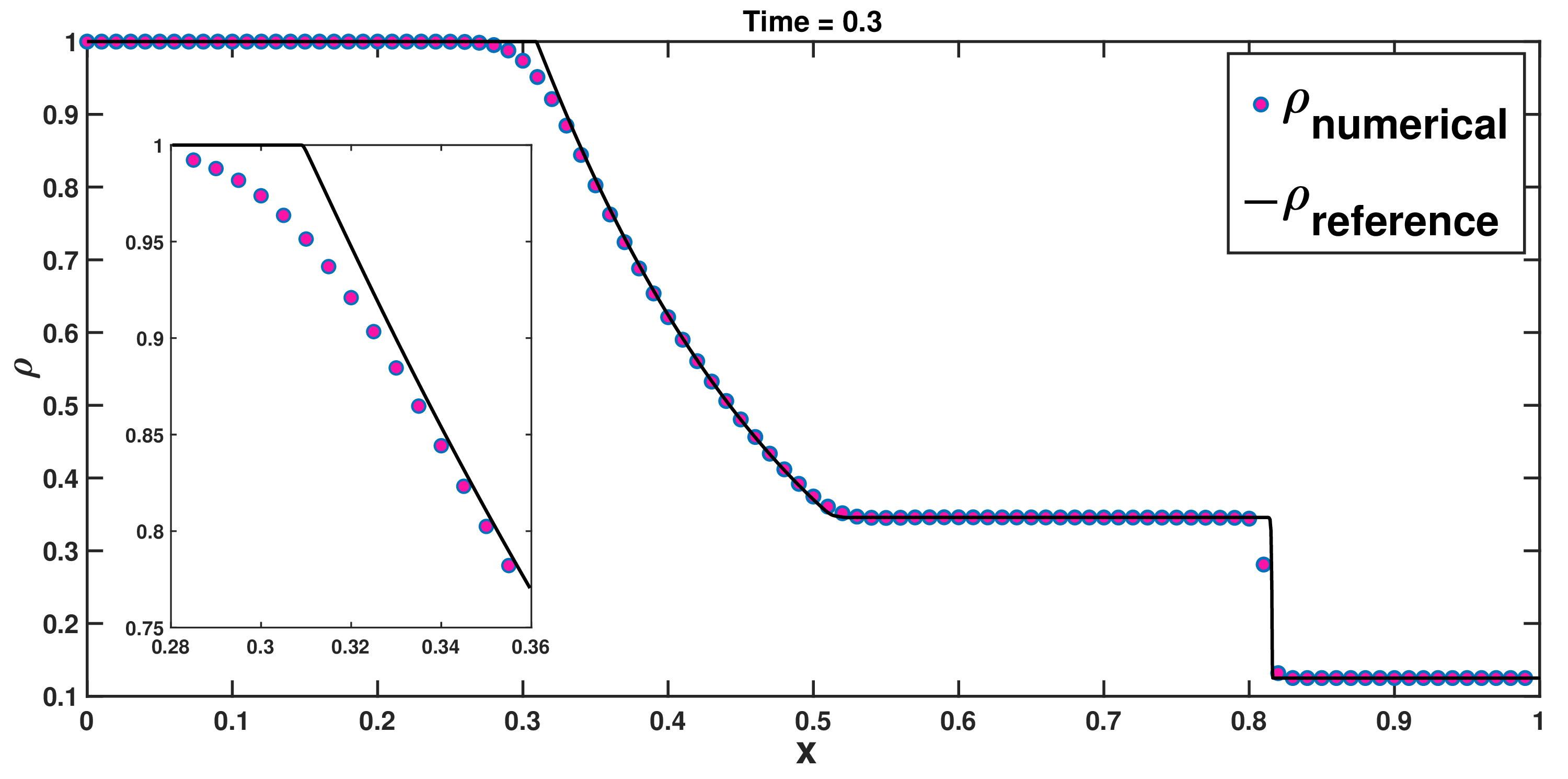}
    \end{minipage}
    \hfill
    \begin{minipage}[b]{0.48\linewidth}
        \includegraphics[width=\linewidth]{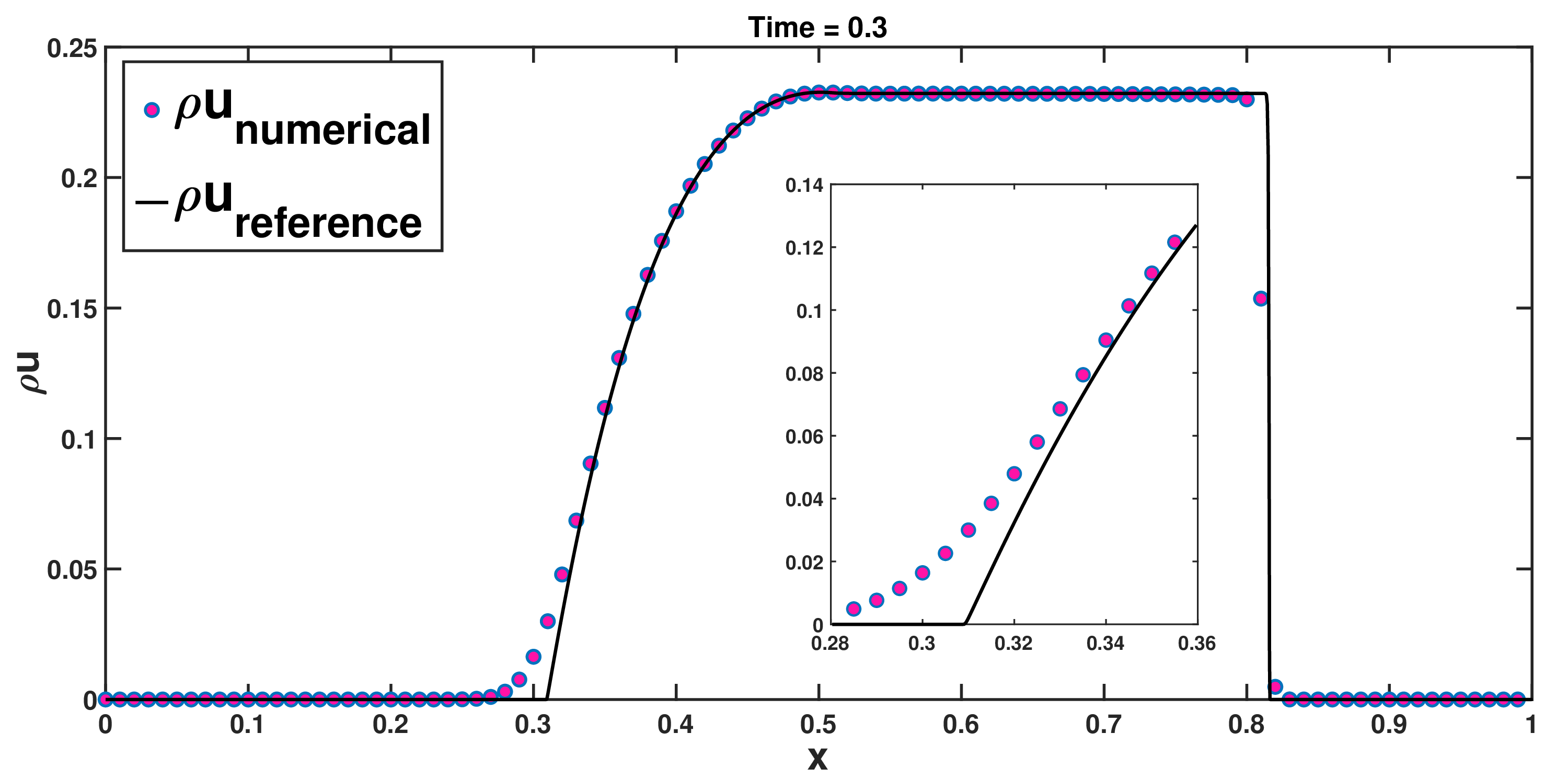}
    \end{minipage}

    \vspace{0.5cm} 

    \centering
    \includegraphics[width=0.5\linewidth]{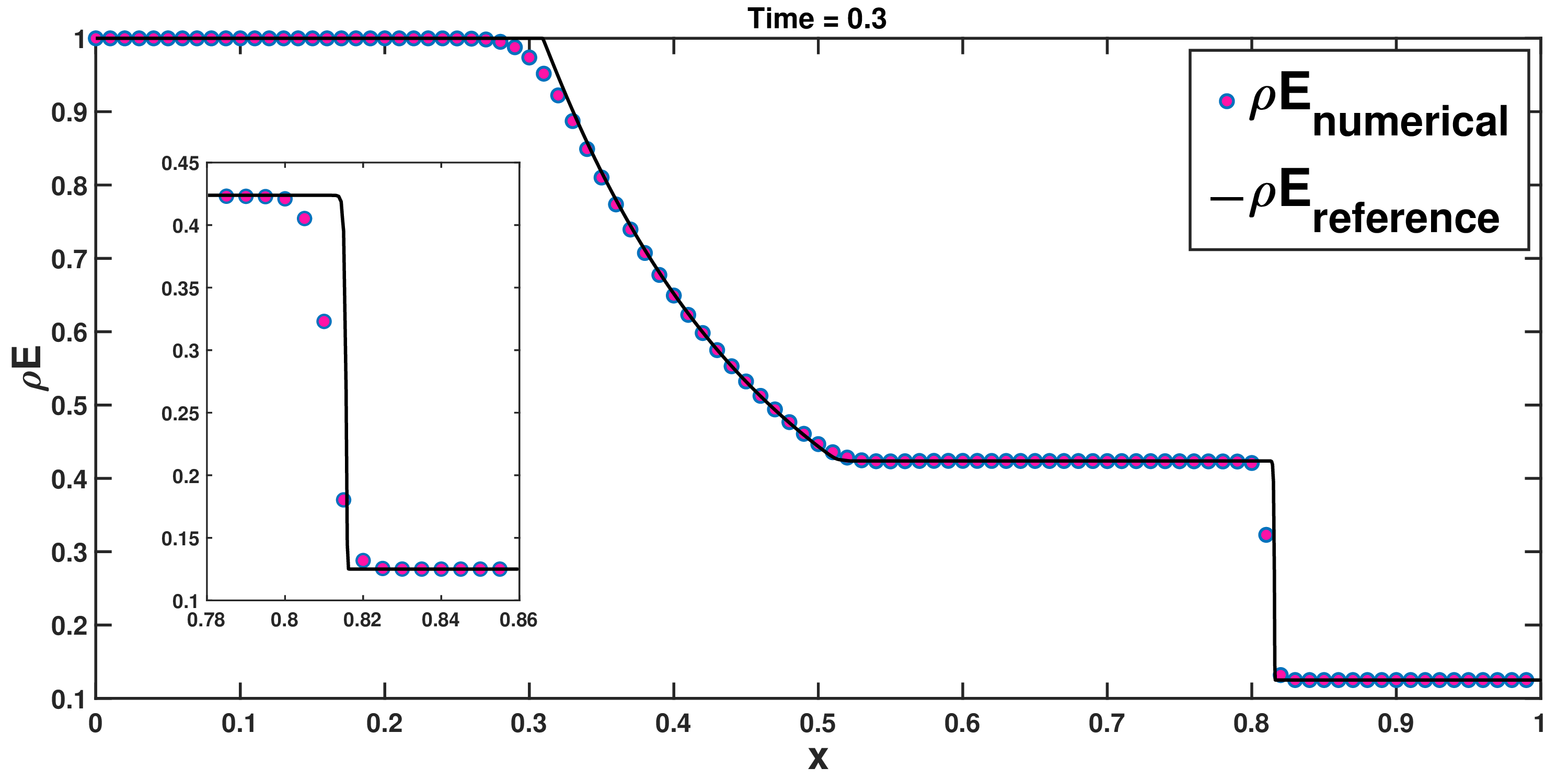}
    
    \caption{Euler heat transfer model with non-smooth case: comparison between Numerical solution $\rho$(left), $\rho u$(right) and $\rho E$(center) and the reference solution (IMEX-RK2) with $\varepsilon = 10^{-8}$, CFL $0.9$ and $N = 200$.}
    \label{(4a)}
\end{figure}
The plotted results demonstrate that the proposed numerical scheme produces accurate and stable solutions, even in the presence of strong stiffness. To further clarify the method’s ability to handle sharp gradients and discontinuities, a zoomed-in view of the solution near the discontinuity point has been included in each plots. This enlarged portion provides a closer look at how the scheme captures the steep transition and maintains consistency across the stiff region. The close agreement with the reference solution in this critical region supports the robustness and reliability of the scheme under challenging conditions.

\subsection{Euler equation with stiff friction}
In this example, we consider the Euler equations of compressible gas dynamics with stiff friction. Two variants are examined: the full Euler system \eqref{eulerGas:stiff_Friction} and the isentropic Euler system with stiff friction \eqref{eulerGas_isentropic:stiff_friction}. In both cases, transmissive boundary conditions are prescribed at the domain boundaries,
\begin{equation*}
\begin{cases}
\rho_t + (\rho u)_x = 0,\\
(\rho u)_t + (\rho u^2 + p)_x = -\alpha \rho u,\qquad \alpha=1/\varepsilon,\\
(\rho E)_t + ((\rho E+p)u)_x = -\alpha \rho u^2.
\end{cases}
\end{equation*}

\textbf{Non-Smooth case:} The initial condition for the conserved variables, density ($\rho$), velocity ($u$) and pressure ($p$) is prescribed by
\begin{eqnarray}
    \begin{aligned}
    (\rho(x,0), u(x,0), p(x,0))=
        \begin{cases}
            (1.65,0,5.039849068), \, \text{if}\,\,\,x\leq0.25\\
            (0.01,0,0.003962233) \, \text{if}\,\,\,x>0.25.
            \end{cases} \\
    \end{aligned}
\end{eqnarray}
The experiments are carried out with parameters $\text{CFL} = 0.9$, $\gamma = 1.4$, and $\alpha = 10^{8}$. The computational domain is set to $[0,1]$, and the solution is evolved up to the final time $T = 2.0$. Numerical results for $\rho$, $\rho u$, and $\rho E$ are obtained using $N = 1000$ grid points. For validation, these results are compared with high-resolution reference solutions computed on a finer grid with $N = 4000$, employing a second-order IMEX Runge–Kutta (IMEX-RK2) scheme. This comparison allows us to assess the accuracy of the proposed scheme at the same resolution.

\subsubsection{Full Euler system with stiff friction}   
 \begin{figure}[!ht]
    \begin{minipage}[b]{0.48\linewidth}
        \includegraphics[width=\linewidth]{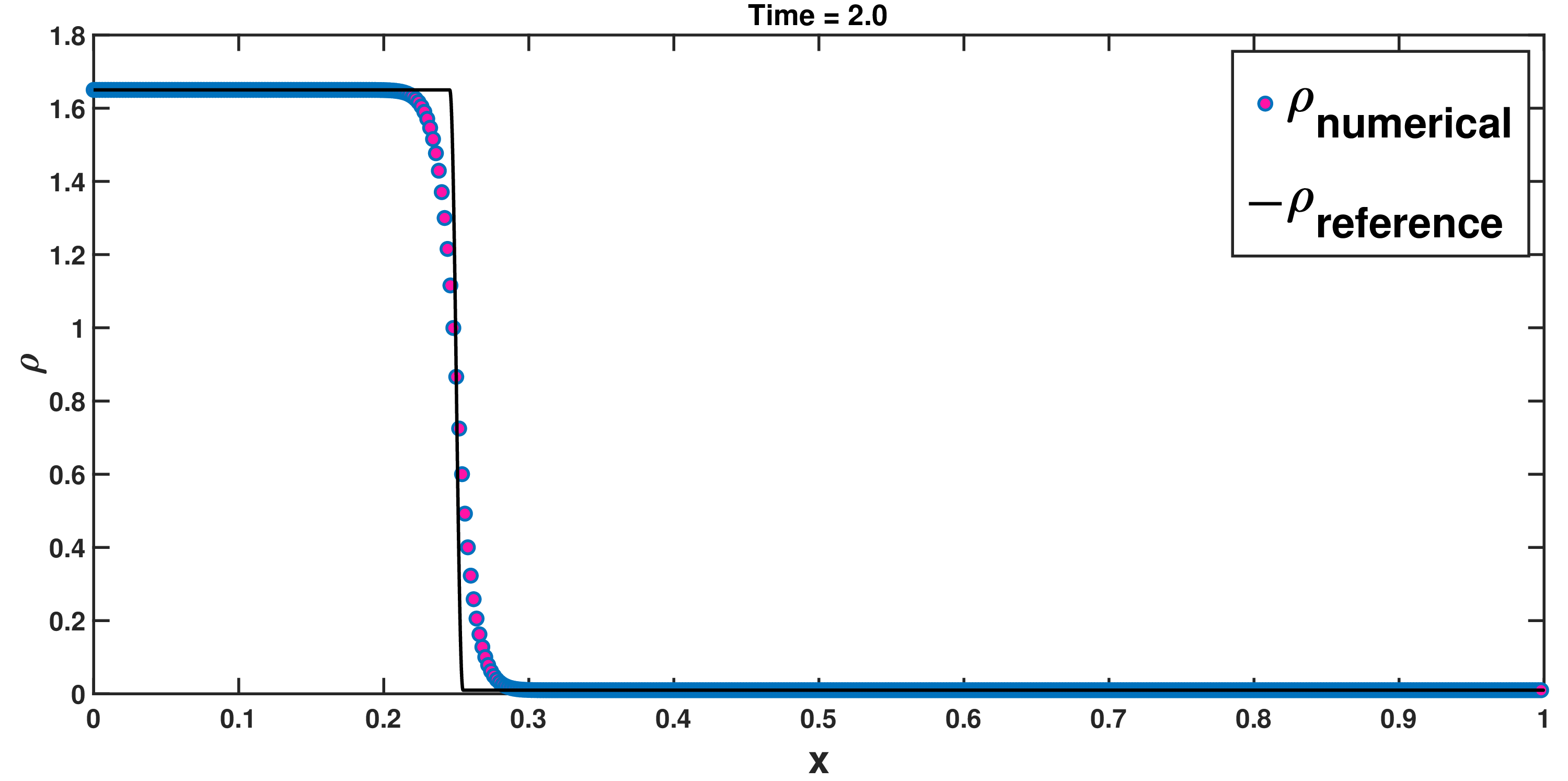}
    \end{minipage}
    \hfill
    \begin{minipage}[b]{0.48\linewidth}
        \includegraphics[width=\linewidth]{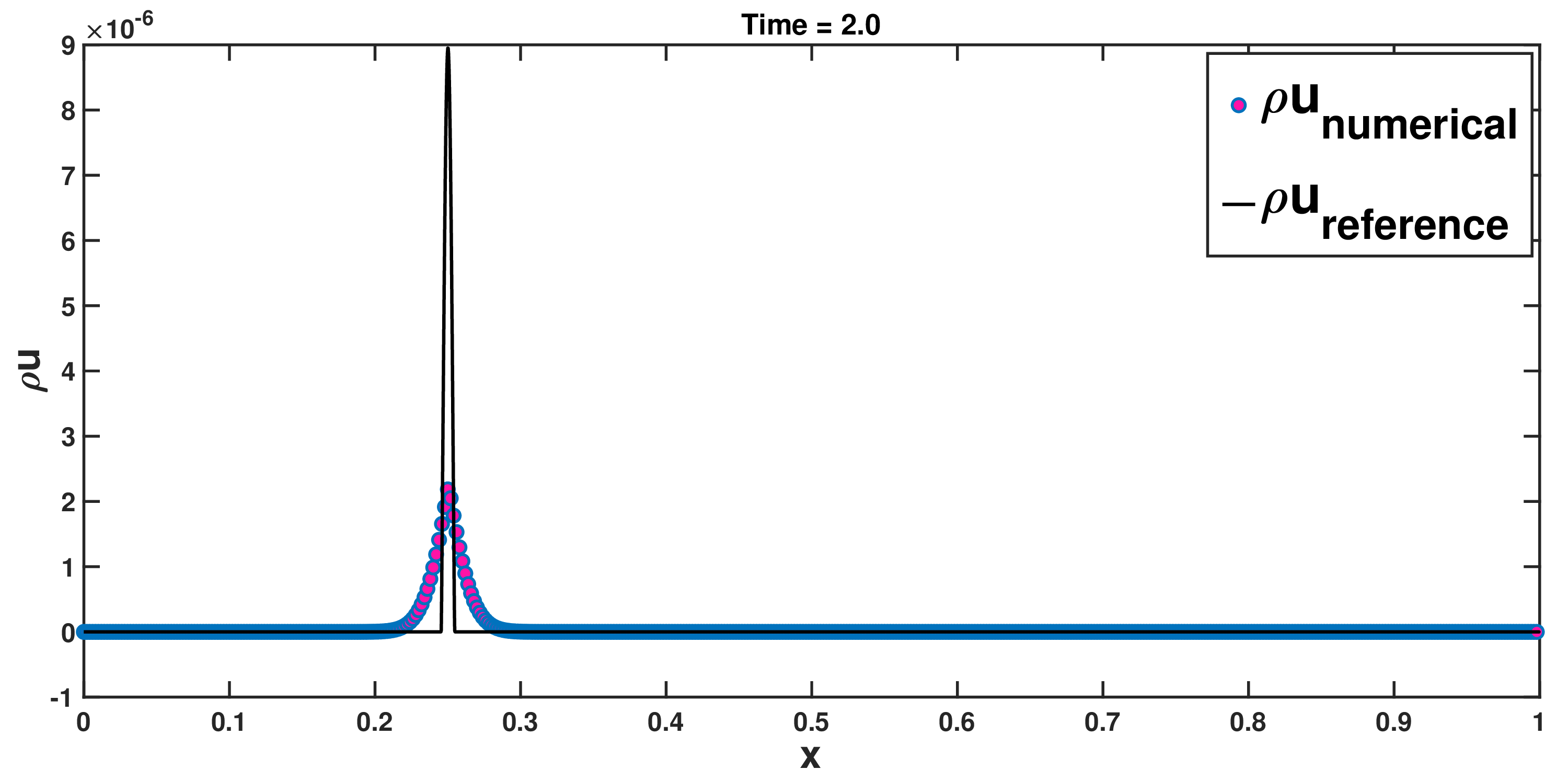}
    \end{minipage}

    \vspace{0.5cm} 

    \centering
    \includegraphics[width=0.5\linewidth]{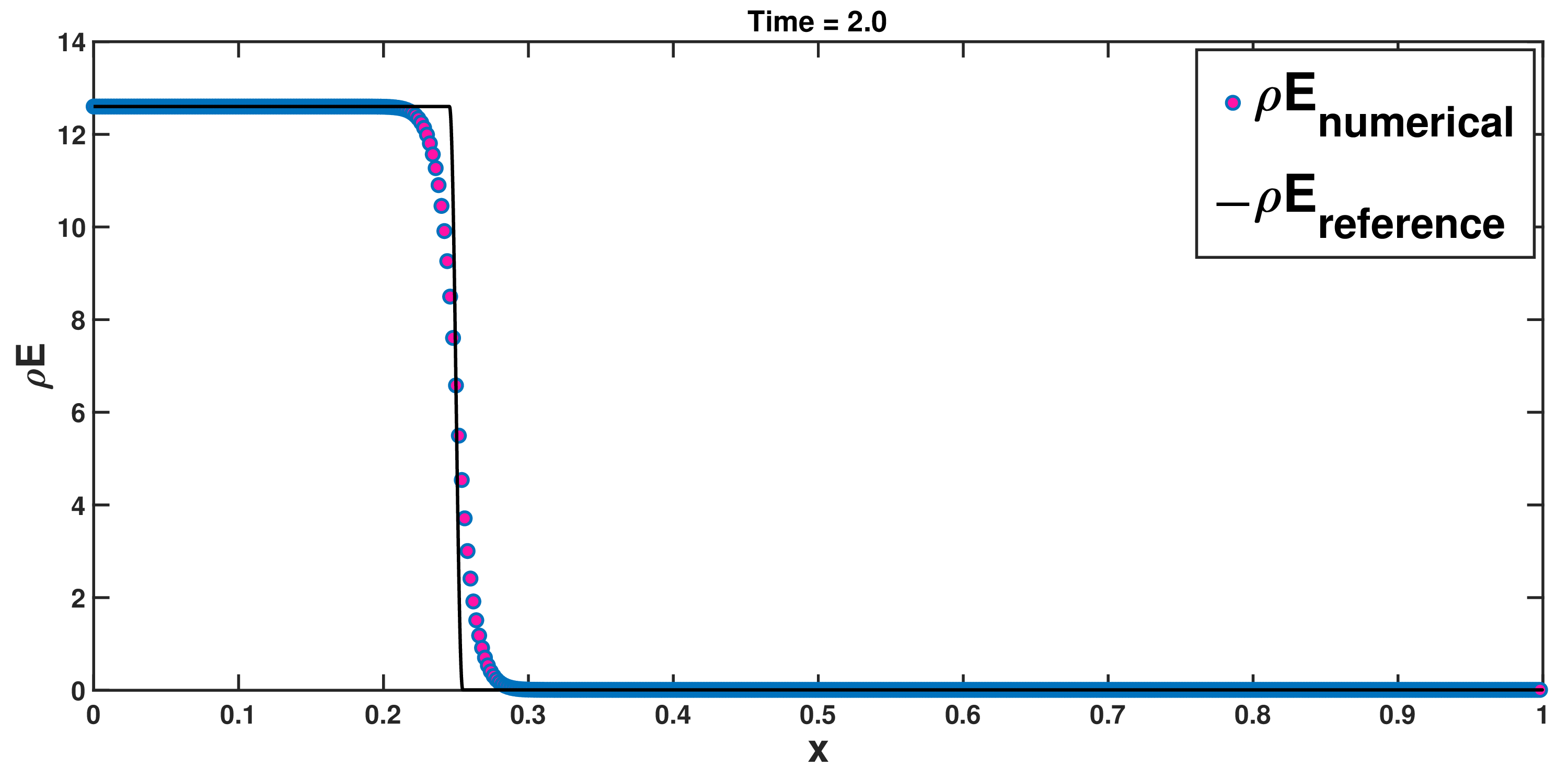}
    
    \caption{Full Euler system with stiff friction model in the non-smooth case: comparison between the numerical and reference (IMEX-RK2) solutions. The top row shows $\rho$ (left), $\rho u$ (right) and $\rho E$ (center) for $\alpha = 10^{8}$, CFL $0.9$, and $N=1000$.}
    \label{Eulergasfig:example1}
\end{figure}
For the full Euler equations, the pressure is given by
\begin{equation*}
p = (\gamma - 1)\left(\rho E - \tfrac{1}{2}\rho u^2\right).
\end{equation*}
The numerical results are presented in Figure~\ref{Eulergasfig:example1}, where $\rho$ (left), $\rho u$ (center), and $\rho E$ (right) are compared with the corresponding high-resolution reference solutions. Overall, the agreement is very good, confirming the robustness of the proposed method. As anticipated, the most delicate quantity is the momentum $\rho u$, for which the scheme has more difficulty in reproducing the sharp peak. This behavior can be attributed both to the coarser resolution used in the test case ($N=1000$ compared to $N=4000$ in the reference solution) and, more significantly, to the larger time step dictated by the choice $\text{CFL}=0.9$, whereas the IMEX-RK2 reference solution was computed with a smaller effective step corresponding to $\text{CFL}=0.5$.
\subsubsection{Isentropic Euler system with stiff friction}
For the isentropic model, the pressure is determined by the equation of state
\begin{equation*}
p = k \rho^{\gamma}, \qquad k = 1.
\end{equation*}
\begin{figure}[!ht]
     \centering
     \begin{subfigure}[b]{0.5\textwidth}
         \centering
         \includegraphics[width=\textwidth]{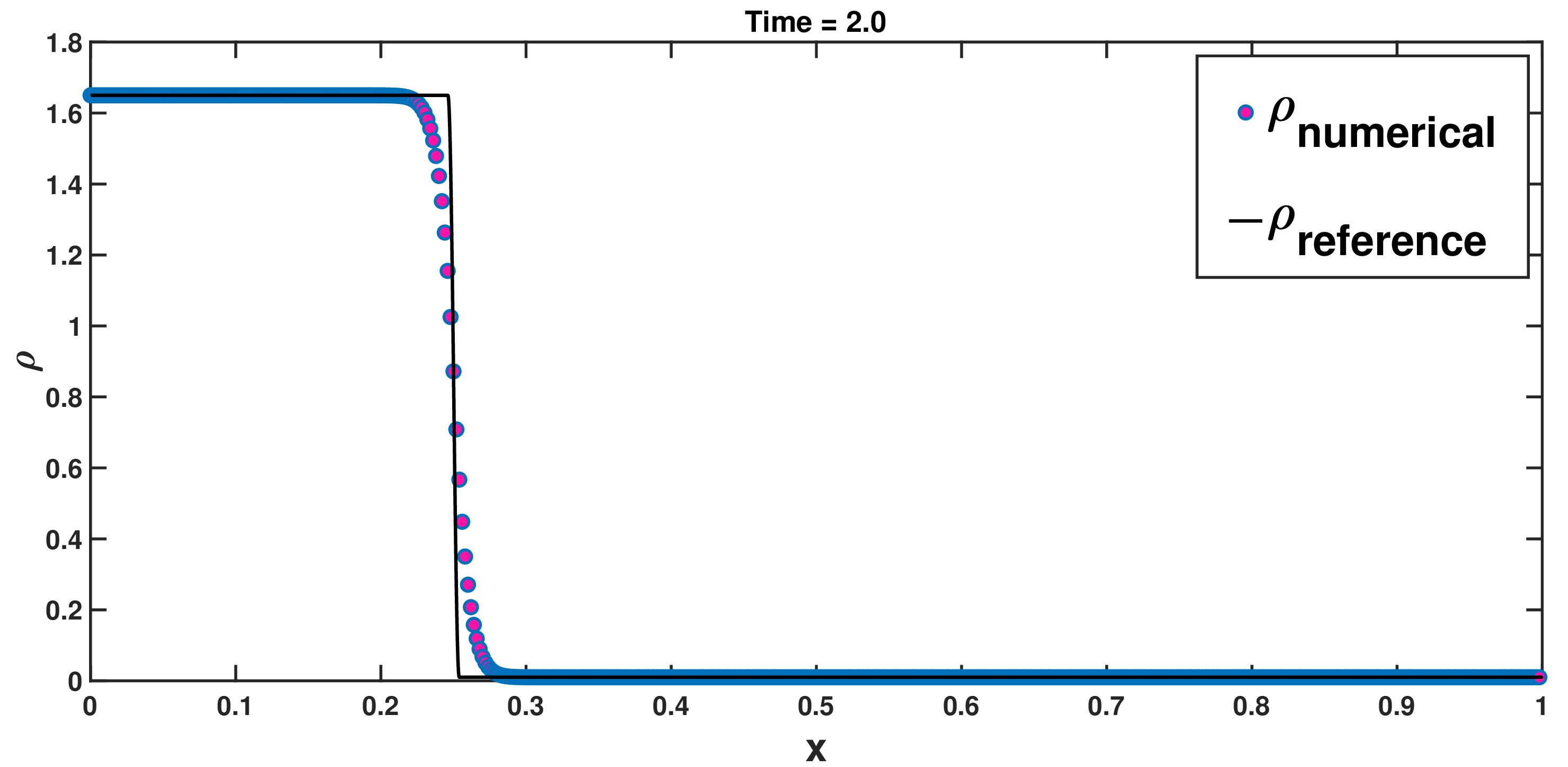}
         \caption{$\rho$ and $\alpha = 10^{8}$}
         \label{euler_isentropic:1}
     \end{subfigure}
     \hfill
     \begin{subfigure}[b]{0.48\textwidth}
         \centering
         \includegraphics[width=\textwidth]{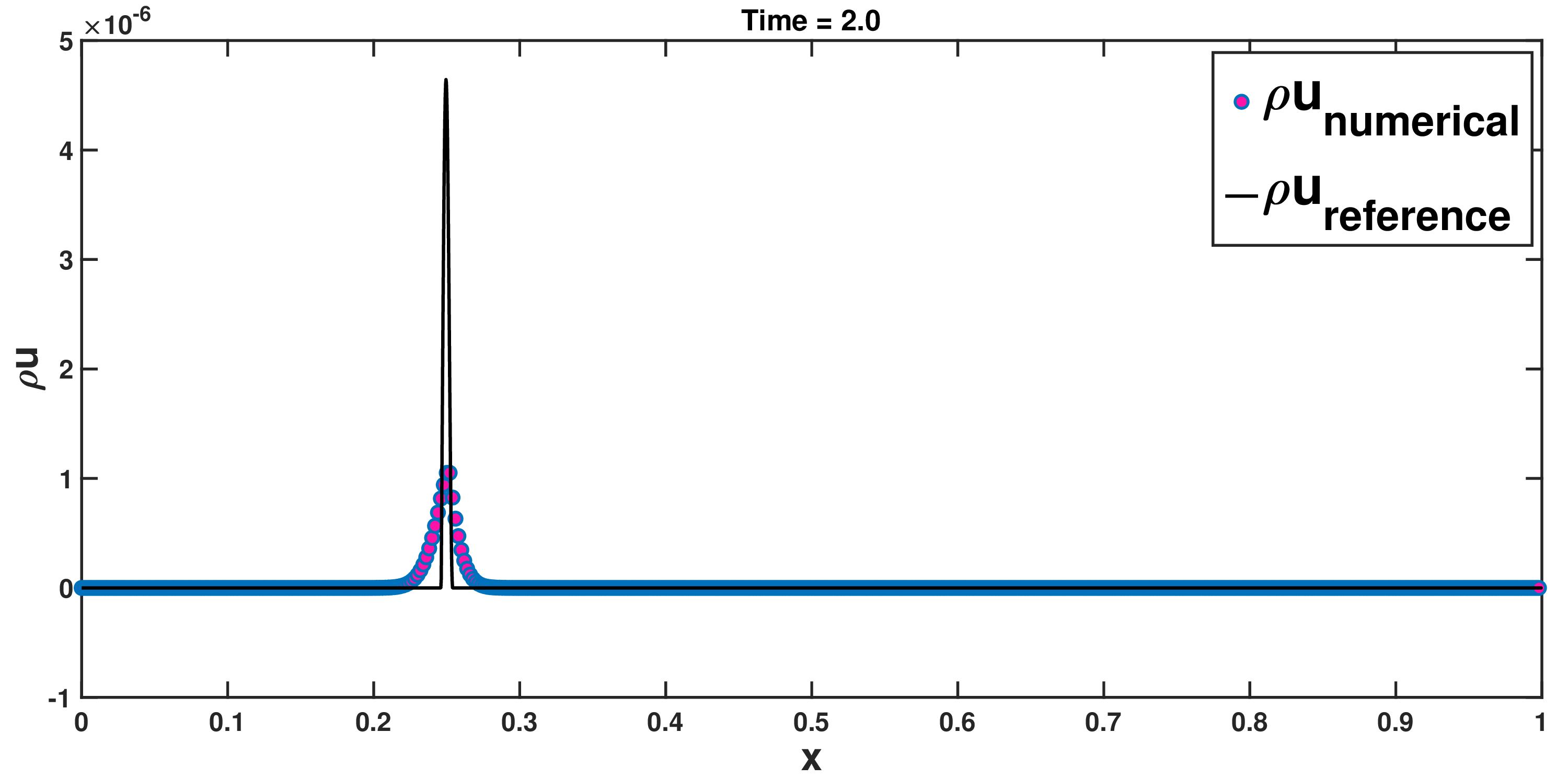}
         \caption{$\rho u$ and $\alpha = 10^{8}$}
         \label{euler_isentropic:2}
     \end{subfigure}
     \caption{Isentropic Euler system with stiff friction model with an non-smooth case: comparison between Numerical solution $\rho$ (left) and $\rho u$ (right) and reference solution with CFL $0.9$ and $N=1000$.}
        \label{euler_isentropic:example1}
\end{figure}
Using the same numerical setup as in the previous test, the results are reported in Figure~\ref{euler_isentropic:example1}, where the density $\rho$ (left) and the momentum $\rho u$ (right) are compared with the corresponding high-resolution reference solutions computed with IMEX-RK2 on a grid of $N=4000$ points. The overall agreement is very good, confirming once again the reliability of the proposed approach. As in the full Euler case, the variable $\rho u$ is more challenging to capture accurately, especially in the peak region. This behavior is consistent with the observations made for the full Euler equations and is mainly due to the coarser discretization ($N=1000$ vs.~$N=4000$ in the reference) and the larger time step associated with $\text{CFL} = 0.9$, compared to the smaller $\text{CFL} = 0.5$ employed in the IMEX-RK2 scheme.
\subsection{Two-Dimensional Hyperbolic Relaxation Model}
\begin{figure}[!ht]
     \centering
     \begin{subfigure}[b]{0.32\textwidth}
         \centering  
         \includegraphics[scale=0.21, trim=12.5cm 0cm 11.5cm 0cm, clip]{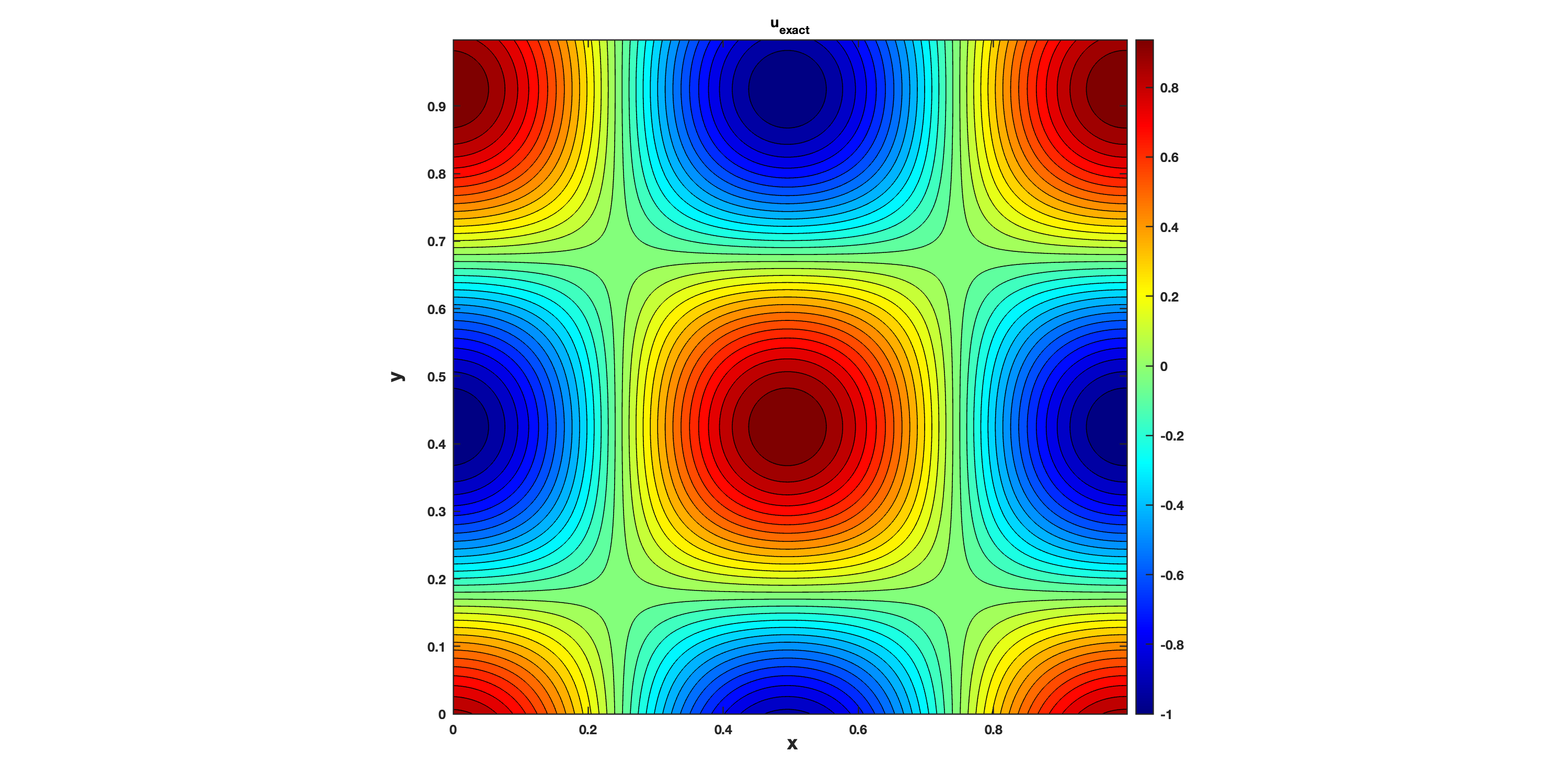}
      \caption{$u$ exact}
         \label{2dxinu:exact}
     \end{subfigure}
     \hfill
     \begin{subfigure}[b]{0.32\textwidth}
         \centering   
         \includegraphics[scale=0.215, trim=12.5cm 0cm 11.5cm 0cm, clip]{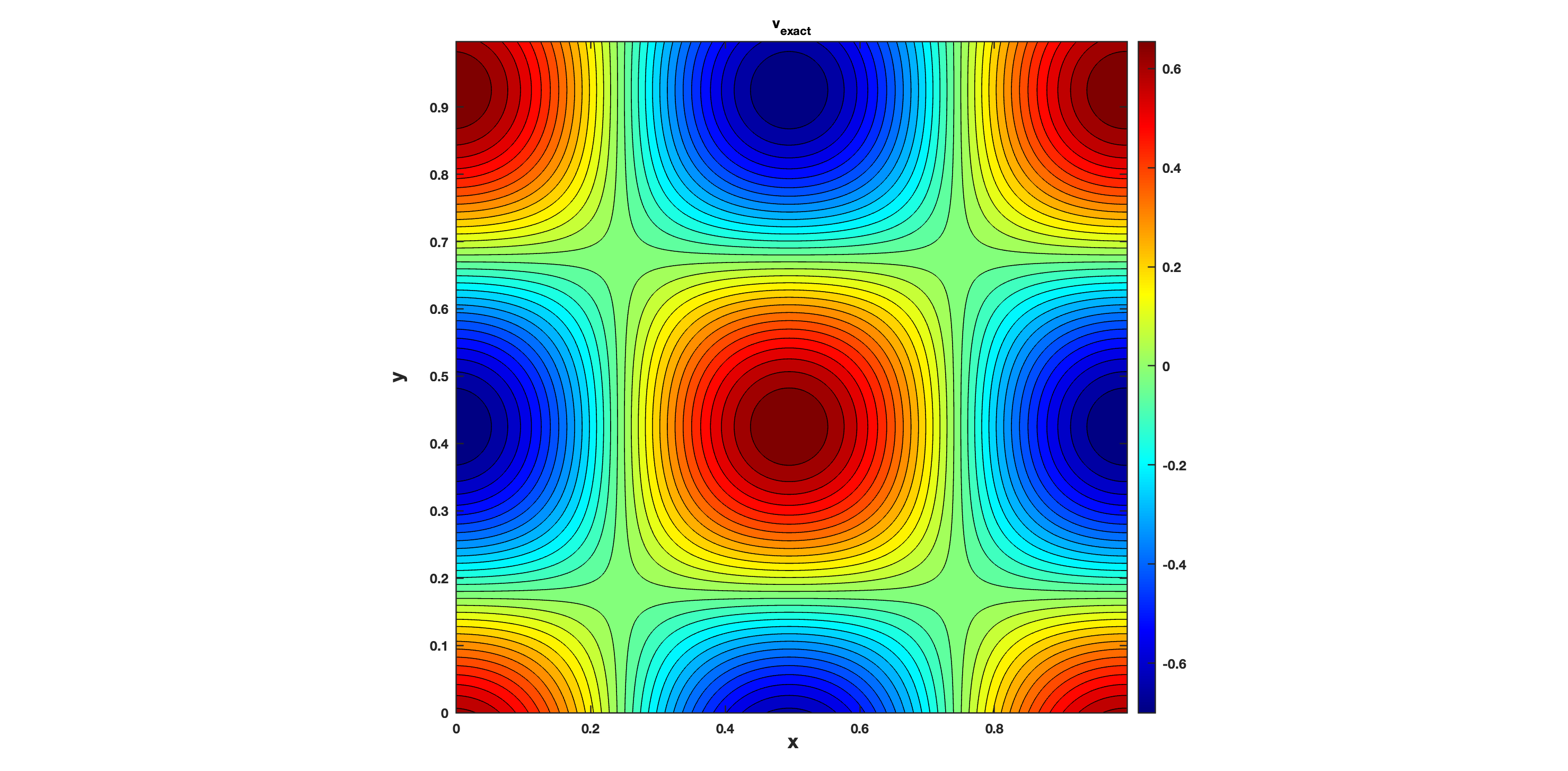}
         \caption{$v$ exact}
         \label{2dxinv:exact}
     \end{subfigure}
     \hfill
    \begin{subfigure}[b]{0.32\textwidth}
         \centering
      \includegraphics[scale=0.215, trim=12.5cm 0cm 11.5cm 0cm, clip]{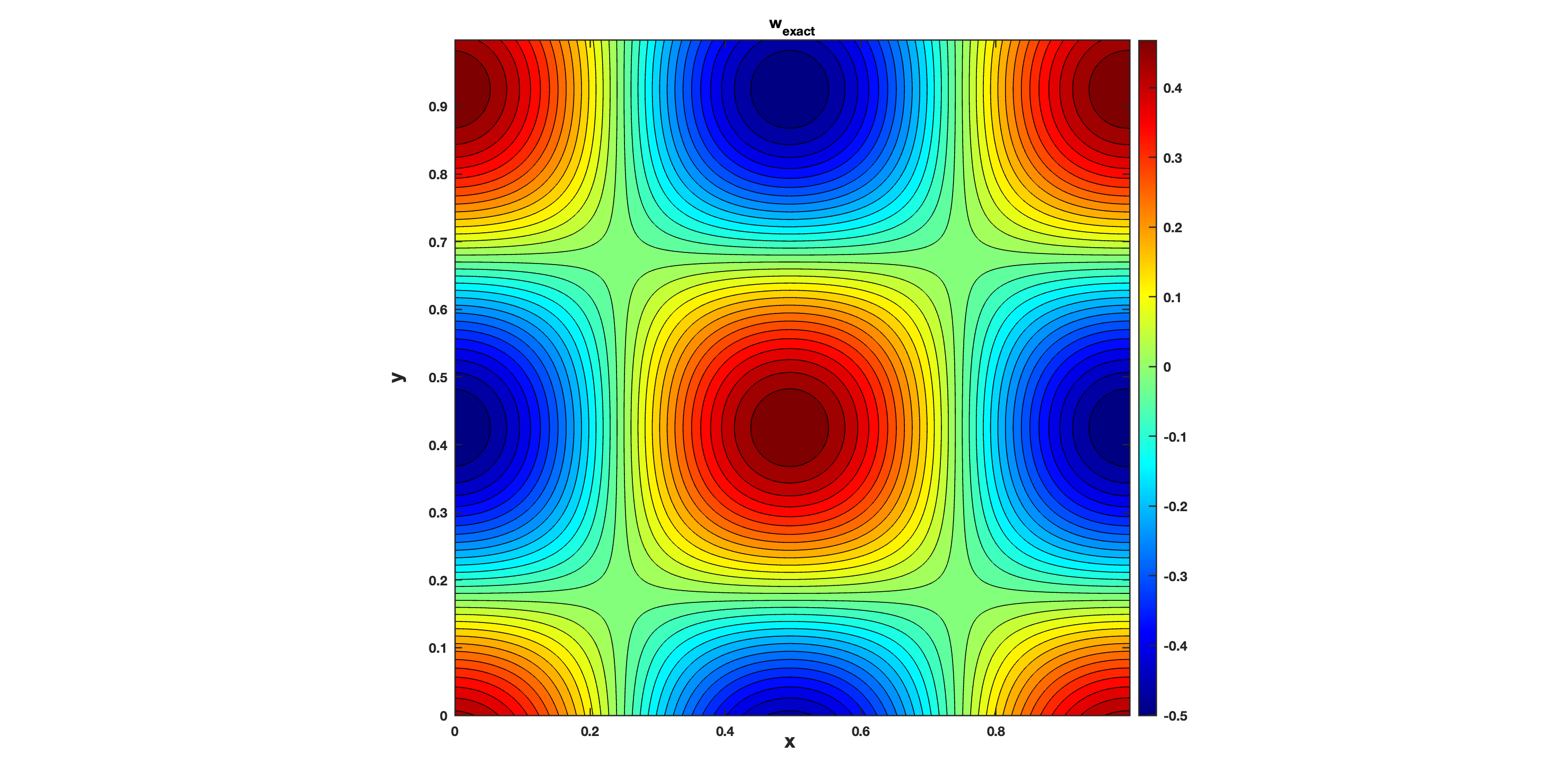}
     \caption{$w$ exact}
         \label{2dxinw:exact}
     \end{subfigure}
     \centering
     \begin{subfigure}[b]{0.32\textwidth}
         \centering  
         \includegraphics[scale=0.215, trim=12.5cm 0cm 11.5cm 0cm, clip]{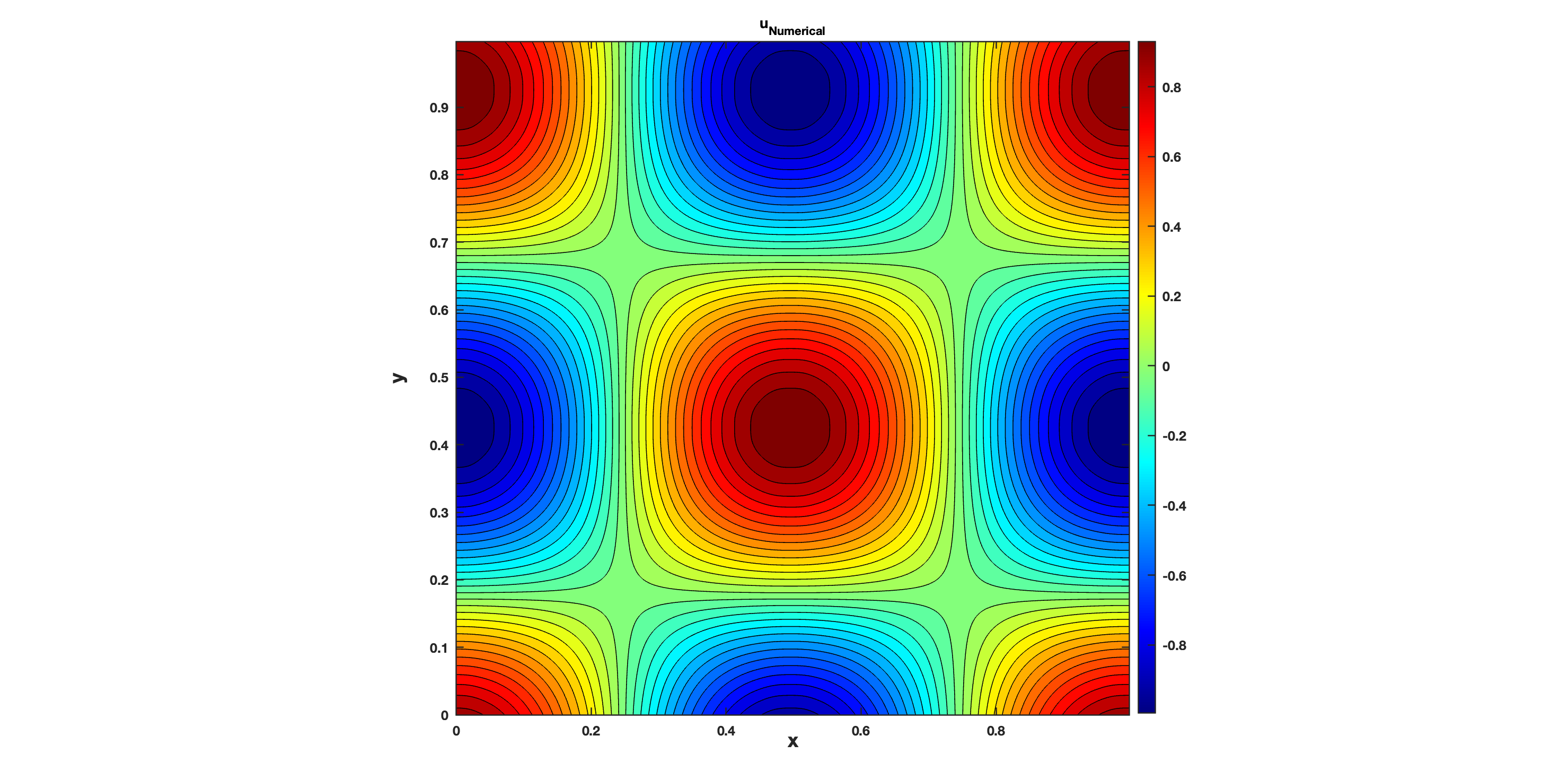}
      \caption{$u$ numerical}
         \label{2dxinu:num}
     \end{subfigure}
     \hfill
     \begin{subfigure}[b]{0.32\textwidth}
         \centering   
         \includegraphics[scale=0.215, trim=12.5cm 0cm 11.5cm 0cm, clip]{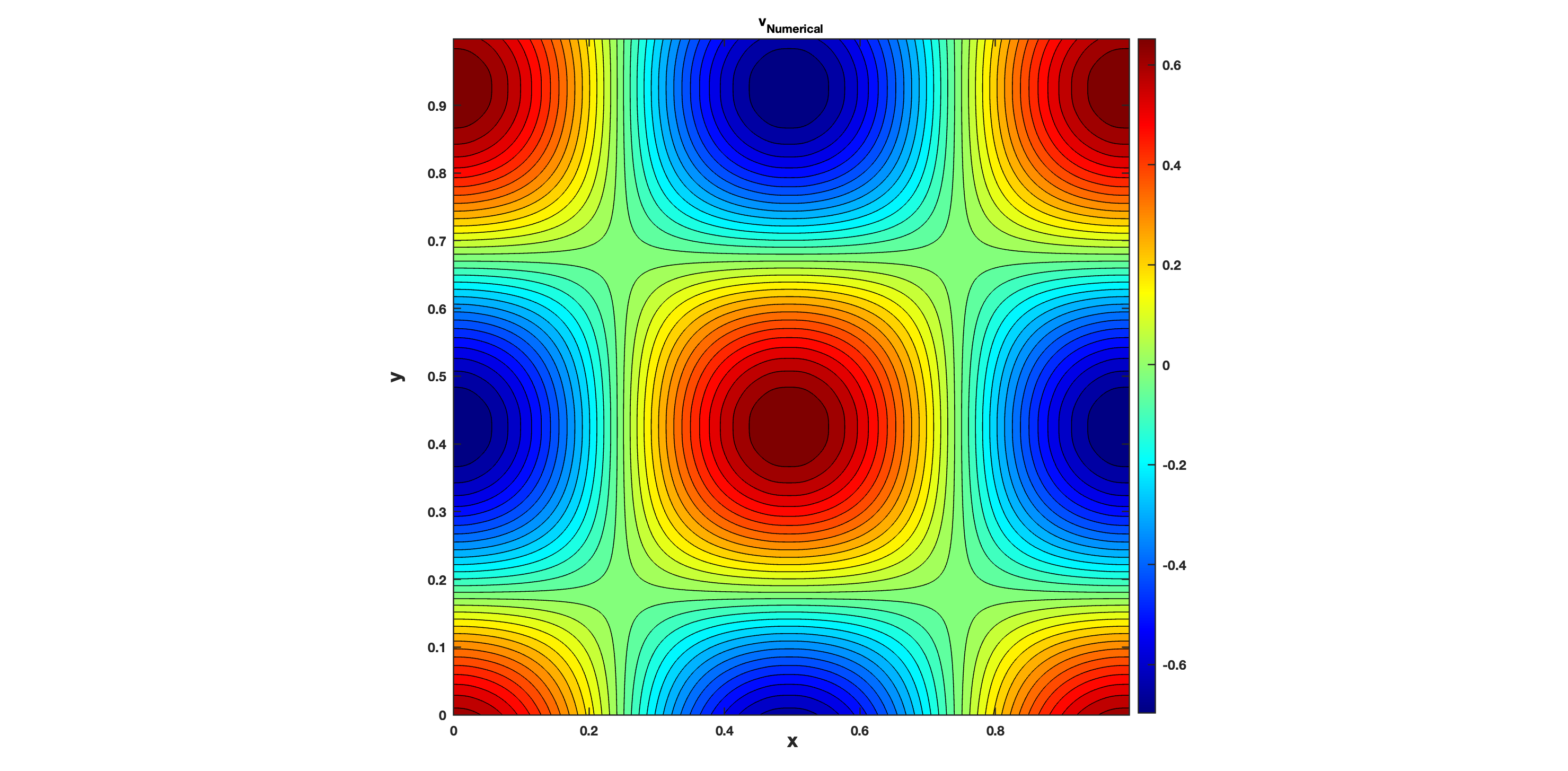}
         \caption{$v$ numerical}
         \label{2dxinv:num}
     \end{subfigure}
     \hfill
    \begin{subfigure}[b]{0.325\textwidth}
         \centering
      \includegraphics[scale=0.215, trim=12.5cm 0cm 11.7cm 0cm, clip]{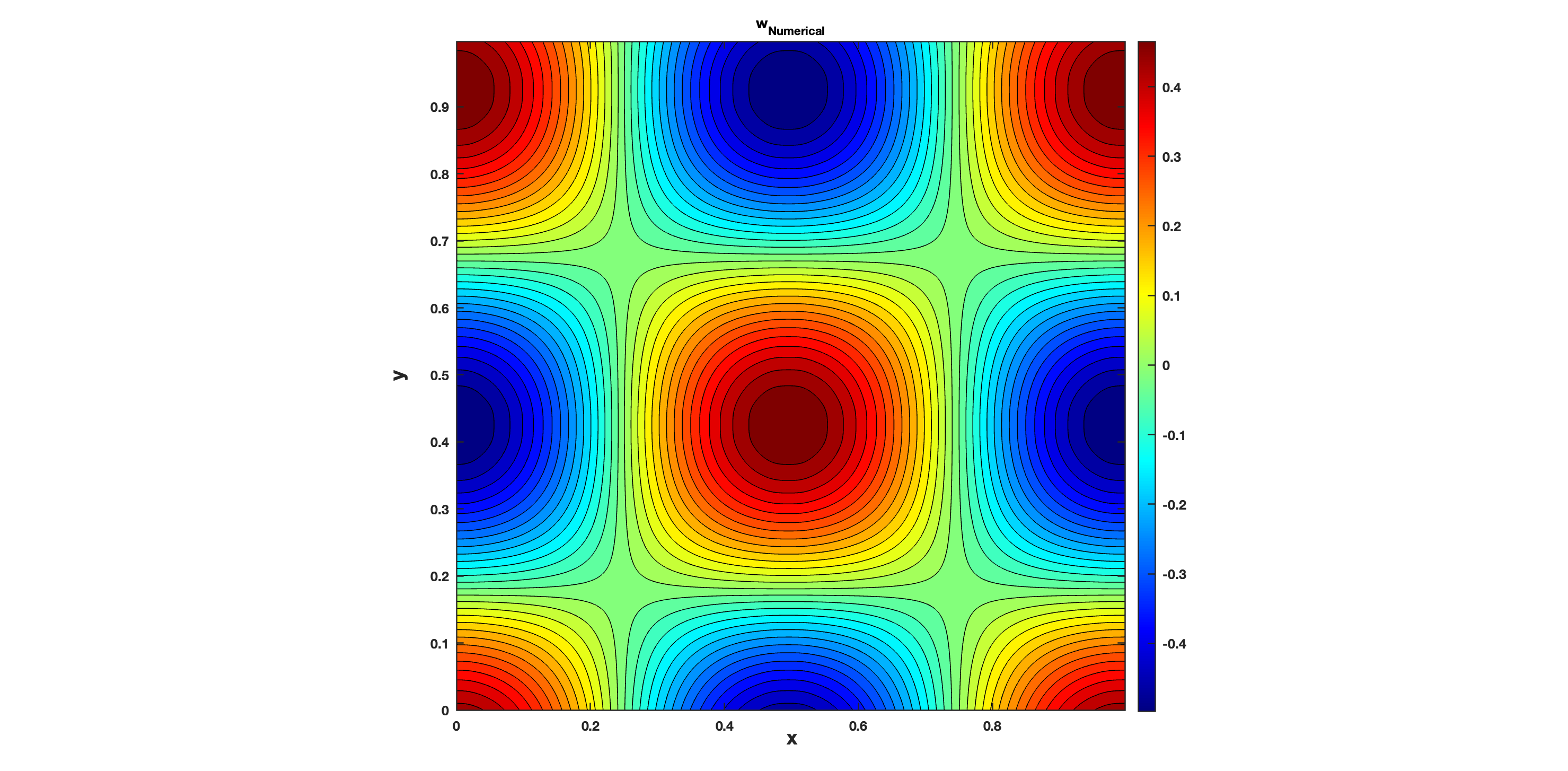}
     \caption{$w$ numerical}
         \label{2dxinw:num}
     \end{subfigure}
         \caption{Exact and numerical contour plots with $30$ levels of the variables $u$, $v$, and $w$ for the two-dimensional smooth test problem computed using the Jin-Xin relaxation system. The computations are carried out on a uniform $400 \times 400$ Cartesian grid with CFL number $0.9$, relaxation parameter $\varepsilon = 10^{-10}$, and final time $T = 0.35$.}
    \label{2dxin:num}
\end{figure}
\subsubsection{2d Jin-Xin Model}
We begin our analysis with the two-dimensional Jin-Xin relaxation system\begin{equation}\label{Xin_Jin2d_model}
\begin{cases}
u_t + v_x +w_y = 0,\\[2mm]
v_t + u_x = -\dfrac{1}{\varepsilon}\left(v - au\right),\\[2mm]
w_t + u_y = -\dfrac{1}{\varepsilon}\left(w - bu\right),\\
\end{cases}
\end{equation}
where $\varepsilon > 0$ denotes the relaxation parameter, and $a$ and $b$ are given constants. The system is considered on the spatial domain $[0,1]\times[0,1]$ with periodic boundary conditions.\\
\noindent\textbf{Smooth case (Well-Prepared):}
The smooth initial data is as follows
\begin{equation}
\label{smoothdata:XinJin2d}
\begin{cases}
u(x,y,0) = \sin(2\pi x)\sin(2\pi y), \\[2mm]
v(x,y,0) = a\,u(x,y,0), \\[2mm]
w(x,y,0) = b\,u(x,y,0),
\end{cases}
\end{equation}
where the parameters are fixed as $a = 0.7$ and $b = 0.5$. With this choice, the initial data lies on the equilibrium manifold, leading to a smooth evolution of the solution. This configuration provides an appropriate setting to assess the accuracy and stability of the proposed numerical scheme in the absence of discontinuities. The spatial domain $[0,1]\times[0,1]$ is discretized using a uniform Cartesian mesh with $N_x = N_y = 400$ grid points in the $x$- and $y$-directions. The numerical solution is computed up to the final time $T = 0.35$ with a CFL number of $0.9$. Figure~\ref{2dxin:num} compares the exact and numerical solutions of the variables $u$, $v$, and $w$ for $\varepsilon = 10^{-10}$. The top row presents the exact solutions, whereas the bottom row shows the numerical results. In each row, the left, center, and right panels correspond to $u$, $v$, and $w$, respectively. All contour plots use 30 levels. The exact solutions are obtained analytically using the Fourier transform method. Tables~\ref{Tabu}, \ref{Tabv}, and \ref{Tabw} report the experimental order of convergence for 
$u$, $v$, and $w$, respectively.
\begin{table}[ht!]
\caption{$L^{1}$-Error and order of convergence for $u$ (2D Jin-Xin model)}
\centering
\begin{tabular}{*{7}{c}}
\toprule
\multirow{2}{*}{N} & 
\multicolumn{2}{c}{$\varepsilon=10^{-10}$} &
\multicolumn{2}{c}{$\varepsilon=10^{-8}$} & 
\multicolumn{2}{c}{$\varepsilon=10^{-7}$} \\
\cmidrule(lr){2-3}
\cmidrule(lr){4-5}
\cmidrule(lr){6-7}
 & $L^{1}$-Error & Order 
 & $L^{1}$-Error & Order 
 & $L^{1}$-Error & Order \\
\midrule
32  & 1.4418e-03 & - 
    & 1.4418e-03 & - 
    & 1.4418e-03 & - \\
64  & 3.5003e-04 & 2.04
    & 3.5003e-04 & 2.04
    & 3.5003e-04 & 2.04 \\
128 & 8.2978e-05 & 2.08
    & 8.2978e-05 & 2.08
    & 8.2977e-05 & 2.08 \\
256 & 1.9824e-05 & 2.07
    & 1.9824e-05 & 2.07
    & 1.9824e-05 & 2.07 \\
\bottomrule
\end{tabular}
\label{Tabu}
\end{table}

\begin{table}[ht!]
\caption{$L^{1}$-Error and order of convergence for $v$ (2D Jin-Xin model)}
\centering
\begin{tabular}{*{7}{c}}
\toprule
\multirow{2}{*}{N} & 
\multicolumn{2}{c}{$\varepsilon=10^{-10}$} &
\multicolumn{2}{c}{$\varepsilon=10^{-8}$} & 
\multicolumn{2}{c}{$\varepsilon=10^{-7}$} \\
\cmidrule(lr){2-3}
\cmidrule(lr){4-5}
\cmidrule(lr){6-7}
 & $L^{1}$-Error & Order 
 & $L^{1}$-Error & Order 
 & $L^{1}$-Error & Order \\
\midrule
32  & 1.0093e-03 & - 
    & 1.0093e-03 & - 
    & 1.0093e-03 & - \\
64  & 2.4502e-04 & 2.04
    & 2.4502e-04 & 2.05
    & 2.4502e-04 & 2.05 \\
128 & 5.8085e-05 & 2.08
    & 5.8084e-05 & 2.08
    & 5.8085e-05 & 2.08 \\
256 & 1.3877e-05 & 2.07
    & 1.3878e-05 & 2.07
    & 1.3878e-05 & 2.07 \\
\bottomrule
\end{tabular}
\label{Tabv}
\end{table}

\begin{table}[ht!]
\caption{$L^{1}$-Error and order of convergence for $w$ (2D Jin-Xin model)}
\centering
\begin{tabular}{*{7}{c}}
\toprule
\multirow{2}{*}{N} & 
\multicolumn{2}{c}{$\varepsilon=10^{-10}$} &
\multicolumn{2}{c}{$\varepsilon=10^{-8}$} & 
\multicolumn{2}{c}{$\varepsilon=10^{-7}$} \\
\cmidrule(lr){2-3}
\cmidrule(lr){4-5}
\cmidrule(lr){6-7}
 & $L^{1}$-Error & Order 
 & $L^{1}$-Error & Order 
 & $L^{1}$-Error & Order \\
\midrule
32  & 7.2092e-04 & - 
    & 7.2092e-04 & - 
    & 7.2091e-04 & - \\
64  & 1.7502e-04 & 2.04
    & 1.7502e-04 & 2.04
    & 1.7502e-04 & 2.04 \\
128 & 4.1489e-05 & 2.08
    & 4.1489e-05 & 2.08
    & 4.1489e-05 & 2.08 \\
256 & 9.9122e-06 & 2.07
    & 9.9122e-06 & 2.07
    & 9.9122e-06 & 2.07 \\
\bottomrule
\end{tabular}
 \label{Tabw}
\end{table}
\subsubsection{2D Euler Jin-Xin-type relaxation Model}

\begin{figure}[htbp]
    \begin{subfigure}{0.48\linewidth}
        \includegraphics[scale=0.335, trim=12.3cm 0cm 12.3cm 0cm, clip]{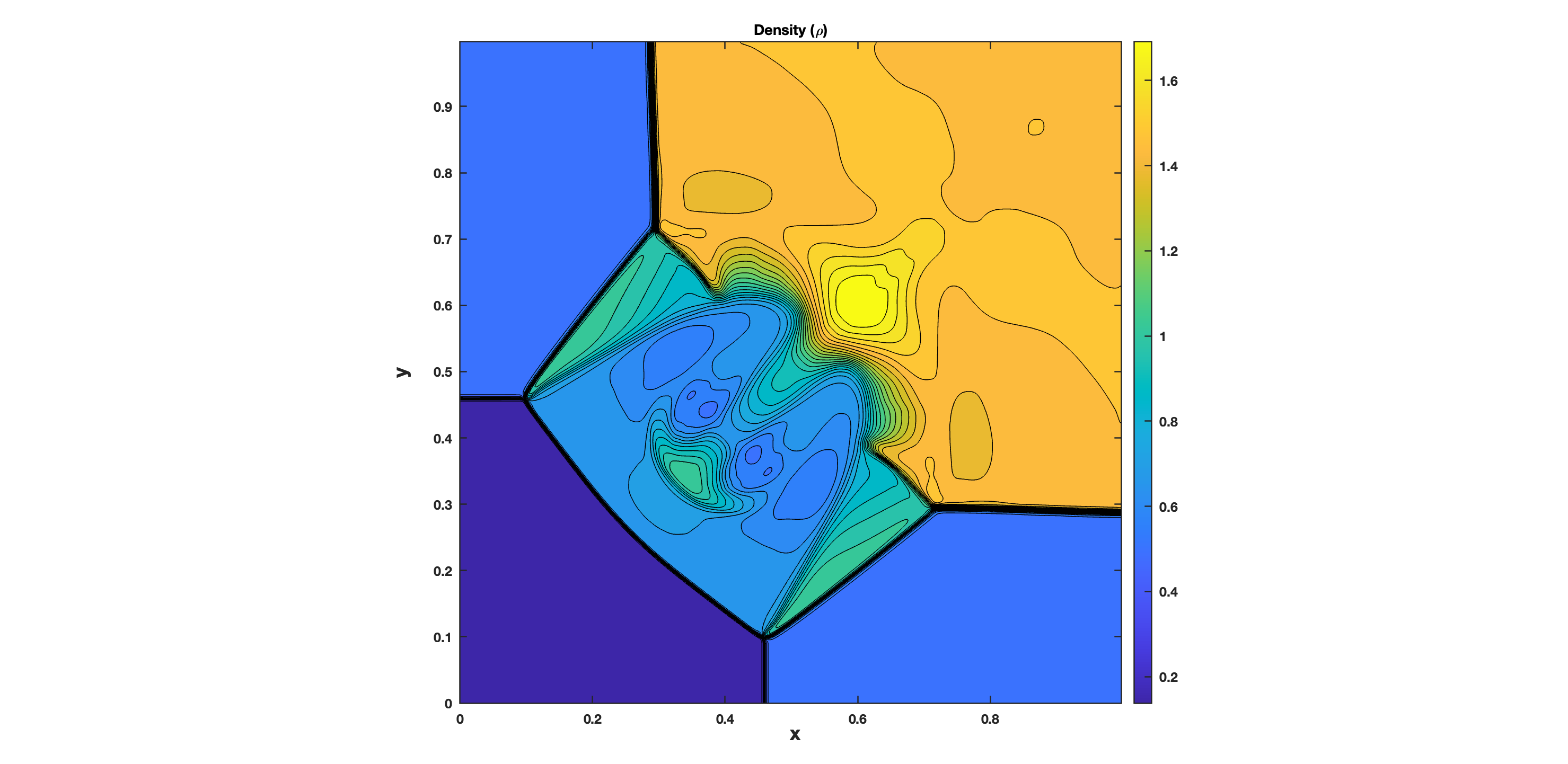}
        \caption{Configuration 3}
        \label{2drp1}
    \end{subfigure}
    \hfill
    \begin{subfigure}[b]{0.48\linewidth}
        \includegraphics[scale=0.324, trim=12.3cm 0cm 12.3cm 0cm, clip]{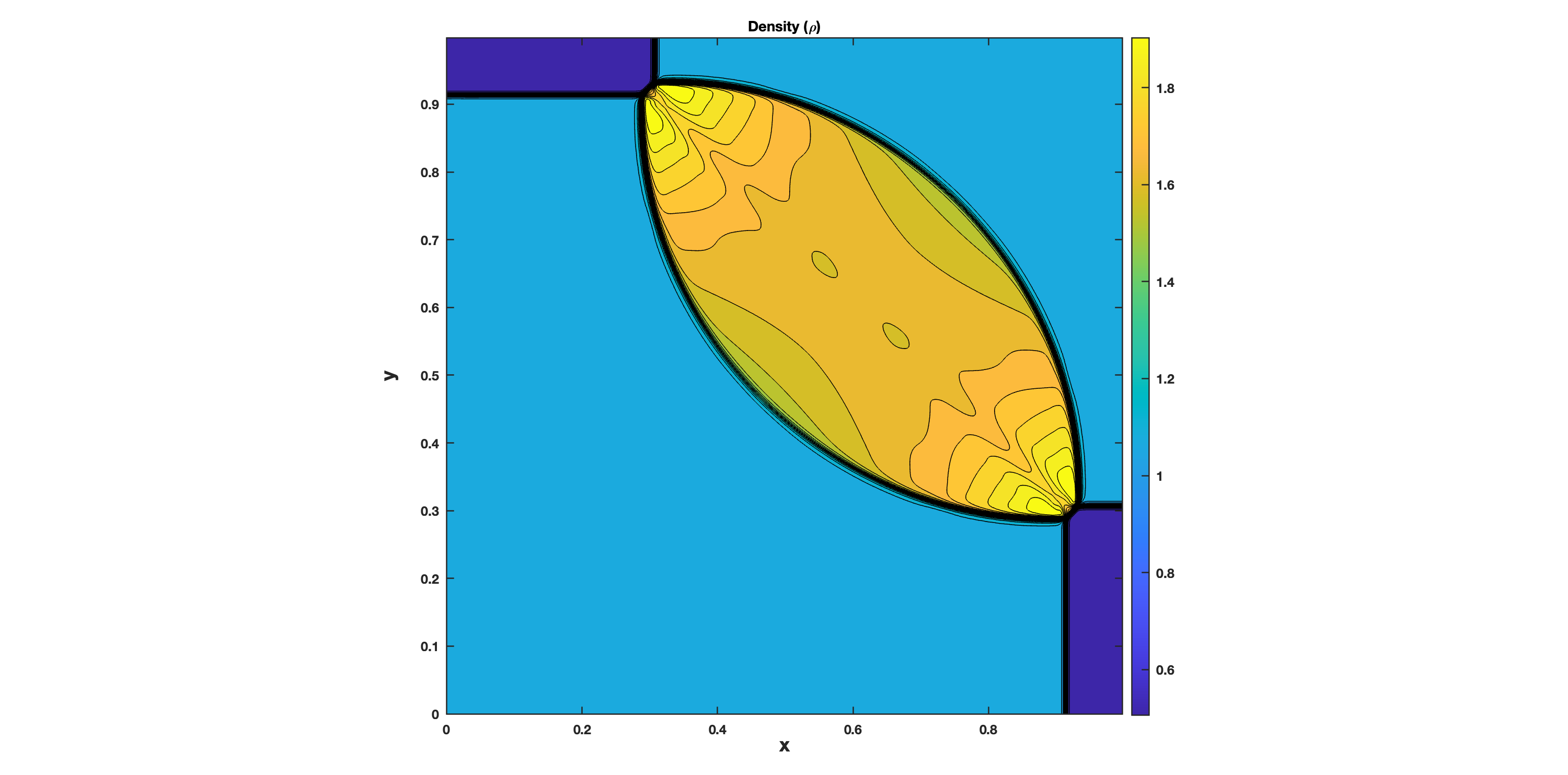}
        \caption{Configuration 4}
        \label{2drp2}
    \end{subfigure}

    \vspace{0.5cm} 

    \centering
    \begin{subfigure}[b]{0.48\linewidth}
        \includegraphics[scale=0.335, trim=12.3cm 0cm 12.3cm 0cm, clip]{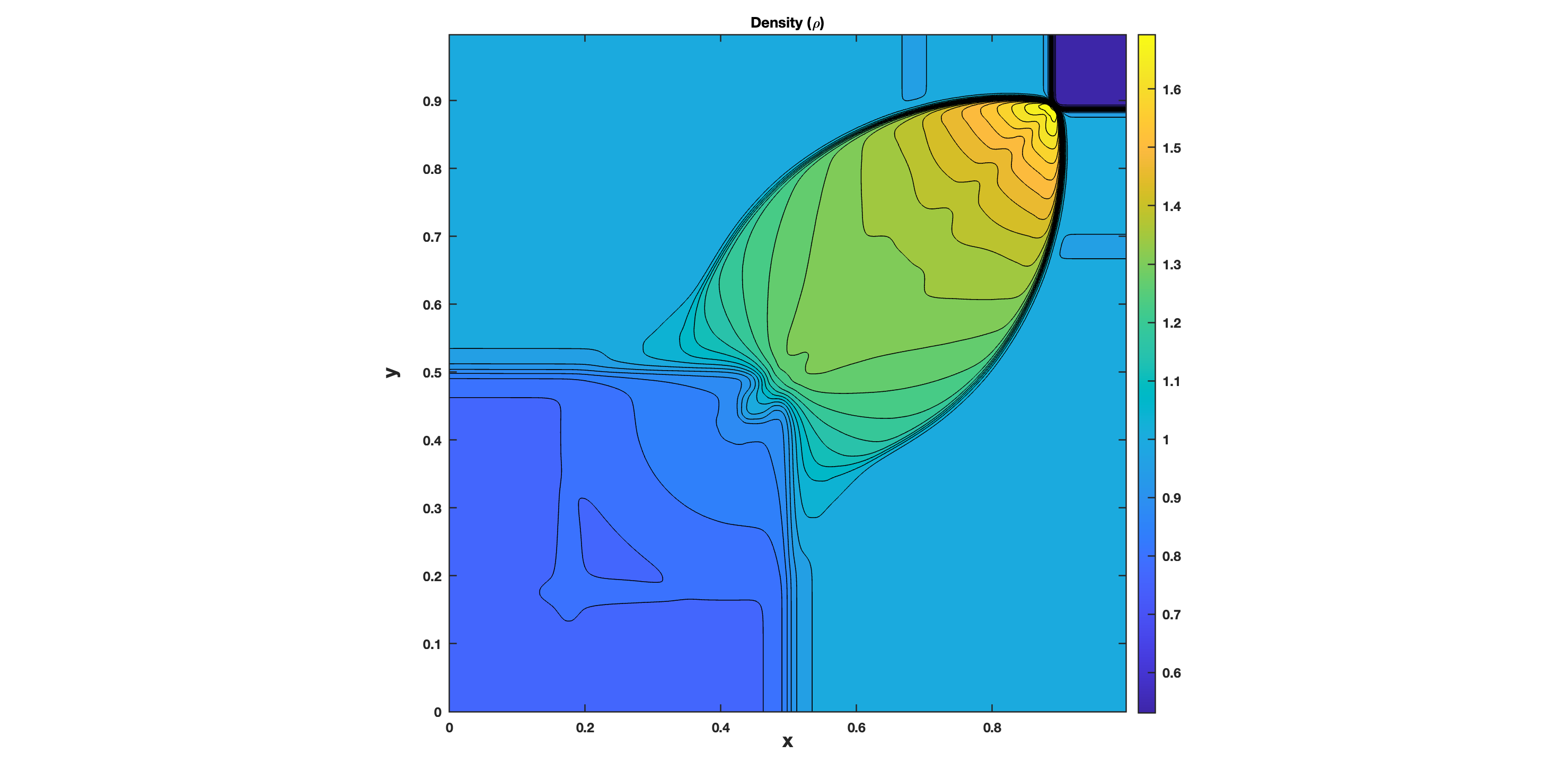}
        \caption{Configuration 12}
        \label{2drp3}
    \end{subfigure}
    
    \caption{Density ($\rho$) contour plots with $30$ levels for three two-dimensional Riemann problems computed using the Euler Jin-Xin-type relaxation system. The simulations are performed on a uniform $400 \times 400$ Cartesian grid with CFL number $0.9$ and relaxation parameter $\varepsilon = 10^{-12}$. The corresponding final times are $T = 0.8$, $T = 0.25$, and $T = 0.25$, respectively.}
    \label{2drp}
\end{figure}

The two-dimensional Euler Jin-Xin type relaxation system can be written as
\begin{equation}
\begin{cases}
U_t + V_x + W_y = 0,\\
V_t + A U_x = -\dfrac{1}{\varepsilon}\left(V - F_1(U)\right),\\
W_t + B U_y = -\dfrac{1}{\varepsilon}\left(W - F_2(U)\right),
\end{cases}
\end{equation}
where $\varepsilon > 0$ is the relaxation parameter. The matrices $A$ and $B$ are diagonal, with eigenvalues $a_i$ and $b_i$, respectively. In this work, we choose
\[
a_1 = a_2 = a_3 = a_4 = 9.75,
\qquad
b_1 = b_2 = b_3 = b_4 = 7.62.
\]
We consider three types of two-dimensional Riemann problems \cite{KurganovTadmor2002} as initial conditions.
\textbf{Configuration 3:}
The initial data for $(\rho, u, v, p)$ is prescribed as
\[
(\rho, u, v, p)(x,y,0) =
\begin{cases}
(1.5, 0, 0, 1.5), & x \ge 0.8,\; y \ge 0.8,\\
(0.5323, 1.2060, 0, 0.3), & x < 0.8,\; y \ge 0.8,\\
(0.1380, 1.2060, 1.2060, 0.0290), & x < 0.8,\; y < 0.8,\\
(0.5323, 0, 1.2060, 0.3), & x \ge 0.8,\; y < 0.8.
\end{cases}
\]
with
\[
V(x,y,0) = F_1(U(x,y,0)), 
\qquad 
W(x,y,0) = F_2(U(x,y,0)).
\]
\textbf{Configuration 4:}
\[
(\rho, u, v, p)(x,y,0) =
\begin{cases}
(1.1, 0, 0, 1.1), & x \ge 0.5,\; y \ge 0.5,\\
(0.5026, 0.8939, 0, 0.35), & x < 0.5,\; y \ge 0.5,\\
(1.1, 0.8939, 0.8939, 1.1), & x < 0.5,\; y < 0.5,\\
(0.5065, 0, 0.8939, 0.35), & x \ge 0.5,\; y < 0.5.
\end{cases}
\]
with
\[
V(x,y,0) = F_1(U(x,y,0)), 
\qquad 
W(x,y,0) = F_2(U(x,y,0)).
\]
\textbf{Configuration 12:}
\[
(\rho, u, v, p)(x,y,0) =
\begin{cases}
(0.5313, 0, 0, 0.4), & x \ge 0.5,\; y \ge 0.5,\\
(1, 0.7276, 0, 1), & x < 0.5,\; y \ge 0.5,\\
(0.8, 0, 0, 1), & x < 0.5,\; y < 0.5,\\
(1, 0, 0.7276, 1), & x \ge 0.5,\; y < 0.5.
\end{cases}
\]
with
\[
V(x,y,0) = F_1(U(x,y,0)), 
\qquad 
W(x,y,0) = F_2(U(x,y,0)).
\]
Figure~\ref{2drp} presents the density contours (30 contour levels) for the three two-dimensional Riemann problems computed using the Euler Jin-Xin relaxation system. The computations are performed on a uniform $400 \times 400$ grid over the domain $[0,1] \times [0,1]$. Transmissive boundary conditions are imposed in both spatial directions. 
The CFL number is set to $0.9$, and the relaxation parameter is taken as $\varepsilon = 10^{-12}$. The final simulation times are $T = 0.8$ for Configuration~3 and $T  = 0.25$ for Configuration~4 and~12, in which we observed the symmetric resolution of the numerical solution with the proposed approach.

\section{Conclusion} 
In this study, we have addressed the numerical approximation of hyperbolic balance laws endowed with stiff relaxation source terms, with numerical investigations conducted in both one and two spatial dimensions. By generalizing the classical Nessyahu-Tadmor (NT) central scheme to incorporate the nonhomogeneous component through an implicit trapezoidal update and a backward semi–implicit Taylor expansion, we introduced the CS-EBT2 method. A detailed consistency analysis confirmed that the proposed scheme retains second order accuracy as the relaxation parameter $\varepsilon \to 0$, while the stability study revealed a significantly enlarged stability region, allowing larger CFL numbers compared to conventional central schemes.

Our numerical investigation encompassed five canonical model problems: the linear Jin-Xin relaxation model, the shallow water model, the nonlinear Broadwell equations, the Euler system with a stiff heat transfer term, and the Euler system with stiff friction, which includes both the full and the isentropic versions. Moreover, multidimensional numerical experiments were conducted for two-dimensional Jin-Xin and Jin-Xin-type Euler relaxation systems, confirming the effectiveness of the proposed scheme in higher-dimensional regimes. Through the Jin-Xin model with smooth initial data, we verified the theoretical order of accuracy and presented corresponding error and convergence tables for both the Jin-Xin and Broadwell models. The scheme was further tested for various initial conditions smooth as well as discontinuous across all models with relatively large CFL numbers. Comparisons with reference IMEX-RK2 solutions on finer grids demonstrated that the CS-EBT2 method accurately captures sharp discontinuities and intricate relaxation dynamics without producing spurious oscillations, even as $\varepsilon$ approaches zero. Its semi–implicit treatment effectively enhances both stability and performance, ensuring reliable results under strongly stiff regimes.

Beyond its immediate performance gains, CS–EBT2 offers a flexible framework in which the semi-implicit treatment can be extended to more complex geometries, higher dimensions, or coupled multi-physics systems. The successful two-dimensional simulations further confirm the robustness and multidimensional applicability of the proposed approach. Future investigations will focus on adaptive mesh refinement strategies to further optimize computational effort in localized regions of stiffness, as well as the incorporation of positivity-preserving limiters to guarantee adherence to physical constraints in multi-component flows. Overall, the CS-EBT2 scheme is a reliable method for simulating stiff hyperbolic phenomena, providing stability, accuracy, and straightforward implementation. Extending the present approach to higher-order variants (order greater than two) will be investigated in future work.
\section*{Acknowledgments}
This research has received funding from the European Union’s NextGenerationUE – Project: Centro Nazionale HPC, Big Data e Quantum Computing, “Spoke 1” (No. CUP E63C22001000006). E. Macca was partially supported by GNCS No. CUP E53C24001950001 Research Project "Soluzioni Innovative per Sistemi Complessi: Metodi Numerici e Approcci Multiscala"; PRIN 2022 PNRR “FIN4GEO: Forward and Inverse Numerical Modeling of hydrothermal systems in volcanic regions with application to geothermal energy exploitation”, No. P2022BNB97; PRIN 2022 “Efficient numerical schemes and optimal control methods for time-dependent partial differential equations”, No. 2022N9BM3N - Finanziato dall’Unione europea - Next Generation EU – CUP: E53D23005830006; PRIN Project 2022  (2022KA3JBA, entitled “Advanced numerical methods for time dependent parametric partial differential equations and applications”). E. Macca is a member of the INdAM Research group GNCS. The authors Sudipta Sahu and Rathan Samala is supported by NBHM, DAE, India (Ref. No. 02011/46/2021 NBHM(R.P.)/R and D II/14874).

\bibliographystyle{ieeetr}
\bibliography{biblio}

\end{document}